\documentclass[9pt]{article}

\usepackage{amssymb, amsmath, latexsym, amscd, amsthm}
\usepackage{pdfsync}
\usepackage{graphicx, color, epsfig}  
\usepackage{epstopdf}
\usepackage[all,cmtip]{xy}

\def\sref#1{Section~\ref{#1}}

\def\tref#1{Theorem~\ref{#1}}
\def\dref#1{Definition~\ref{#1}}

\def\rref#1{Remark~\ref{#1}}
\def\cref#1{Corollary~\ref{#1}}
\def\pref#1{Proposition~\ref{#1}}
\def\lref#1{Lemma~\ref{#1}}

\def\fgref#1{Figure~\ref{#1}}
\def\aref#1{Appendix~\ref{#1}}

\def\nref#1{Notation~\ref{#1}}

\newtheoremstyle{thm}{1em}{1em}{\itshape}{}{\bfseries}{.}{.5em}{} 
\newtheoremstyle{rem}{1em}{1em}{}{}{\bfseries}{.}{.5em}{}           

\theoremstyle{thm}                
\newtheorem{thm}{Theorem}[section]
\newtheorem{prop}[thm]{Proposition}

\newtheorem{cor}[thm]{Corollary}
\newtheorem{lem}[thm]{Lemma}

\newtheorem{dfn}[thm]{Definition}

\theoremstyle{rem} 
\newtheorem{rmk}[thm]{Remark}
\newtheorem{rmks}[thm]{Remarks}

\newtheorem{nota}[thm]{Notation}
\newtheorem{notas}[thm]{Notations}

\newcommand{\intro}[1]
{\renewcommand{\thesection}{\fnsymbol{section}}
\setcounter{section}{-1}
\section{#1}
\renewcommand{\thesection}{\arabic{section}}
}

\DeclareMathOperator{\Int}{Int}
\DeclareMathOperator{\Ker}{Ker}
\DeclareMathOperator{\im}{Im}

\DeclareMathOperator{\id}{id}
\DeclareMathOperator{\End}{End}
\DeclareMathOperator{\Def}{def}
\DeclareMathOperator{\In}{Int}

\DeclareMathOperator{\Gr}{Gr}
\DeclareMathOperator{\Op}{Op}

\newcommand{\Pf}{\noindent{\bf Proof}. }
\newcommand{\cqfd}
{%
\mbox{}%
\nolinebreak%
\hfill%
\rule{2mm}{2mm}%
\medbreak%
\par%
}
 
\renewcommand{\b}{{\mathcal B}{}}
\newcommand{\fb}{{\mathfrak b}{}}

\renewcommand{\d}{{\mathcal D}{}}  
\newcommand{\e}{{\mathcal E}{}}
\newcommand{\ee}{{\bf e}{}}
\newcommand{\F}{{\mathcal F}{}}

\newcommand{\h}{{\mathcal H}{}} 
\newcommand{\I}{{\mathcal I}{}} 
\newcommand{\ii}{{\bf i}{}}

\newcommand{\EL}{{\mathcal L}{}} 
\newcommand{\fL}{{\mathfrak L}{}} 

\newcommand{\EN}{{\mathcal N}{}}

\newcommand{\OO}{{\mathcal O}{}} 
\newcommand{\p}{{\mathcal P}{}} 
\newcommand{\PP}{{\rm Pair}{}}

\newcommand{\R}{\mathbb R}

\newcommand{\fs}{{\mathfrak s}{}}
\newcommand{\T}{{\mathcal T}{}}

\newcommand{\U}{{\mathcal U}{}}

\newcommand{\V}{{\mathcal V}{}}
\newcommand{\W}{{\mathcal W}{}}
\newcommand{\X}{{\mathcal X}{}}
\newcommand{\IX}{{\mathfrak X}{}}
\newcommand{\IY}{{\mathfrak Y}{}}
\newcommand{\Y}{{\mathcal Y}{}}
\newcommand{\ZE}{{\mathcal Z}{}}

\newcommand{\eps}{\varepsilon }
\newcommand{\ol}{\overline}
\newcommand{\ul}{\underline}

\newcommand{\rp}{respectively }

\begin{document} 

\title{Affine connections and symmetry jets}

\author{M\'elanie Bertelson$^\ast$, Pierre Bieliavsky$^\dagger$} 

\date{\today}

\maketitle

\begin{center}
$^\ast$ Chercheur Qualifi\'e F.N.R.S. \\
D\'epartement de Math\'ematique, C.P. 218 \\
Universit\'e Libre de Bruxelles \\
 Boulevard du Triomphe \\
1050 Bruxelles -- Belgique\\
{\tt mbertels@ulb.ac.be} \\
\ \\
$^\dagger$ Institut de Recherche en Math\'ematique et Physique \\
Universit\'e Catholique de Louvain \\
2, Chemin du Cyclotron \\
1348 Louvain-la-Neuve -- Belgique \\
{\tt Pierre.Bieliavsky@uclouvain.be}
\end{center}

\begin{abstract} We establish a bijective correspondence between affine connections and a class of semi-holonomic jets of local diffeomorphisms of the underlying manifold called symmetry jets in the text. The symmetry jet corresponding to a torsion free connection consists in the family of $2$-jets of the geodesic symmetries. Conversely, any connection is described in terms of the geodesic symmetries by a simple formula involving only the Lie bracket of vector fields. We then formulate, in terms of the symmetry jet, several aspects of the theory of affine connections and obtain geometric and intrinsic descriptions of various related objects involving the gauge groupoid of the frame bundle. In particular, the property of uniqueness of affine extension admits an equivalent formulation as the property of existence and uniqueness of a certain groupoid morphism. Moreover, affine extension may be carried out at all orders and this allows for a description of the tensors associated to an affine connections, namely the torsion, the curvature and their covariant derivatives of all orders, as obstructions for the affine extension to be holonomic. In addition this framework provides a nice interpretation for the absence of other tensors.
\end{abstract}

\tableofcontents

\intro{Introduction}

The heart of this paper consists in the observation that any torsionless affine connection $\nabla$ on a smooth manifold $M$ can be described in terms of its geodesic symmetries $(s_x)_{x\in M}$ (cf.~paragraph preceding \lref{thelemma}) by the formula~:

\begin{equation}\label{conn2}
\Bigl(\nabla_X Y\Bigr)_x = \frac{1}{2} \Bigl[X, Y + (s_x)_*Y\Bigr]_x,
\end{equation}
where $x \in M$ and $X$, $Y$ are vector fields on $M$. \\

When the ambient manifold is endowed with the structure of a symplectic symmetric space and $\nabla$ is the unique connection admitting the symmetries as affine transformations, a similar formula, involving the symplectic structure as well as the symmetries, appeared in \cite{B}  (see \rref{sss} for more details). \\

D'Atri, in \cite{D'A} gave a proof of the consequence of the previous formula that the collection of geodesic symmetries determines the connection. More precisely, he wrote the expression (\ref{Christoffel}) that defines the Christoffel symbols in terms of the geodesic symmetries and proved that the $\Gamma$'s so defined transform as they are expected to under a change of coordinates. D'Atri raised the question of characterizing germs of geodesic symmetries amongst germs near the diagonal of maps $s : M \times M \to M$ for which each $s(x, \cdot)$ is an involutive diffeomorphism of a neighborhood of $x$ having $x$ as an isolated fixed point. \sref{geod-sym&loc-sym-sp} provides a partial answer to that question. It is partial in that we do not provide an algebraic characterization of the geodesic symmetries. Rather, they appear as weak integral submanifolds of a (non-necessarily integrable) distribution constructed from the data of the family $(j^2_xs_x)_{x \in M}$, itself arbitrary (except for the condition  $j^1_xs_x = -I$ for all $x \in M$).\\

The expression (\ref{conn2}) relies only on the family of second order jets $(j^2_xs_x)_{x\in M}$ of the geodesic symmetries. Notice that the first order part of the jet $j^2_xs_x$ coincides with the linear map $-I_x: T_xM \to T_xM : X \to -X$. It is also true that a family of $2$-jets whose first order coincides with the family $-I = (-I_x)_{x \in M}$, called in the text a {\bf holonomic symmetry jet}, induces a torsionless affine connection. Moreover this correspondence between torsionless affine connections and holonomic symmetry jets is one-to-one (\tref{conn<-->sym-jet}). \\

This procedure may be enlarged to encompass affine connections with torsion. One only needs to consider semi-holonomic symmetry jets as well. Recall that a \emph{non-holonomic} $(1,1)$-jet $j^1_xb$ is the first jet at a point $x$ of a family $b(x') = j^1_{x'}f_{x'}$ of $1$-jets (cf.~\cite{Ehresmann}, \cite{Lib}, see also \aref{jet-of-bis}). It is said to be \emph{semi-holonomic} when $b(x) = j^1_xb^0$ with $b^0(x') = f_{x'}(x')$ and  \emph{holonomic} when $j^1_xb$ is a genuine $2$-jet. As stated in \pref{symm-jet<-->conn}, there is a bijective correspondence between semi-holonomic symmetry jets and affine connections that extends the above-mentioned correspondence for torsionless affine connections. \\

One of our purpose from thereon as been to understand how the various actors of the theory of affine connections are related to the symmetry jet. Often so, although not always, the results appearing in the text are known, but their formulation is geometric and differs from the usual one. Moreover, all proofs are thoroughly intrinsic, shading light on the true nature of the objects involved and their relationship to one another. A very simple instance thereof appears in \sref{UAE}, where affine transformations are described as leaves of a certain distribution on a groupoid naturally associated to $M$. Also, some results become quasi-obvious and consist simply in \emph{looking}, like the fact that for a locally symmetric space, the geodesic symmetries are affine transformations (cf.~\sref{lss}). \\

Another aim is to render precise the intuition that at order two an affine manifold $(M, \nabla)$ should be a locally symmetric space in the sense that at each point $x$ in $M$, there should exist a canonical germ of locally symmetric space osculatory to $M$ at $x$. The treatment of affine connections carried on in the text seems to bring us closer to that goal.\\

Coming back to the content of the paper, two important actors in the theory of connections are the torsion and the curvature tensors. A natural and short path from symmetry jets to these goes by the property of uniqueness of affine extension and its reformulation in terms of groupoids. Let us recall the classical statement. \\

{\bf Uniqueness of affine extension~:} On a manifold $M$ endowed with a torsionless affine connection $\nabla$, any linear isomorphism $\xi : T_xM \to T_yM$ may be uniquely extended to an affine $2$-jet. That is, there exists a $2$-jet $j^2_xf$ that satisfies $j^1_xf = \xi$ and 
$$f_{*_x} \Bigl(\nabla_{X_x} Y \Bigr) = \nabla\bigl._{f_{*_x} X_x} f_* Y.$$
We show that this property is valid for affine connections with torsion as well and moreover admits a nice reformulation and proof in terms of groupoids. To clarify our point, we need a couple of paragraphs. \\ 

The general linear groupoid of the tangent bundle of $M$, also called the gauge groupoid of the frame bundle and commonly denoted by $GL(TM) \rightrightarrows M$, is the set of all linear isomorphisms $\xi : T_xM \to T_yM$ between tangent spaces to $M$ (for a good treatment of groupoids, we recommend either \cite{Mck-87}, \cite{Mck-05} or \cite{dSW}). The source and target of $\xi$, denoted respectively by $\alpha(\xi)$ and $\beta(\xi)$, are $x$ and $y$. Observe that it is indeed a groupoid as only those pairs $(\xi, \xi')$ of linear isomorphisms for which the source of $\xi$ coincides with the target of $\xi'$ can be composed. All properties of Lie groupoids are satisfied by $GL(TM) \rightrightarrows M$ (see \cite{Mck-05}, 1.1.12). \\

The first jet extension of a Lie groupoid $G \rightrightarrows M$, denoted by $\b^{(1)}(G)\rightrightarrows M$, is defined to be the collection of all $1$-jets of local bisections of $G$ (see \aref{jet-of-bis}). It is naturally endowed with the structure of a Lie groupoid over $M$. The notion of $k$-jet extension is defined similarly. This procedure can be iterated and the groupoids $\b^{(l)}(\b^{(k)}(G))$,  $\b^{(m)}(\b^{(l)}(\b^{(k)}(G)))$, ... are denoted by $\b^{(l,k)}(G)$, $\b^{(m,l,k)}(G)$,... respectively. There are some noticeable inclusions like~: $\b^{(2)}(G) \subset \b^{(1,1)}(G)$. Now if $G$ is the pair groupoid $\PP(M) = M \times M \rightrightarrows M$, then $\b^{(1)}(G)$ coincides with $GL(TM)$ and a semi-holonomic symmetry jet is a section $\fs : M \to \b^{(1,1)}(\PP(M))$ whose first order part is $-I$. 
 \\

Now we can reformulate the property of uniqueness of affine extension in the way that is useful for us. \\

\noindent
{\bf Uniqueness of affine extension bis~:} If $\fs : M \to \b^{(1,1)}(\PP(M))$ is a symmetry jet there exists a unique groupoid morphism 
$$S : \b^{(1)}(\PP(M)) \to \b^{(1,1)}(\PP(M))$$ 
that satisfies $S \circ -I = {\mathfrak s}$. Moreover the image of $S$ consists of all affine jets. \\

The map $S$ plays a central role here. It can be thought of as a distribution on $\b^{(1)}(\PP(M))$ through the simple and well-known observation that the data of a $(1,1)$-jet $j^1_xb$ is equivalent to that of the plane $D(j^1_xb) = b_{*_ x} (T_xM)$ tangent to $\b^{(1)}(\PP(M))$. Thus the map $S$ is, in other guise, the distribution $\d$ on $\b^{(1)}(\PP(M))$ defined by $\d_\xi = D(S(\xi))$ and called the affine distribution. Now affine transformations of $M$ appear as leaves of $\d$ that are bisections of $\b^{(1)}(\PP(M))$.\\

If the connection is torsionless, an affine jet is necessarily holonomic, i.e.~is a genuine $2$-jet. If the connection has torsion, this is not anymore true but some affine jets may still be holonomic. A semi-holonomic $(1,1)$-jet is holonomic when it is fixed under the canonical involution of $\b^{(1,1)}(\PP(M))$ that permutes the derivations. The torsion appears as a measure of the defect of holonomy as the next result makes clear~: \\

{\bf \pref{prop-torsion}} The jet $S(\xi)$ is holonomic exactly when $\xi$ preserves the torsion $T$. More precisely, if  $X$ and $Y$ are vector fields on $M$ and $x\in M$, setting $\X = Y_{*_x} X_x \in T_{Y_x}TM$, the relation
\begin{equation}\label{intro-torsion}
S(\xi) \cdot \X - \kappa(S(\xi)) \cdot \X = \xi\Bigl(T(X_x, Y_x)\Bigr) - T\Bigl(\xi(X_x), \xi(Y_x)\Bigr)
\end{equation}
holds, where the $\cdot$ refers to the natural action of $\b^{(1,1)}(\PP(M))$ on $T^2M$ (derived from the action of $\b^{(1)}(\PP(M))$ on $TM$). \\

Particularizing (\ref{intro-torsion}) to $\xi = -I_x$, we obtain a formula for the torsion in terms of the symmetry jet~${\mathfrak s}$~:
$$T (X_x, Y_x) = \frac{1}{2} \Bigr(\kappa\bigl({\mathfrak s}(x)\bigr) \cdot \X - {\mathfrak s}(x) \cdot \X \Bigr).$$


Since the curvature involves the second order covariant derivative, it seems natural to investigate affine extensions of the next order. It is not difficult to prove inductively that affine extensions of all order exist (\pref{uaebis} and first paragraph of \sref{why-there-are-no}). In particular, there exists a groupoid morphism 
$${\mathbb S} : \b^{(1)}(\PP(M)) \to \b^{(1,1,1)}(\PP(M))$$ 
whose image is the collection of all affine $(1,1,1)$-jets, i.e.~jets that preserve the second covariant derivative with respect to $\nabla$ (\dref{affine111jet}). Here again an affine $(1,1,1)$-jet need not be a $3$-jet and two natural involutions on $\b^{(1,1,1)}(\PP(M))$, denoted $\kappa$ and $\kappa_*$, generate the group $S_3$ of all permutations of the derivations. The $3$-jets are the fixed points of both $\kappa$ and $\kappa_*.$ \\

Now an affine $(1,1,1)$-jet ${\mathbb S}(\xi)$, for which $S(\xi)$ is holonomic, is $\kappa$-invariant if and only if $\xi$ preserves the curvature tensor $R$. More precisely, the following formula appears in \pref{prop-curvature}~:
\begin{equation}\label{curvature}
{\mathbb S}(\xi) \cdot {\mathfrak X} - \kappa ({\mathbb S}(\xi)) \cdot {\mathfrak X} = \xi \Bigl(R(X_x, Y_x)Z_x\Bigr) - R(\xi X_x, \xi Y_x) \; \xi Z_x,
\end{equation}
where $\IX = Z_{**_{Y_x}} Y_{*_x}X_x \in T_{Z_{*_x}Y_x}T^2M$ for vector fields $X$, $Y$, $Z$ on $M$ and some point $x$ in $M$. Here $\cdot$ stands for the action of $\b^{(1,1,1)}(\PP(M))$ on $T^3M$. \\

How about invariance under $\kappa_*$~? What is the obstruction~? The first covariant derivative of the torsion tensor as the following formula, a direct consequence of (\ref{intro-torsion}), implies~:
$${\mathbb S}(\xi) \cdot {\mathfrak X} - \kappa_*({\mathbb S}(\xi)) \cdot {\mathfrak X} = \xi \Bigl( (\nabla_{Z_x}T^\nabla) (X_x, Y_x) \Bigr) - \Bigl(\nabla_{\xi Z_z}T^\nabla\Bigr)(\xi X_x, \xi Y_x).$$

To summarize, an affine $(1,1,1)$-jet is a genuine $3$-jet if and only if its first order part preserves the torsion, its first covariant derivative and the curvature. \\

Setting $\xi = a I$, with $a \neq 1, -1$, in the first of these relations, we obtain a description of the curvature of a torsion free connection $\nabla$ in terms of the affine extension map~:
$$R(X_x, Y_x)\,Z_x = \frac{1}{a (1 - a^2)} \Bigl({\mathbb S}(aI_x) \cdot {\mathfrak X} - \kappa ({\mathbb S}(aI_x)) \cdot {\mathfrak X}\Bigr).$$ 
The following alternative description, whose advantage is that it involves only the symmetry jet, is proved in \tref{curvthm}~:
$$R(X_x, Y_x) Z_x = \frac{1}{4} \Bigl(\kappa(j^1_x {\mathfrak s}) \cdot j^1_x {\mathfrak s} \cdot \IX -  \;j^1_x {\mathfrak s} \cdot \kappa(j^1_x {\mathfrak s}) \cdot \IX \Bigr).
$$
It is tempting (and even true provided the symbols are suitably defined) to write 

$$R = \frac{1}{4}\bigl[ \kappa(j^1\fs), j^1\fs\bigr].$$

We said earlier that affine extensions exist at all order. A natural question therefore is~: ``what happens at order $4$ and higher~?" The  $\kappa$-invariance of an affine jet of order at least $4$ is automatic. A simple explanation for that phenomenon, that relies on the formalism developed previously, is exposed in \sref{why-there-are-no}. It provides in particular an elementary justification for the well-known fact that there is no other natural tensor in affine geometry than the torsion and the curvature tensors (as well as their covariant derivatives of all orders of course). \\


Given a family of $TM$-valued covariant tensors $\{Q_i; i \in I\}$ on $M$, the closed set of $1$-jets that preserve all of them is a Lie subgroupoid of $\b^{(1)}(\PP(M))$ denoted hereafter by $\b(\{Q_i; i \in I\})$. Thus, considering the torsion and curvature tensors $T$ and $R$ of an affine connection, a collection of Lie subgroupoids of $\b^{(1)}(\PP(M))$ naturally appears~: $\b(T)$, $\b(T, \nabla T, R)$, etc.... Their intersection 
$$\b_o = \b(T, R, ..., \nabla^kT, \nabla^kR, ...)$$ 
is the largest subset of $\b^{(1)}(\PP(M))$ entirely foliated by leaves (of maximal dimension) of the associated affine distribution $\d$ and containing all of them.\\ 

In particular, when the affine connection is locally symmetric, that is satisfies $T = 0$ and $\nabla R=0$, the groupoid $\b_o$ coincides with $\b(R)$. Thus through any $1$-jet $\xi$ that preserves $R$ passes a leaf of $\d$ which is necessarily the $1$-jet extension of a partial affine transformation $\varphi_\xi$ of $(M, \nabla)$. As proved in \sref{geod-sym&loc-sym-sp}, the maps $\varphi_\xi$ is the geodesic extension of $\xi$, that is

$$\varphi_\xi = \exp_y \circ \; \xi \circ \exp_x^{-1}.$$

As a by-product, we can reprove existence and uniqueness of the Levi-Civita connection of a pseudo-Riemannian metric (\pref{LC}). The specificity of pseudo-Riemannian metrics, compared to other tensors like symplectic forms, is that the collection of $2$-jets extending $-I_x$ has the same dimension as the set of $1$-jets of symmetric tensors extending a given non-degenerate symmetric tensor~$g_x$. \\

We also compare our approach with Kobayashi's correspondence between torsionless affine connections and admissible sections (\cite{K}). The latter can easily be extended to connections with torsion by enlarging the class of admissible sections and we show how an ``admissible'' section naturally induces a symmetry jets and vice-versa. In the same spirit, affine connections are a class of Lie algebroid connections and we describe the correspondance between the latter and symmetry jets. \\

The paper is organized as follows. There is a large appendix that contains all the relevant material with our preferred notation about groupoids and groupoids of jets of bisections. The detailed structure of $\b^{(1,1)}(\PP(M))$ and $\b^{(1,1,1)}(\PP(M))$, as well as that of the second and third tangent bundles $T^2M$ and $T^3M$ is investigated there. The idea is to read the appendix alongside with the main text. As to the latter, we begin, in Section 1, with the correspondence between torsionless affine connections and holonomic symmetry jets. The case of affine connections with torsion is treated in Section 2. Section 3 contains a reformulation of the property of uniqueness of affine extensions in terms of groupoids and introduces the distribution $\d$. The correspondence between the defect of holonomy of the affine extension and the defect of invariance of the torsion is treated in Section 4. The description (\ref{curv-j1s}) of the curvature in terms of the first jet extension of the symmetry jet is established in Section 5. The  property of uniqueness of affine extension at order $3$ is handled in Section 6. Section 7 proves the correspondence between the defect of $\kappa$-invariance (\rp $\kappa_*$-invariance) of the affine extension of order $3$ of a $1$-jet $\xi$ and the defect of invariance of the curvature (\rp first covariant derivative of the torsion) under $\xi$ and describes geometrically the integrability locus of $\d$. Section 8 investigates the question of existence and holonomy of order $4$ affine extensions. It is proven there that the $\kappa$-invariance does automatically hold, explaining why no new tensor appears. Section 9 describes the horizontal distribution on an associated bundle induced from an affine connection. An alternative proof of existence and uniqueness of the Levi-Civita connection is presented in Section 10. A geometric correspondence  between a symmetry jet and the Lie algebroid connection associated to the induced affine connection is shown in Section 11. The geometric description of parallel transport and geodesics appears in Section 12. Whence follows, in Section 13, the construction in terms of $\d$ of the $1$-jet extension of the map $\varphi_\xi = \exp_y \circ \; \xi \circ \exp_x$ for some linear isomorphism $\xi : T_xM \to T_yM$. That section also contains the  proof of the property that for a locally symmetric space, through any $\xi$ in $\b^{(1)}(\PP(M))$ that preserves the curvature passes a leaf of $\d$, which is then necessarily $j^1\varphi_\xi$. Section 14 deals with Kobayashi's correspondence between torsionless affine connections and admissible sections of the second order frame bundle. 

\section*{Acknowledgments} 
We wish to warmly thank Wolfgang Bertram, Wilhelm Klingenberg, the late Shoshichi Kobayashi, Alan Weinstein and Joe Wolf for very stimulating conversations regarding the content of the paper. We are grateful to Yannick Voglaire for pointing out D'Atri's article to us. Last but not least, Kirill MacKenzie has made numerous thoughtful suggestions to improve the presentation of the paper. We thank him heartily. 

\section{Torsionless affine connections as symmetry jets}

An affine connection on a smooth manifold $M$ is commonly defined to be a $\R$-bilinear map
$$\nabla: \IX(M) \times \IX(M) \to \IX(M) : (X, Y) \mapsto \nabla_X Y$$
which is $C^\infty(M)$-linear in the first argument and satisfies the Leibniz rule in the second argument, that is, for all $X$, $Y$ in $\IX(M)$ and all $f \in C^\infty(M)$, we have
\begin{enumerate}
\item[-] $\nabla_{fX}Y = f \nabla_X Y$,
\item[-] $\nabla_X(fY) = Xf \; Y + f \nabla_X Y$.
\end{enumerate}
The previous relations suggest to think of $(\nabla_X Y)_x$ as a derivative of $Y$ in the direction of $X_x$. Such a derivative exists already~: it is simply $Y_{*_x}X_x$ ($Y$ is thought of as a section $M \to TM$ of the tangent bundle). The point is that $Y_{*_x}X_x$ lies in the second tangent bundle $T^2M$, while we would like to have an element of $TM$. In fact, an affine connection is really a projection from $T^2M$ to $TM$ and the associated horizontal distribution on $TM$ is its ``kernel". A precise statement that relies on notations introduced in \aref{second-t-b} is the content of the next lemma. 

\begin{lem}\label{tilde-nabla} An affine connection on the manifold $M$ is  a map
$$\widetilde{\nabla} : T^2M  \to TM : \X \to \widetilde{\nabla}(\X)$$
such that
\begin{enumerate} 
\item[$\bullet$] $\widetilde{\nabla}$ is a morphism of vector bundles over the map $p : TM \to M$ when $T^2M$ is endowed with the vector bundle structure associated to either $p$ or $p_*$.
\item[$\bullet$] $\widetilde{\nabla} (i^p_{0_x}(V_x)) = V_x$.
\end{enumerate}
The correspondence with the classical definition of affine connection goes through the relation~:
\begin{equation}\label{alt-conn}
\nabla_{X_x}Y = \widetilde{\nabla}(Y_{*_x}X_x).
\end{equation}
\end{lem}

\Pf Given an affine connection $\nabla$, it induces a map $\widetilde{\nabla} : T^2M \to TM$ defined by (\ref{alt-conn}) on non vertical vectors which can be extended by $\widetilde{\nabla}(i^p_{0_x}(Y_x)) = Y_x$ on vertical vectors. 
It is not difficult to verify the linearity conditions except for the case of two vectors $\X_1$ and $\X_2$ in the same $p$-fiber and whose $p_*$-projections are linearly dependent (as we may not assume that $\X_1$ and $\X_2$ are of type $Y_{*}(X)$ for a same local vector field $Y$). Suppose then that $\X_1 = {Y_1}_{*}(X_x)$ and $\X_2 = {Y_2}_*(a X_x)$ for some $a \in \R_0$. Then it is not difficult to verify that 
$$\X_1 + \X_2 = (1+a) \frac{d}{dt} \Bigl( \frac{1}{(1+a)} Y_1^t + \frac{a}{(1 + a)} Y_2^t \Bigr)\Bigr|_{t = 0}.$$
Hence $\X_1 + \X_2 = Y_*((1+a)X_x)$ with 
$$Y = \frac{1}{1+a} Y_1 + \frac{a}{1 + a} Y_2.$$ 
Verifying now that $\widetilde{\nabla}(\X_1 + \X_2) = \widetilde{\nabla}(\X_1) + \widetilde{\nabla}(\X_2)$ is easy. \\

Conversely, given a map $\widetilde{\nabla}$ as in the statement of the lemma, defining $\nabla$ through (\ref{alt-conn}) yields an affine connection. Indeed, the Leibniz rule is the only point that might not seem to follow immediately. It is a direct consequence of \rref{Leibniz} according to which 
$$(fY)_{*_x} (X_x) = X_xf \; Y_x + m_{f(x)*}\Bigl( Y_{*_x}X_x \Bigr),$$
where $X_x f \; Y_x$ really means $I(f(x)Y_x, X_xf Y_x) = i(f(x)Y_x) +_* i_{0_x}(X_xf \; Y_x)$.
\cqfd

\begin{rmk}\label{hor-dist}
The horizontal distribution $\h = \h^{\nabla}$ associated to the connection $\nabla$ is the ``kernel" of $\widetilde{\nabla}$, that is 
$$\h = \widetilde{\nabla}^{-1}(0_{TM}).$$ 
\end{rmk}

Now, let $(s_x)_{x \in M}$ be a smooth family of smooth local diffeomorphisms $M$ such that $s_x$ is defined near $x$ and satisfies $s_x(x) = x$ and ${s_{x_{*_x}}} = - \id$. We consider the bilinear map $\nabla : \IX(M) \times \IX(M) \to \IX(M) : (X, Y) \mapsto \nabla_XY$ defined by the formula
\begin{equation}\label{connection}
\Bigl(\nabla_X Y\Bigr)_x = \frac{1}{2} \Bigl[X, Y + (s_x)_*(Y)\Bigr]_x, 
\end{equation}
where $x \in M$ and $X, Y \in {\mathfrak X}(M)$. To be precise, the vector field $Y + (s_x)_*(Y)$ achieves the value $Y_{x'} + (s_x)_{*_{s_x^{-1}(x')}}(Y_{s_x^{-1}(x')})$ at the point $x'$. 

\begin{prop} Formula (\ref{connection}) defines a torsionless affine connexion on $M$. 
\end{prop}

\Pf Let us verify that the three conditions defining a torsionless affine connection are satisfied. First of all $\nabla_{fX}Y = f\nabla_XY$ because the vector field $Y + (s_x)_*(Y)$ vanishes at $x$. \\

To prove the condition $\nabla_{X} fY = Xf \; Y + f \nabla_X Y$, observe that $(s_x)_*(fY) = (f \circ s_x^{-1}) (s_x)_*Y$. Hence 
$$[X, (s_x)_*(fY)] = X(f \circ s_x^{-1}) (s_x)_*Y + (f \circ s_x^{-1})[X, (s_x)_*Y],$$ 
which evaluated at $x$ yields 
$$[X, (s_x)_*(fY)]_x = X_x f Y_x + f(x) [X, (s_x)_*Y]_x.$$ 
Then 
$$\begin{array}{ccl}
2\Bigl(\nabla_XfY \Bigr)_x & = & \Bigl[X, fY + (s_x)_*(fY)\Bigr]_x\\
& = & X_xf Y_x + f(x) \Bigl[X, Y\Bigr]_x + X_x f Y_x + f(x) \Bigl[X, (s_x)_*(Y)\Bigr]_x\\
& = & 2 X_xf \; Y_x +  2 f(x) \Bigl(\nabla_XfY \Bigr)_x 
\end{array}$$

Finally the torsion $T^\nabla(X,Y) = \nabla_X Y - \nabla_Y X - [X, Y]$ vanishes because
$$\begin{array}{ccl}
\Bigl(\nabla_X Y - \nabla_Y X \Bigr)_x & = & \frac{1}{2} \Bigl[X, Y + (s_x)_*Y\Bigr]_x + \frac{1}{2} \Bigl[X + (s_x)_*X, Y\Bigr]_x \\
& = & \Bigl[X, Y]_x + \frac{1}{2} \Bigl[X + (s_x)_*(X), Y + (s_x)_*(Y)\Bigr]_x \\
& = & \Bigl[X, Y\Bigr]_x.
\end{array}$$
\cqfd

\begin{rmk}
The Christoffel symbols of the connexion $\nabla$ with respect to local coordinates $x^1, ..., x^n$ around $x$ are
\begin{equation}\label{Christoffel}
\Gamma_{ij}^k (x) =  - \frac{1}{2} \frac{\partial^2s_x^k}{\partial x^i \partial x^j} (x).
\end{equation}
\end{rmk}

Observe that, of the whole family $(s_x)_{x \in M}$, only $(j^2_xs_x)_{x \in M}$ consisting of the second-order jet of each $s_x$ at $x$ plays a role. In other words, any section 
$${\mathfrak s} : M \to \b^{(2)}(\PP(M))$$ of the bundle $\alpha : \b^{(2)}(\PP(M)) \to M$ of $2$-jets of local diffeomorphisms of $M$ (cf.~\nref{b1} in \aref{jet-of-bis}) that projects onto the section 
$$-I : M \to \b^{(1)}(\PP(M)) : x \mapsto [- I_x : X \mapsto -X]$$ 
via the canonical projection $p : \b^{(2)}(\PP(M)) \to \b^{(1)}(\PP(M))$ determines a connexion $\nabla^{\mathfrak s}$. 

\begin{dfn} A section ${\mathfrak s} : M \to \b^{(2)}(\PP(M))$ such that $p  \circ  {\mathfrak s} = -I$ is called hereafter a holonomic symmetry jet.
\end{dfn}

\begin{rmk}\label{sss} As described in \cite{B}, the canonical connection of a symplectic symmetric space $(M, (s_x)_{x \in M}, \omega)$ admits an expression similar to (\ref{connection}). More precisely, if $X$, $Y$ and $Z$ are vector fields on $M$, then the following expression defines both the unique $s_x$-invariant symplectic connection on the symplectic symmetric space $(M, (s_x)_{x \in M}, \omega)$ and the connection $\nabla^{\mathfrak s}$ associated to the symmetry jet ${\mathfrak s}(x) = j^2_xs_x$~:

\begin{equation}
\omega_x\Bigl(\nabla_X Y, Z\Bigr) = \frac{1}{2} X_x\Bigl(\omega(Y + (s_x)_*Y, Z)\Bigr).
\end{equation} 
Indeed, it is a consequence of the following short computation
$$\begin{array}{ccl}
0 & = & \bigl(\nabla_X\omega\bigr)_x\Bigl(Y + (s_x)_*Y, Z\Bigr) \\
& = & X_x\Bigl(\omega(Y + (s_x)_*Y, Z)\Bigr) - \omega_x\Bigl(\nabla_X(Y + (s_x)_*Y), Z\Bigr) \\
& = & X_x\Bigl(\omega(Y + (s_x)_*Y, Z)\Bigr) - \omega_x\Bigl([X, Y + (s_x)_*Y], Z\Bigr).
\end{array}$$
\end{rmk}

\begin{dfn} A diffeomorphism $\varphi$ of a manifold $M$ endowed with an affine connection $\nabla$ is said to be affine if 
$$\varphi_* \Bigl(\nabla_XY\Bigr) = \nabla_{\varphi_*X} \varphi_*Y$$ 
holds for all $X$, $Y$ in $\IX(M)$. Likewise the $2$-jet $j^2_x\varphi$ of a local diffeomorphism $\varphi : U \subset M \to V\subset M$ at a point $x$ of its domain is said to be affine if the previous relation holds at $x$, that is~:
$$\varphi_{*_x} \Bigl(\nabla_{X_x}Y\Bigr) = \nabla_{\varphi_{*_x}X_x} \varphi_*Y.$$ 
\end{dfn}

\begin{lem}\label{inv} The connexion $\nabla^{\mathfrak s}$ admits ${\mathfrak s}(x)$ as affine $2$-jet.
\end{lem}

\Pf $$\begin{array}{ccl}
2\Bigl(\nabla^{\mathfrak s}_X Y\Bigr)_x & = & \Bigl[X, Y + (s_x)_*Y\Bigr]_x \\
& = & - (s_x)_{*_x} \Bigl[X, Y + (s_x)_*Y\Bigr]_x \\
& = & - \Bigl[(s_x)_*X, (s_x)_*Y + (s_x)_* \circ  (s_x)_*Y\Bigr]_x \\
& = & - \Bigl[- X, (s_x)_*Y + (s_x)_* \bigl( (s_x)_*Y \bigr) \Bigr]_x \\
& = & 2 \Bigl(\nabla^{\mathfrak s}_X (s_x)_*Y\Bigr)_x,
\end{array}$$
\cqfd

\begin{rmk}
Given a whole family $(s_x)_{x \in M}$ of smooth diffeomorphisms of $M$ satisfying $s_x(x) = x$ and $(s_x)_{*_x} = -\id$, there is not a global $s_x$-invariance of $\nabla^{\mathfrak s}$, for ${\mathfrak s}(x) = j^2_xs_x$, unless the $s_x$'s satisfy additional relations of a global nature. A typical example is a symmetric space, where the symmetries $s_x$ satisfy $s_x \circ s_x = \id$ and $s_y \circ s_x \circ s_y = s_{s_y(x)}$, In that case, the associated connection is globally $s_x$-invariant. Indeed, 
$$\begin{array}{ccl}
\Bigl(\nabla_{(s_y)_*X} (s_y)_*Y\Bigr)_x 
& = & \frac{1}{2} \Bigl[(s_y)_*X, (s_y)_*Y + (s_x)_* \circ (s_y)_*Y\Bigr]_x \\
& = & \frac{1}{2} \Bigl[(s_y)_*X, (s_y)_*Y + (s_y)_* \circ (s_{s_y(x)})_*Y\Bigr]_x \\
& = & \frac{1}{2} (s_y)_{*_{s_y(x)}}\Bigl[X, Y + (s_{s_y(x)})_*Y\Bigr]_{s_y(x)} \\
& = & (s_y)_{*_{s_y(x)}} \Bigl(\nabla_{X} Y\Bigr)_{s_y(x)}  
\end{array}$$
 \end{rmk}

We have obtained so far a connection from a symmetry jet. On the other hand, a connection induces a family of local diffeomorphisms, its geodesic symmetries. More precisely, let $\exp : \OO \subset TM \to TM$ be the exponential map associated to the connection $\nabla$, that is the map that sends a tangent vector $X$ to the time-one geodesic tangent to $X$. Here $\OO$ is assumed to be some neighborhood of the $0$-section in $TM$ on which $\exp$ is defined and such that if $\OO_x$ denotes the intersection of $\OO$ with $T_xM$, then $-\OO_x = \OO_x$ and the restriction of $\exp$ to $\OO_x$ --- denoted by $\exp_x$ --- is a diffeomorphism onto some open subset $\U_x$ of $M$. The geodesic symmetry at $x$ associated to $\nabla$ is the local involutive diffeomorphism $a^\nabla_x : \U_x \to \U_x : y \to \exp(- \exp^{-1} (y))$. As can be expected, the $2$-jet at $x$ of the geodesic symmetry of the connection $\nabla^{\mathfrak s}$ coincides with ${\mathfrak s}$. 

\begin{lem}\label{thelemma} If $a_x$ denotes the geodesic symmetry at $x$ induced by the connection $\nabla^{\mathfrak s}$ associated to the symmetry jet ${\mathfrak s}$, then $$j^2_x a_x = {\mathfrak s}(x).$$
\end{lem}

\Pf Let $\gamma(t) = \exp_x(tX_x)$ be a geodesic with tangent vector field $X_{\gamma(t)} = \frac{d\gamma}{dt}(t)$. The latter satisfies  
$$0 = \nabla_{X_x} X = \frac{1}{2} \Bigl[X_x, X + (s_x)_*X\Bigr]_x.$$
Equivalently 
$$0 = X_x\Bigl((X +  (s_x)_*X)(f)\Bigr) = X_x(Xf) + X_x\Bigl( (s_x)_*X(f)\Bigr),$$
for all $f \in C^\infty(M)$. Developing the right hand side of the previous equality yields
$$\begin{array}{rcl}
0 & = & \displaystyle{X_x(Xf) - (s_x)_{*_x}X_x \Bigl((s_x)_*X(f)\Bigr)} \\
\vspace{0.2cm}
& = & - \displaystyle{\frac{d}{dt}\exp_x(-tX_x)\Bigl|_{t=0} (Xf) - \frac{d}{dt} s_x \circ \exp_x(tX_x)\Bigl|_{t=0} (s_x)_*X(f)} \\
\vspace{0.2cm}
& = & \displaystyle{\frac{d}{dt} -X_{\exp_x(-tX_x)} (f) \Bigl|_{t=0} - \;\frac{d}{dt}( (s_x)_*X)_{s_x \circ \, \exp_x(tX_x)}(f)\Bigl|_{t=0}} \\
\vspace{0.2cm}
& = & \displaystyle{\frac{d}{dt} \frac{d}{ds} f \circ \exp_x(-sX_x)\Bigl|_{s=t} \Bigl|_{t=0} - \;\frac{d}{dt}(s_x)_*(X_{\exp_x(tX_x)})(f)\Bigl|_{t=0}} \\
\vspace{0.2cm}
& = & \displaystyle{\frac{d^2}{dt^2} f \circ \exp_x(-tX_x)\Bigl|_{t=0} - \;\frac{d}{dt} X_{\exp_x(tX_x)} (f \circ s_x)\Bigl|_{t=0}} \\
\vspace{0.2cm}
& = & \displaystyle{\frac{d^2}{dt^2} f \circ \exp_x \circ -I_x (tX_x)\Bigl|_{t=0} - \;\frac{d}{dt} \frac{d}{ds} f \circ s_x \circ \exp_x (sX_x)\Bigl|_{s=t}\Bigl|_{t=0}} \\
\vspace{0.2cm}
& = & \displaystyle{\frac{d^2}{dt^2} f \circ \exp_x \circ -I_x (tX_x)\Bigl|_{t=0} - \;\frac{d^2}{dt^2}f \circ s_x\circ \exp_x(tX_x)\Bigl|_{t=0}}.
\end{array}$$
We claim that this implies that the maps $f \circ \exp_x \circ -I_x$ and $f \circ s_x \circ \exp_x$ coincide up to order $2$. Indeed, observe that the differential at $0_x \in T_xM$ of these two maps coincide. Hence their second differential (cf.~\nref{f**} in \aref{b11})
$$f_{**_x} \circ {\exp_x}_{**_{0_x}} \circ (-I_x)_{**_{0_x}} \quad \mbox{and} \quad f_{**_x} \circ {s_x}_{**_x} \circ {\exp_x}_{**_{0_x}}$$ 
belong to a same fiber of $p \times p_* : \b^{(2)}(\PP(M)) \to \b^{(1)}(\PP(M)) \times \b^{(1)}(\PP(M))$, so that their difference is a symmetric bilinear map $B$ from $T_xM \times T_xM$ to $T_yM$ (cf.~\rref{affine-str}) which is determined by its values on pairs $(X, X) \in T_xM \times T_xM$. Now, the previous calculation shows that $B(X, X)$ vanishes for all $X$. Whence the result.
\cqfd 

\begin{rmk}
So starting from a symmetry jet ${\mathfrak s}$, we obtain in a canonical way, via the affine connection $\nabla^{\mathfrak s}$, a smooth family $(a^{\nabla^{\mathfrak s}}_x)_{x \in M}$ of local involutive diffeomorphisms which integrate pointwise the section~${\mathfrak s}$.
\end{rmk}

Given a torsionless affine connection $\nabla$, let us denote by ${\mathfrak a}^\nabla$ the symmetry jet defined by 
$${\mathfrak a}^\nabla (x) = j^2_x a^\nabla_x.$$ 

\begin{thm}\label{conn<-->sym-jet}
The two correspondences ${\mathfrak s} \rightsquigarrow \nabla^{\mathfrak s}$ and $\nabla \rightsquigarrow {\mathfrak a}^\nabla$ are inverse to one another. In particular it is true that any affine connection $\nabla$ is associated to the symmetry jet consisting of the family of 2-jets  $a^\nabla$ of its  geodesic symmetries, through the relation 
\begin{equation}\label{conn/sym}
\Bigl(\nabla_X Y\Bigr)_x = \frac{1}{2} \Bigl[X, Y + (a^\nabla_x)_*Y\Bigr]_x, \quad X, Y \in {\mathfrak X}(M)
\end{equation}
\end{thm}

\Pf Half of the \pref{conn<-->sym-jet}, namely the fact that ${\mathfrak a}^{\nabla^{\mathfrak s}} = {\mathfrak s}$, has been proven in \pref{thelemma}. The other half, that is, 
$$\nabla^{a^\nabla} = \nabla,$$ 
is easily seen once we know that the geodesic symmetries are affine up to order $2$. Indeed, suppose $(\nabla_XY)_x = (\nabla_X(a^\nabla_x)_*Y)_x$, then
$$\begin{array}{ccl}
\Bigl(\nabla_XY\Bigr)_x & = & \frac{1}{2} \Bigl( \nabla_X Y \Bigr)_x + \frac{1}{2} \Bigl( \nabla_X (a^\nabla_x)_*Y \Bigr)_x \\
& = & \frac{1}{2} \Bigl( \nabla_X (Y + (a^\nabla_x)_*Y) \Bigr)_x \\
& = &  \frac{1}{2} \Bigl[X, Y + (a^\nabla_x)_*Y \Bigr]_x 
\end{array}$$

We prove now that the geodesic symmetries are affine up to order $2$. A linear connection $\nabla$ induces a horizontal distribution $\h = \h^\nabla$ on the tangent bundle, described in terms of $\widetilde{\nabla}$ (cf.~\lref{tilde-nabla}) as its kernel, that is, 
$$\h = \widetilde{\nabla}^{-1}(0_{TM}).$$

Therefore it is sufficient to prove that $(a_x)_{**_{Y_x}}$ maps $\h_{Y_x}$ onto $\h_{-Y_x}$. Suppose $X$ and $Y$ are vector fields defined near $x$ tangent to $\h$ at $X_x$ and $Y_x$. Then $X + Y$ is also tangent to $\h$ at $X_x + Y_x$ because the connection is linear and $[X, Y]_x = 0$ because the connection is torsionless. Besides, a vector field $Z$ that is tangent to $\h$ at $Z_x$ is also tangent to the velocity vector field of the geodesic $\exp_x(tZ_x)$, which implies that $(a_x)_*Z$ is also tangent to the velocity vector field of $\exp_x(-tZ_x)$. Thus $\nabla_{Z_x} (a_x)_*Z = 0$ and 

$$\begin{array}{ccl}
\Bigl(\nabla_X (a_x)_*Y\Bigr)_x & = & \Bigl(\nabla_{X + Y} (a_x)_*(X + Y)\Bigr)_x - \Bigl(\nabla_X (a_x)_*X\Bigr)_x \\
& & - \Bigl(\nabla_Y(a_x)_*Y\Bigr)_x - \Bigl(\nabla_Y(a_x)_*X\Bigr)_x\\
& = & - \Bigl(\nabla_Y(a_x)_*X\Bigr)_x \\
& = & \Bigl(\nabla_{(a_x)_*Y}(a_x)_*X\Bigr)_x\\
& = & \Bigl(\nabla_{(a_x)_*X}(a_x)_*Y\Bigr)_x + \Bigl[(a_x)_*Y, (a_x)_*X\Bigr]_x \\
& = & - \Bigl(\nabla_X(a_x)_*Y\Bigr)_x - \Bigl[Y,X\Bigr]_x \\
& = & - \Bigl(\nabla_X(a_x)_*Y\Bigr)_x.
\end{array}$$
Hence $\nabla_{X_x} (a_x)_*Y = 0$ for all $X_x \in T_xM$, which implies that $(a_x)_{**_x}$ preserves the horizontal distribution $\h$ along $T_xM$. 
\cqfd

\begin{rmk}\label{Kandinsky}
Notice in particular that for $f$ to be a local involutive diffeomorphism requires no condition on its second order derivative when $f_{*_x} = -I_x$. In other terms, when $j^1_xf = -I_x$, we have $j^2_xf^{-1} = j^2_xf$. (See also \cref{invertibility}).
\end{rmk}

\begin{rmk}
A alternative proof of \pref{conn<-->sym-jet} will be provided when handling the torsion case (cf.~\pref{symm-jet<-->conn}).
\end{rmk}

\section{Connections with torsion}\label{conn-torsion}

For a connection with torsion, the $2$-jet of a geodesic symmetry is not any more affine. Indeed, on the one hand, if a map is affine up to order $2$ at a point $x$, its differential at $x$ must preserve the torsion. On the other hand, the torsion being a $3$-tensor, cannot be preserved by $-I_x$ unless it vanishes at $x$. Nevertheless by relaxing slightly the notion of symmetry jet as in the following definition, one may establish a bijective correspondence between symmetry jets and arbitrary affine connections. Regarding notation, we refer to \aref{b11} and \dref{bh}.

\begin{dfn} A symmetry jet is a section 
$${\mathfrak s} : M \to \b_h^{(1,1)}(M)$$ 
of $\alpha : \b_h^{(1,1)}(M) \to M$ whose first order part $p \circ {\mathfrak s} = p_* \circ {\mathfrak s} : M \to \b^{(1)}(\PP(M))$ coincides with~$-I$. A symmetry jet is said to be holonomic if it takes its values in $\b^{(2)}(\PP(M))$ and semi-holonomic otherwise.
\end{dfn}

\begin{prop}\label{symm-jet<-->conn}
Given a symmetry jet ${\mathfrak s}$, the formula 
\begin{equation}\label{sym-jet-->conn1}
\Bigl(\nabla^{\mathfrak s}_X Y\Bigr)_x  = \frac{1}{2} \Bigl[X, Y + s_x (Y)\Bigr]_x,
\end{equation}
where ${\mathfrak s}(x) = j^1_x s_x$, for some local bisection $s_x : U_x \subset M \to \b^{(1)}(\PP(M))$, defines an affine connection. Moreover, this induces a bijective correspondence between symmetry jets and affine connections. 
\end{prop}
The proof will show that the condition that ${\mathfrak s}(x)$ belongs to $\b^{(1,1)}(\PP(M))$ rather than $\b^{(1,1)}_{nh}(\PP(M))$ ensures that the Leibniz rule is satisfied and cannot be relaxed. Of course the symmetry jet is holonomic, or $\kappa$-invariant (cf.~\lref{kappa-prop}), 
if and only if the connection is torsionless.
\begin{rmk}
Observe the following alternative expression for $\nabla^{\mathfrak s}$~:
\begin{equation}\label{sym-jet-->conn2}
\nabla^{\mathfrak s}_{X_x}Y = \frac{1}{2} \pi \Bigl(Y_{*_x} X_x, \pmb{-} {\mathfrak s}(x) \cdot Y_{*_x}X_x \Bigr),
\end{equation}
where the thick minus $\pmb{-}$ denotes the composition of the scalar multiplication by $-1$ with respect to one vector bundle structures of $T^2M$ over $TM$ with scalar multiplication by $-1$ with respect to the other one (cf.~\aref{second-t-b}), i.e.~$\pmb{-} = m_1 \circ m_{1*}$.  The expression ${\mathfrak s}(x) \cdot Y_{*_x}X_x$ stands for the action of the $(1,1)$-jet ${\mathfrak s}(x)$ on $Y_{*_x}X_x$ (cf.~(\ref{(1,1)-action})). The two vectors $Y_{*_x} X_x$ and $\pmb{-} {\mathfrak s}(x) \cdot Y_{*_x}X_x$ are in the same fiber of the affine fibration $T^2M \to TM \oplus TM$ (see the figure below). Whence their difference yields an element of $TM$ (cf.~(\ref{i}) and (\ref{def-i}) in \aref{second-t-b}). The fact that (\ref{sym-jet-->conn2}) coincides with (\ref{sym-jet-->conn1}) is an easy consequence of the relation between the Lie bracket and $\kappa$ (cf.~\pref{involution-prop})~:
$$\begin{array}{cll}
\Bigl[X, Y + s_x (Y)\Bigr]_x & = & \pi\Bigl((Y + s_x (Y))_{*_x}X_x, \kappa(X_{*_x}(Y + s_x (Y))_x) \Bigr)\\
& = & \pi\Bigl(Y_{*_x}X_x +_* {\mathfrak s}(x) \cdot Y_{*_x}(-X_x), \kappa(0_{X_x}) \Bigr) \quad \mbox{(thanks to (\ref{coucou}))}\\
& = & \pi\Bigl(Y_{*_x}X_x +_*  m_{-1} \bigl( {\mathfrak s}(x) \cdot Y_{*_x}X_x \bigr), {0_*}_{X_x} \Bigr)  \mbox{(thanks to (\ref{prop-kappa}))} \\
& = & \pi\Bigl(Y_{*_x}X_x, \pmb{-} {\mathfrak s}(x) \cdot Y_{*_x}X_x \Bigr) \quad \mbox{(thanks to (\ref{abc})).}
\end{array}$$

\begin{figure}[h!]\label{fig-conn}
\begin{center}
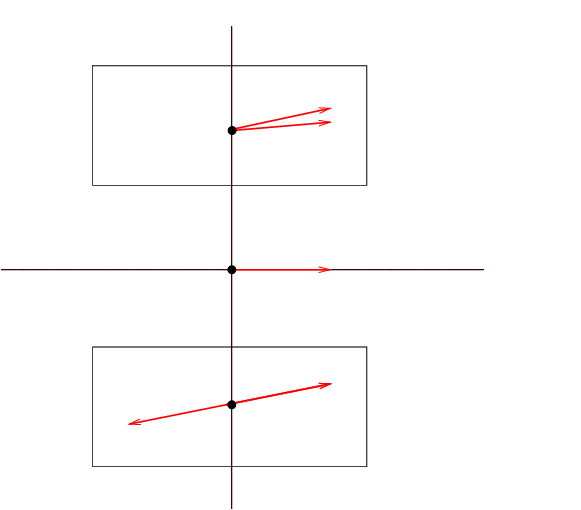 
\end{center}
\end{figure}
\end{rmk}

\noindent
{\bf Proof of \pref{symm-jet<-->conn}} Let ${\mathfrak s} : M \to \b^{(1,1)}(\PP(M))$ be a symmetry jet. Proving that (\ref{sym-jet-->conn2}) defines an affine connexion amounts to showing that the left hand side satisfies the Leibniz rule as the $C^\infty(M)$-linearity in the first argument is easy to see. So let $f$ be a smooth function on $M$, then
$$\nabla^{\mathfrak s}_{X_x}fY = \frac{1}{2} \pi \Bigl((fY)_{*_x} X_x, \pmb{-} {\mathfrak s}(x) \cdot (fY)_{*_x}X_x \Bigr).$$
Recall from (\ref{scalar-mult}) that
$$(fY)_{*_x} X_x = X_xf \; Y_x + m_{f(x)*} (Y_{*_x}X_x),$$
where $X_x f \; Y_x$ really means $I_p(f(x)Y_x, X_xf Y_x) = i(f(x)Y_x) +_* i^p_{0_M}(X_xf \; Y_x)$. So 
$$\begin{array}{ccl}
\pmb{-}{\mathfrak s}(x) \cdot (fY)_{*_x} X_x & = & \pmb{-}{\mathfrak s}(x) \cdot \Bigl(I_p(f(x)Y_x, X_xf Y_x) +  m_{f(x)*} (Y_{*_x} X_x)\Bigr) \\
& = & \pmb{-} \Bigl(I_p(- f(x)Y_x, - X_xf Y_x) +   m_{f(x)*} \bigl( {\mathfrak s}(x) \cdot Y_{*_x}X_x \bigr) \Bigr)\\
& = & \Bigl(I_p(f(x)Y_x, - X_xf Y_x) + \pmb{-}  m_{f(x)*} \bigl( {\mathfrak s}(x) \cdot Y_{*_x}X_x \bigr) \Bigr).
\end{array}$$
Then 

$$\begin{array}{ccl}
\nabla^{\mathfrak s}_{X_x}fY & = & \frac{1}{2} \pi\Bigl(I_p(f(x)Y_x, X_xf Y_x) +  m_{f(x)*} (Y_{*_x} X_x), \\
& & \;\;\;\;\;\;\; I_p(f(x)Y_x, - X_xf Y_x) + \pmb{-}  m_{f(x)*} \bigl( {\mathfrak s}(x) \cdot Y_{*_x}X_x \bigr)\Bigr) \\
& = & \frac{1}{2} \pi \Bigl(I_p(f(x)Y_x, X_xf Y_x), I_p(f(x)Y_x, -X_xf Y_x) \Bigr) +\\
& & \frac{1}{2} \pi \Bigl( m_{f(x)*} (Y_{*_x} X_x), m_{f(x)*} \bigl( \pmb{-} {\mathfrak s}(x) \cdot Y_{*_x}X_x  \bigl) \Bigr) \\
& = & \Bigl. X_xf Y_x + f(x) \nabla_{X_x}Y.
\end{array}$$

We explain now how to associate a symmetry jet ${\mathfrak s}$ to a connection $\nabla$. Extracting ${\mathfrak s}$ from (\ref{sym-jet-->conn2}) yields the following expression~:
\begin{equation}\label{conn-->sym-jet}
{\mathfrak s}(x) \cdot \X = \pmb{-} \X + m_{-1*} \Bigl(I(Y_x, 2 \nabla_{X_x} Y) \Bigr), 
\end{equation}
where $\X = Y_{*_x}X_x$ and $\nabla = \nabla^{\mathfrak s}$. Moreover, for any connection $\nabla$, the relation  (\ref{conn-->sym-jet}) defines a symmetry jet ${\mathfrak s}$ whose associated connection $\nabla^{\mathfrak s}$ is $\nabla$. Indeed, 
$$\begin{array}{ccl}
\nabla^{\mathfrak s}_{X_x}Y & = & \frac{1}{2} \pi \Bigl(Y_{*_x} X_x, \pmb{-} {\mathfrak s}(x) \cdot Y_{*_x}X_x \Bigr) \\
& = & \frac{1}{2} \pi \Bigl(Y_{*_x} X_x, Y_{*_x} X_x - I(Y_x, 2 \nabla_{X_x} Y) \Bigr) \\
& = & \nabla_{X_x} Y.
\end{array}$$
\cqfd

\section{Uniqueness of affine extension revisited}\label{UAE}

This section provides an alternative description of the well-known {\sl property of uniqueness of affine extension}, which states that on a manifold $M$ endowed with a torsionless affine connection $\nabla$, any linear isomorphism $\xi : T_xM \to T_yM$, $x, y \in M$, admits a unique lift to an affine $2$-jet, meaning that there exists a local diffeomorphism $f : U \subset M \to V \subset N$ whose differential at $x$ coincides with $\xi$ and that satisfies for any pair of vector fields $X$ and $Y$
$$f_{*_x} \Bigl( \nabla_{X_x} Y \Bigr) = \nabla_{f_{*_x}(X_x)} \Bigl(f_{*} (Y)\Bigr),$$  
a relation which depends only on the $2$-jet $j^2_xf$ of $f$. The property of uniqueness of affine extension holds for connections with torsion as well, provided the affine extension is allowed to belong to $\b_h^{(1,1)}(M)$ instead of $\b^{(2)}(\PP(M))$.

\begin{dfn}
An affine jet is an element $\xi = j^1_x b$ of $\b_h^{(1,1)}(M)$ such that
$$b(x)(\nabla_{X_x}Y) = \nabla_{b(x)(X_x)}b(Y).$$
\end{dfn}

\begin{rmk}\label{saffine} The image of ${\mathfrak s}$ consists of affine jets. This is a consequence of the fact that any $(1,1)$-jet in $\b^{(1,1)}(\PP(M))$ whose first order part lies in the bisection $-I$ is its own inverse (\pref{e-iota}). Indeed, we know from (\ref{coucou}) that $(s_x(Y))_{*_x} X_x = - {\mathfrak s}(x) \cdot Y_{*_x}X_x$. Thus
$$\begin{array}{ccl}
\nabla^{\mathfrak s}_{X_x}s_x(Y) & = & \frac{1}{2} \pi \Bigl( - {\mathfrak s}(x) \cdot Y_{*_x} X_x, m_{-1*} \bigl( {\mathfrak s}(x) \cdot {\mathfrak s}(x) \cdot Y_{*_x}X_x \bigr)\Bigr) \\
& = & \frac{1}{2} \pi \Bigl(- {\mathfrak s}(x) \cdot Y_{*_x} X_x, m_{-1*} \bigl( Y_{*_x}X_x\bigr) \Bigr) \\
& = & \Bigr. \frac{1}{2}  \nabla^{\mathfrak s}_{X_x}Y.
\end{array}$$
The last equality follows from $\pi(\X, \Y) = - \pi(\Y, \X) = \pi(m_{-1*} \Y, m_{-1*} \X)$.
\end{rmk}

The following proposition consists in the property of uniqueness of affine extension extended to affine connections with torsion and shows also that the map $\b^{(1)}(\PP(M)) \to \b^{(1,1)}(\PP(M))$ that associates to a $1$-jet its unique affine extension, can be characterized as being the unique groupoid morphism and section of $p = p_*$ that extends ${\mathfrak s}$, in the sense that $S \circ -I = {\mathfrak s}$.

\begin{prop}\label{uae} Given an affine connection $\nabla^{\mathfrak s}$ there is a unique Lie groupoid morphism $S : \b^{(1)}(\PP(M)) \to \b_h^{(1,1)}(M)$ (\cite{Mck-05}) such that $S \circ -I = {\mathfrak s}$ and $p \circ S = \id$. Moreover the set of affine jets is precisely the image of $S$. When the symmetry jet is holonomic, the morphism $S$ takes its values in $\b^{(2)}(\PP(M))$.
\end{prop}

\Pf Let $j^1_x b \in \b_h^{(1,1)}(M)$, with $\beta(b(x)) = y$. Then

$$\begin{array}{ccl}
b(x) \bigl( \nabla_{X_x} Y\bigr) - \nabla_{b(x) X_x} b Y 
& = &  b(x) \Bigl(\frac{1}{2} \pi \bigl(Y_{*_x} X_x, \pmb{-} {\mathfrak s}(x) \cdot Y_{*_x}X_x \bigr)\Bigr) \\
& & -  \frac{1}{2} \pi \Bigl((bY)_{*_y} (b(x) X_x), \pmb{-} {\mathfrak s}(y) \cdot (bY)_{*_y} (b(x) X_x) \Bigr)\\
& = & \frac{1}{2} \pi \Bigl(j^1_x b \cdot Y_{*_x} X_x, \pmb{-} j^1_x b \cdot {\mathfrak s}(x) \cdot Y_{*_x}X_x \Bigr) \\
& & - \frac{1}{2} \pi \Bigl(j^1_x b \cdot Y_{*_x} X_x, \pmb{-} {\mathfrak s}(y) \cdot j^1_x b \cdot Y_{*_x}X_x \Bigr) \\
& = & \frac{1}{2} \pi \Bigl( \pmb{-} {\mathfrak s}(y) \cdot j^1_x b \cdot Y_{*_x}X_x, \pmb{-} j^1_x b\cdot {\mathfrak s}(x) \cdot Y_{*_x}X_x \Bigr).
\end{array}$$
This implies that $j^1_x b$ is affine if and only if 
$$j^1_x b \cdot {\mathfrak s}(x) = {\mathfrak s}(y) \cdot j^1_x b.$$
(This statement relies on \lref{111-jetsasmaps}.) Or, equivalently
\begin{equation}\label{equation}
{\mathfrak s}(y) \cdot j^1_x b \cdot {\mathfrak s}(x) = j^1_x b.
\end{equation}
In terms of the associated plane (cf.~\rref{1jets-as-planes} in \aref{jet-of-bis}), the previous equation (\ref{equation}) is satisfied if and only if $D(j^1_xb)$, which lies in $\e_\xi$, (cf.~\rref{h-e}) satisfies 
\begin{equation}\label{equation-bis}
D({\mathfrak s}(y)) \cdot D(j^1_xb) \cdot D({\mathfrak s}(x)) = D(j^1_xb).
\end{equation}

Now, for any $1$-jet $\xi$ in $\b^{(1)}(\PP(M))$, define the map

\begin{equation}\label{psi}
\psi_\xi : T_\xi \b^{(1)}(\PP(M)) \to T_\xi \b^{(1)}(\PP(M)) : X_\xi \mapsto \ol{Y}^{D({\mathfrak s}(y)), \alpha_*} \cdot X_\xi \cdot \ol{X}^{D({\mathfrak s}(x)), \beta_*},
\end{equation}
where $X = \alpha_{*_{\xi}}(X_{\xi})$, $Y = \beta_{*_{\xi}}(X_{\xi})$ and where $\ol{X}^{D({\mathfrak s}(x)), \beta_*}$ (\rp $\ol{Y}^{D({\mathfrak s}(y)), \alpha_*}$) denotes the lift of $X$ (\rp $Y$) via $\beta_*$ (\rp $\alpha_*$) in $D({\mathfrak s}(x))$ (\rp $D({\mathfrak s}(y))$). The dot in the previous formula refers to the differential of the multiplication in the groupoid $\b^{(1)}(\PP(M))$, that is the map 
$$m_{*_{(\xi_2, \xi_1)}} : T_{\xi_2}\b^{(1)}(\PP(M)) \times_{(\alpha_{*_{\xi_2}}, \beta_{*_{\xi_1}})} T_{\xi_1}\b^{(1)}(\PP(M)) \longrightarrow T_{\xi_2 \cdot \xi_1} \b^{(1)}(\PP(M)),$$  
\begin{equation}\label{diff-mult}
m_{*_{(\xi_2, \xi_1)}} (X_{\xi_2}, X_{\xi_1}) \stackrel{\rm not}{=} X_{\xi_2} \cdot X_{\xi_1}.
\end{equation}
    
\begin{figure}[h!] 
\begin{center}
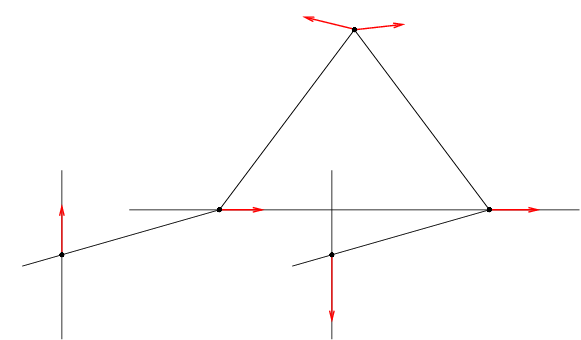 
\caption{The map $\psi_\xi$}
\end{center}
\end{figure} 
The relation (\ref{equation-bis}) can be reformulated in terms of $\psi_\xi$ as follows~:
$$\psi_\xi \Bigl(D(j^1_xb)\Bigl) = D(j^1_xb).$$
The map $\psi_\xi$ is an involutive automorphism of $T_{\xi}\b^{(1)}(\PP(M))$. Indeed, on the one hand, $\alpha_{*_\xi} \circ \psi_\xi = -\alpha_{*_\xi}$ and $\beta_{*_\xi} \circ \psi_\xi = -\beta_{*_\xi}$ and on the other hand,
\begin{equation}\label{x-x=0}
\begin{array}{ccl}
\ol{X}^{D({\mathfrak s}(x)), \beta_*} \cdot \ol{-X}^{D({\mathfrak s}(x)), \beta_*} & = & \ol{X}^{D({\mathfrak s}(x)), \beta_*} \cdot -\ol{X}^{D({\mathfrak s}(x)), \beta_*} \\
& = & \ol{X}^{D({\mathfrak s}(x)), \beta_*} \cdot \iota_{*} \bigl(\ol{X}^{D({\mathfrak s}(x)), \beta_*} \bigr) \\
& = & 0_{x},
\end{array}
\end{equation}
as implied by \pref{e-iota} and, similarly, $\ol{-Y}^{D({\mathfrak s}(y)), \alpha_*} \cdot \ol{Y}^{D({\mathfrak s}(y)), \alpha_*} = 0_y$. Hence $T_\xi \b^{(1)}(\PP(M))$ decomposes into a direct sum of eigenspaces corresponding to the eigenvalues $\pm 1$~:
$$T_\xi \b^{(1)}(\PP(M)) = E^{\psi_\xi}_{+1} \oplus E^{\psi_\xi}_{-1}.$$
Clearly $E^{\psi_\xi}_{+1} = \Ker p_{*_\xi}$ for $p : \b^{(1)}(\PP(M)) \to \PP(M)$. Since $\Ker p_{*_\xi} \subset \e_\xi$, the subspaces $E^{\psi_\xi}_{-1}$ and $\e_\xi$ are transverse and therefore
$$\d_{\xi} = E^{\psi_\xi}_{-1} \cap \e_{\xi}$$
defines a distribution on $\b^{(1)}(\PP(M))$ corresponding to a section 
$$S : \b^{(1)}(\PP(M)) \to \b^{(1,1)}(\PP(M))$$ of $p$
such that $D(S(\xi)) = \d_\xi$. We claim that $S$ is a groupoid morphism whose image consists of the set of affine jets. The first statement is a consequence of the following simple observation~:
$$\psi_{\xi_2 \cdot \xi_1}(X_{\xi_2} \cdot X_{\xi_1}) = \psi_{\xi_2}(X_{\xi_2}) \cdot \psi_{\xi_1}(X_{\xi_1}),$$
itself implied by (\ref{x-x=0}). As to the second statement, along the bisection $-I$, the image of $S$ coincides with ${\mathfrak s}$. Indeed, the relation (\ref{x-x=0}) implies that
$$D(\fs (x)) \subset E^{\psi_{-I_x}}_{-1}.$$ 
Hence $S(-I)$ consists of affine jets (cf.~\rref{saffine}). Thus, for any $\xi : T_xM \to T_yM$ in $\b^{(1)}(\PP(M))$, the property of $S$ to be a groupoid morphism that coincides with ${\mathfrak s}$ on $-I$ implies that $S(\xi) = S(-I_y \cdot \xi \cdot -I_x) = S(-I_y) \cdot S(\xi) \cdot S(-I_x) = {\mathfrak s}(y) \cdot S(\xi) \cdot {\mathfrak s}(x)$ which implies that $S(\xi)$ is affine. \\

Concerning the very last statement of the proposition, it is a consequence of the property that $\kappa$ is a groupoid morphism (cf.~\lref{kappa-prop}) and the first part of the proposition. Indeed, 
$$\kappa (S(\xi)) = \kappa({\mathfrak s}(y) \cdot S(\xi) \cdot {\mathfrak s}(x)) = {\mathfrak s}(y) \cdot \kappa (S(\xi)) \cdot {\mathfrak s}(x)$$ 
implies that $\kappa(S(\xi)) = S(\xi)$.
\cqfd

\begin{dfn}
A section $S : \b^{(1)}(\PP(M)) \to \b^{(1,1)}(\PP(M))$ of $p$ which is a groupoid morphism is called an affine extension or the affine extension of the symmetry jet it is induced from. 
\end{dfn}

\begin{prop}
Affine connections are in bijective correspondence with affine extension.
\end{prop}

Since a $(1,1)$-jet $\xi$ may also appear as a plane $D(\xi)$ tangent to $\b^{(1)}(\PP(M))$ attached to $p(\xi)$ (cf.~\rref{1jets-as-planes} in \aref{jet-of-bis}), the data of the section $S$ is equivalent to that of a distribution, denoted by $\d$ or $\d^{\mathfrak s}$, on $\b^{(1)}(\PP(M))$. It satisfies the following properties~:

\begin{enumerate}
\item[a)] $\d$ is ``transverse" in the sense transverse to both the $\alpha$-fibers and the $\beta$-fibers. 
\item[b)] The fact that $\d$ is induced from a groupoid morphism implies that it coincides with $\varepsilon_{*_x} (T_xM)$ along the identity bisection, is invariant under $\iota_*$ and is preserved by multiplication, that is the map $m_* : TG \times_{(\alpha_*, \beta_*)} TG \to TG$ maps $\d \times_{(\alpha_*, \beta_*)} \d$ onto $\d$.
\item[c)] $\d \subset \e$ (see \dref{def-e}). Equivalently, the bouncing map 
$$\fb : \b^{(1)}(\PP(M)) \to \End(TM, TM) : \xi \mapsto \beta_{*_\xi} \circ \alpha_{*_\xi}|_{\d_\xi}^{-1} = \fb(\d_\xi)$$ induced by $\d$ is the identity.
\end{enumerate}
We use the word ``transverse" to describe condition a) because ``horizontal" is quite dubious in the groupoid case as $\d$ is, at least along the bisection $-I$, rather vertical with respect to our standard picture of a groupoid, since the bouncing map along $-I$ coincides with $-I$. 

\begin{figure}[h!]  
\begin{center}
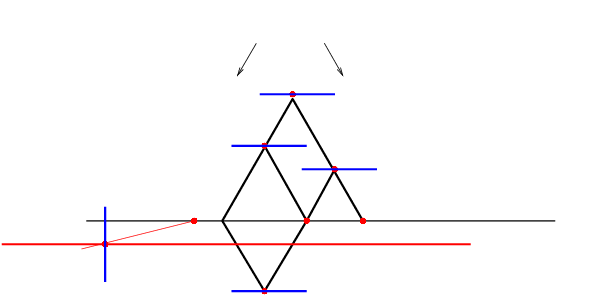 
\caption{The distribution $\d$}
\end{center}
\end{figure}

\begin{dfn}
A distribution on a groupoid $G$ satisfying the previous condition is called an {\sl invariant transverse distribution}.
\end{dfn}

One could rephrase what has been said in this section into the following proposition

\begin{prop}
Symmetry jets are in one-to-one correspondance with affine extensions as well as with invariant transverse distributions.
\end{prop}

Now, affine local or global transformation appear as leaves of the distribution~$\d$.
 
\begin{prop} Let ${\mathfrak s}$ be symmetry jet on the manifold $M$, let $\nabla$ denote the induced affine connection and $\d$ the associated distribution on $\b^{(1)}(\PP(M))$. Through the $1$-jet extension map $\varphi \to j^1\varphi$, affine (local) diffeomorphism of $(M, \nabla^{\mathfrak s})$ correspond to (local) bisections of $\b^{(1)}(\PP(M))$ that are leaves of $\d$.
\end{prop}

\Pf Since $\d$ is contained in $\e$, a leaf of $\d$ is necessarily locally a $1$-jet extension $j^1\varphi$. The latter is then affine since tangent to $\d$.
\cqfd

\section{Torsion and holonomy of affine jets} 

Given a symmetry jet ${\mathfrak s}$ on a manifold $M$, the defect for an affine $(1,1)$-jet $S(\xi)$ to be holonomic, or fixed under the involution $\kappa$ (cf.~\lref{kappa-prop}), coincides with the defect of invariance of the torsion under $\xi$. A precise statement is the content of the next proposition.
\begin{prop}\label{prop-torsion} Let $S(\xi) = j^1_x b$ be an affine jet, then, for any $\X \in T^2M$ with $p(\X) = Y_x$ and $p_*(\X) = X_x$, we have  
\begin{equation}\label{torsion}
\pi \Bigl(S(\xi) \cdot \X, \kappa(S(\xi)) \cdot \X \Bigr) = \xi\Bigl(T^\nabla(X_x, Y_x)\Bigr) - T^\nabla \Bigl(\xi(X_x), \xi(Y_x)\Bigr).
\end{equation}
In particular, the affine $(1,1)$-jet $S(\xi)$ extending $\xi$ is a $2$-jet if and only if $\xi$ preserves the torsion. 
\end{prop}

\Pf Supposing that $X$ and $Y$ are vector fields on $M$ extending $X_x$ and $Y_x$ respectively and such that $\X = \kappa (X_{*_x}Y_x)$, the right hand side of (\ref{torsion}) equals

$$\begin{array}{cl}
& \displaystyle{\xi \Bigl(\nabla_{X_x}Y - \nabla_{Y_x}X - [X,Y]_x\Bigr)} \\
& \displaystyle{- \Bigl( \nabla_{\xi(X_x)} bY + \nabla_{\xi(Y_x)} bX + [bX, bY]_y \Bigr)} \\
= & \displaystyle{\Bigl(\xi(\nabla_{X_x}Y) - \nabla_{\xi(X_x)} bY\Bigr) - \Bigl(\xi(\nabla_{Y_x}X) - \nabla_{\xi(Y_x)} bX\Bigr)} \\
& \displaystyle{- \Bigl(\xi([X,Y]_x) - [bX, bY]_y\Bigr)} \\
= & \displaystyle{- \; \xi \circ \pi \Bigr( Y_{*_x}X_x, \kappa (X_{*_x}Y_x)\Bigr) + \pi \Bigl((bY)_{*_y} \xi(X_x), \kappa\bigl((bX)_{*_y} \xi(Y_x)\bigr)\Bigr).} \\
= & \displaystyle{- \; \pi \Bigr(S(\xi) \cdot Y_{*_x}X_x, S(\xi) \cdot \kappa (X_{*_x}Y_x)\Bigr)} \\
& \displaystyle{+ \; \pi \Bigr(S(\xi) \cdot Y_{*_x}X_x, \kappa(S(\xi) \cdot X_{*_x}Y_x) \Bigr)}\\
= & \displaystyle{\pi \Bigl(S(\xi) \cdot \kappa(X_{*_x}Y_x), \kappa(S(\xi)) \cdot \kappa(X_{*_x}Y_x) \Bigr)} \\
= & \displaystyle{\pi \Bigl(S(\xi) \cdot \X, \kappa(S(\xi)) \cdot \X \Bigr),}
\end{array}$$
where we have used \pref{involution-prop}, as well as relation (\ref{coucou}) from \aref{b11}.
\cqfd
Equation (\ref{torsion}) yields a geometric interpretation of the torsion of an affine connection in terms of its symmetry jet.

\begin{cor} Let ${\mathfrak s}$ be a symmetry jet, and let $\nabla$ be the corresponding affine connection. Then
$$T^{\nabla} (X_x, Y_x) = \frac{1}{2} \pi\Bigr(\kappa\bigl({\mathfrak s}(x)\bigr) \cdot \X, {\mathfrak s}(x) \cdot \X \Bigr),$$
for any $\X \in T^2M$ with $p(\X) = Y_x$ and $p_*(\X) = X_x$.
\end{cor}

In terms of the distribution $\d^\fs$, this means that for any $X_x \in T_xM$, the endomorphism $T^\nabla(X_x, \cdot)$ of $T_xM$ is the difference between the lifts $X_1$ and $X_2$ of $X_x$ in $\d_{-I_x}$ and $\kappa(\d_{-I_x})$ respectively with respect to $\alpha_*$. Indeed, the vector $X_2 - X_1$ lies in  $T_{-I_x}(\b^{(1)}(\PP(M))_{x,x})$ which is naturally identified with $\End(T_xM, T_xM)$. 

\begin{figure}[h!]  
\begin{center}
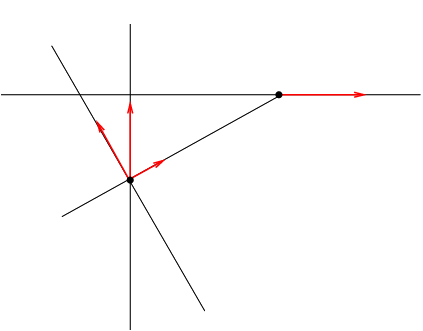 
\caption{The torsion}
\end{center}
\end{figure}

\section{The curvature in terms of the symmetry jet} 

In this section we present a very simple expression for the curvature of an affine connection in terms of the first jet of the symmetry jet. More precisely, we prove the following statement

\begin{thm}\label{curvthm} Let ${\mathfrak s}$ be a symmetry jet on the manifold $M$. Then the curvature tensor $R$ of the associated affine connection admits the following expression 
\begin{equation}\label{curv-j1s}
R(X_x, Y_x) Z_x = \frac{1}{4} \Pi \Bigl(\kappa(j^1_x {\mathfrak s}) \cdot j^1_x {\mathfrak s} \cdot \IX, \;j^1_x {\mathfrak s} \cdot \kappa(j^1_x {\mathfrak s}) \cdot \IX \Bigr),
\end{equation}
where $X$, $Y$ and $Z$ are vector fields on $M$ and $\IX$ stands for $Z_{**_{Y_x}} Y_{*_x} X_x \in T^3M$. 
\end{thm}

\begin{rmk} One could also write, with a slight abuse of notation
$$R = \frac{1}{4}\Pi([\kappa(j^1_x{\mathfrak s}), j^1_x{\mathfrak s}]).$$
\end{rmk}

\begin{rmk} Observe that $\kappa(j^1_x {\mathfrak s})$ is not a $(1,1,1)$-jet since $p(j^1_x {\mathfrak s}) = {\mathfrak s}(x)$ does not coincide with $p_*(j^1_x {\mathfrak s}) = m_{-1*}$ (cf.~\rref{kappa-rem-bis}). Nevertheless, $\kappa(j^1_x {\mathfrak s})$ is an element in $\EL(T^3M)$ (cf.~\dref{el-def}) whose action on an element $\IX$ in $T^3M$ is defined by
$$\kappa(j^1_x {\mathfrak s}) \cdot \IX = \kappa(j^1_x {\mathfrak s} \cdot \kappa(\IX)).$$
Moreover,
\begin{enumerate}
\item[-] $p(\kappa(j^1_x {\mathfrak s})) = m_{-1*}$,
\item[-] $p_*(\kappa(j^1_x {\mathfrak s})) = {\mathfrak s}(x)$,
\item[-] $p_{**}(\kappa(j^1_x {\mathfrak s})) = m_{-1}$.
\end{enumerate}
\end{rmk}

\Pf As a first step, let us compute $\nabla_{X_x} \nabla_Y Z$ in terms of the symmetry jet. 
$$\nabla_{X_x} \nabla_Y Z = \displaystyle{\frac{1}{2} \pi \Bigl((\nabla_Y Z)_{*_x} (X_x), \pmb{-} {\mathfrak s}(x) \cdot (\nabla_Y Z)_{*_x} (X_x)\Bigr)}$$
Let $\gamma : (-\eps, \eps) \to M$ be a path in $M$ tangent to $X_x$ at $t = 0$. Then

$$\begin{array}{ccl}
(\nabla_Y Z)_{*_x} (X_x) & = &\displaystyle{\frac{d}{dt} (\nabla_Y Z)_{\gamma(t)} \Bigl|_{t=0}} \\
& = & \displaystyle{\frac{d}{dt} \frac{1}{2} \pi \Bigl( Z_{*_{\gamma(t)}}Y_{\gamma(t)}, \pmb{-} {\mathfrak s}(\gamma(t)) \cdot Z_{*_{\gamma(t)}}Y_{\gamma(t)} \Bigl)\Bigl|_{t=0}} \\
& = & m_{\frac{1}{2}*} \Bigl(\pi_{*_{(\ZE, \pmb{-} {\mathfrak s}(x) \cdot \ZE)}} \Bigl( \IX, m_{-1*} \circ m_{-1**} \bigl( j^1_x{\mathfrak s} \cdot \IX\bigr) \Bigr)\Bigr),
\end{array}$$
where $\ZE = p(\IX) = Z_{*_x}Y_x$. Since $\IX^1 = \IX$ and $\IX^2 = m_{-1*} \circ m_{-1**} ( j^1_x{\mathfrak s} \cdot \IX)$ belong to the same $({\mathcal P}_2 = p_* \times p_{**})$-fiber, there exists a $\U \in T^{X_x}TM$ such that
$$\IX^1 = \IX^2 +_* \Bigl(e_*(\IX^2) +_{**} (i^p_{0_M})_*(\U)\Bigr) = \IX^2 +_{**} \Bigl(e_{**}(\IX^2) +_{*} (i^p_{0_M})_*(\U)\Bigr).$$
(This follows from (\ref{proj123}) in \aref{third-der}). It is not difficult to verify that $\pi_*(\IX^1, \IX^2) = \Pi_2(\IX^1, \IX^2) = \U$. Now, we claim that for an element $\xi \in \EL^{(1,1,1)}(T^3M)$, 
\begin{equation}\label{pi*}
\pi_*(\xi \cdot \IX^1, \xi \cdot \IX^2) = p_*(\xi) \cdot \pi_*(\IX^1, \IX^2).
\end{equation}
This is verified as follows~:
$$\begin{array}{ccl}
\xi \cdot \IX^1 & = & \xi \cdot \Bigl(\IX^2 +_* \bigl(e_*(\IX^2) +_{**} (i^p_{0_M})_*(\U)\bigr)\Bigr) \\
& = & \xi \cdot \IX^2 +_* \bigl(\xi \cdot e_*(\IX^2) +_{**} \xi \cdot (i^p_{0_M})_*(\U)\bigr) \\
& = & \xi \cdot \IX^2 +_* \bigl(e_*(\xi \cdot \IX^2) +_{**} (i^p_{0_M})_*(p_*(\xi) \cdot \U)\bigr) 
\end{array}$$
The last equality follows from \lref{111-jetsasmaps} and holds for any $\xi \in \EL(T^3M)$ such that $\xi \cdot (i^p_{0_M})_*(\U) = (i^p_{0_M})_*(p_*(\xi) \cdot \U)$. In particular for $\xi = \kappa(j^1_x \fs)$. Indeed, 
$$\begin{array}{cll}
\kappa(j^1_x \fs) \cdot \Bigl( (i^p_{0_M})_*(\U) \Bigr) & = & \kappa\Bigl( j^1_x \fs \cdot \kappa\bigl( (i^p_{0_M})_*(\U) \bigr)\Bigr) \\
& = & \kappa\Bigl(j^1_x \fs \cdot i^{p_*}_{{0_*}_{TM}}(\U) \Bigr) \\
& = & \kappa \Bigl( i^{p_*}_{{0_*}_{TM}} \bigl( \fs(x) \cdot \U \bigr)\Bigr) \\
& = & (i^p_{0_M})_* \Bigl(\fs(x) \cdot \U \Bigr).
\end{array}$$
Thus 
$$\nabla_{X_x} \nabla_Y Z = \displaystyle{\frac{1}{2} \pi \Bigl( m_{\frac{1}{2}*} \bigl[\pi_{*} (\IX^1, \IX^2) \bigr], \pmb{-} m_{\frac{1}{2}*} \bigl[\pi_{*} ( \kappa(j^1_x \fs) \cdot \IX^1, \kappa(j^1_x \fs) \cdot \IX^2 ) \bigr] \Bigr)}.$$
To pursue this computation, observe that $\pi((m_a)_*\X, (m_a)_*\Y) = a \pi (\X, \Y)$ for two elements $\X$ and $\Y$ in a same $(p \times p_*)$-fiber and a real $a$. Thus  
$$\nabla_{X_x} \nabla_Y Z = \displaystyle{\frac{1}{4} \pi \Bigl( \pi_{*} (\IX^1, \IX^2), \pmb{-} \; \pi_{*} (\kappa(j^1_x \fs) \cdot \IX^1, \kappa(j^1_x \fs) \cdot\IX^2) \Bigr)}.$$
Moreover $m_{-1} (\pi_{*} (\IY^1,\IY^2)) = \pi_{*} (m_{-1} (\IY^1), m_{-1} (\IY^2))$ and $m_{-1*} ( \pi_{*} (\IY^1,\IY^2))$ coincides with both $\pi_{*} (m_{-1*} (\IY^1), m_{-1*} (\IY^2))$ and $\pi_{*} (m_{-1**} (\IY^1), m_{-1**} (\IY^2))$ for any $\IY^1$, $\IY^2$ in $T^3M$. Thus 
\begin{equation}\label{nablacarre}
\begin{array}{ccl}
\nabla_{X_x} \nabla_Y Z & = & \displaystyle{\frac{1}{4} \pi \Bigl( \pi_{*} \bigl(\IX^1, \IX^2 \bigr), }\\
& & \displaystyle{\pi_{*} \bigl( m_{-1} \circ m_{-1**} ( \kappa(j^1_x {\mathfrak s}) \cdot \IX^1 ), m_{-1} \circ m_{-1**} ( \kappa(j^1_x {\mathfrak s}) \cdot \IX^2) \bigr) \Bigr)}
\end{array}
\end{equation}
Let us define 
\begin{enumerate}
\item[-] $\widetilde{\IX^1} = m_{-1} \circ m_{-1**} ( \kappa(j^1_x {\mathfrak s}) \cdot \IX^1)$, 
\item[-] $\widetilde{\IX^2} = m_{-1} \circ m_{-1**} (\kappa(j^1_x {\mathfrak s}) \cdot \IX^2) = m_{-1} \circ m_{-1*} (\kappa(j^1_x {\mathfrak s}) \cdot j^1_x{\mathfrak s} \cdot \IX)$. 
\end{enumerate}
Now $\U = \pi_{*} (\IX^1, \IX^2)$, $\widetilde{\U} = \pi_{*} (\widetilde{\IX^1}, \widetilde{\IX^2})$  and $U = \pi(\pi_{*} (\IX^1, \IX^2), \pi_{*} (\widetilde{\IX^1}, \widetilde{\IX^2})) = \nabla_{X_x} \nabla_Y Z$ satisfy the relations
$$\begin{array}{ccl}
\IX^1 & = & \IX^2 +_{**} \Bigl(e_{**}(\IX^2) +_* (i^p_{0_M})_*(\U)\Bigr)\\
\widetilde{\IX^1} & = & \widetilde{\IX^2} +_{**} \Bigl(e_{**}(\widetilde{\IX^2}) +_* (i^p_{0_M})_*(\widetilde{\U})\Bigr) \\
\U & = & \widetilde{\U} +_* \bigl( e_*(\U) + i^p_{0_M}(U) \bigr).
\end{array}$$
Besides, all four elements $\IX^1$, $\IX^2$, $\widetilde{\IX^1}$ and $\widetilde{\IX^2}$ belong to the same $p_{**}$-fiber and 
$$\begin{array}{ccl}
&&\IX^1 -_{**} \IX^2 -_{**} \widetilde{\IX^1} +_{**} \widetilde{\IX^2} \\
& = & \Bigl(e_{**}(\IX^2) +_{*} (i^p_{0_M})_* (\U) \Bigl) -_{**} \Bigl( e_{**}(\widetilde{\IX^2}) +_{*} (i^p_{0_M})_*(\widetilde{\U})\Bigr) \\
& = & e_{**}(\IX^2) +_* \Bigl( (i^p_{0_M})_* (\U) -_* (i^p_{0_M})_*(\widetilde{\U})\Bigr) \\
& = & e_{**}(\IX^2) +_* \Bigl( (i^p_{0_M})_* (\U -_*\widetilde{\U})\Bigr) \\
& = & e_{**}(\IX^2) +_* \Bigl( (i^p_{0_M})_* \bigl(e_*(\U) + i^p_{0_M} (U)\bigr)\Bigr)\\
& = & e_{**}(\IX^2) +_* \Bigl( (i^p_{0_M})_* \circ i_* \circ p_* (\U) + (i^p_{0_M})_* \circ i^p_{0_M} (U)\Bigr)\\
& = & e_{**}(\IX^2) +_* \Bigl( i_* \circ i_* \circ p_* \circ p_{**} (\IX^2) + I (U)\Bigr)\\
& = & e_{**}(\IX^2) +_* \Bigl( i_* \circ p_* \circ i_{**} \circ p_{**} (\IX^2) + I (U)\Bigr)\\
& = & e_{**}(\IX^2) +_* \Bigl( e_* (e_{**} (\IX^2)) + I (U)\Bigr).\\
\end{array}$$
This computation shows, in terms of the piece of notation introduced in (\ref{pi-bis}) that 
$$\nabla_{X_x} \nabla_Y Z = \frac{1}{4} \Pi \Bigl( \IX^1 -_{**} \IX^2 -_{**} \widetilde{\IX^1} +_{**} \widetilde{\IX^2} \Bigr).$$
In other terms
$$\begin{array}{c}
\nabla_{X_x} \nabla_Y Z = \\
\frac{1}{4} \Pi \Bigl( \IX +_{**} \; (m_{-1})_* j^1_x{\mathfrak s} \cdot \IX +_{**} \; (m_{-1}) \; \kappa(j^1_x{\mathfrak s}) \cdot \IX +_{**} \; (m_{-1}) \circ (m_{-1})_{*} \; \kappa(j^1_x{\mathfrak s}) \cdot j^1_x{\mathfrak s} \cdot \IX \Bigr).
\end{array}$$

Now we can tackle the curvature. Without loss of generality, we may assume that $[X,Y]_x = 0$. Then
$$\begin{array}{ccl}
R(X_x, Y_x) Z_x & = & \nabla_{X_x} \nabla_Y Z - \nabla_{Y_x} \nabla_X Z \\
& = & \frac{1}{4} \Pi \Bigl( \IX^1 -_{**} \IX^2 -_{**} \widetilde{\IX^1} +_{**} \widetilde{\IX^2} \Bigr) - \\
& & \frac{1}{4} \Pi \Bigl( \IY^1 -_{**} \IY^2 -_{**} \widetilde{\IY^1} +_{**} \widetilde{\IY^2} \Bigr),
\end{array}$$
with 
\begin{enumerate}
\item[-] $\IY^1 = \IY = Z_{**_{X_x}} X_{*_x} Y_x$, 
\item[-] $\IY^2 = (m_{-1})_* \circ (m_{-1})_{**} j^1_x{\mathfrak s} \cdot \IY$, 
\item[-] $\widetilde{\IY^1} = (m_{-1}) \circ (m_{-1})_{**} \; \kappa(j^1_x{\mathfrak s}) \cdot \IY$, 
\item[-] $\widetilde{\IY^2} = (m_{-1}) \circ (m_{-1})_{*} \; \kappa(j^1_x{\mathfrak s}) \cdot j^1_x{\mathfrak s} \cdot \IY$. 
\end{enumerate}
Observe that (\ref{i-kappa-bis}) implies that 
$$\begin{array}{ccl}
\Pi \Bigl( \IY^1 -_{**} \IY^2 -_{**} \widetilde{\IY^1} +_{**} \widetilde{\IY^2} \Bigr) & = & \Pi \circ \kappa \Bigl( \IY^1 -_{**} \IY^2 -_{**} \widetilde{\IY^1} +_{**} \widetilde{\IY^2} \Bigr)\\
& = & \Pi \Bigl( \kappa(\IY^1) -_{**} \kappa(\IY^2) -_{**} \kappa(\widetilde{\IY^1}) +_{**} \kappa(\widetilde{\IY^2}) \Bigr).
\end{array}$$
Moreover, 
\begin{enumerate}
\item[-] $\kappa(\IY^1) = \kappa (Z_{**_{X_x}} X_{*_x} Y_x) = Z_{**_{X_x}} \kappa(X_{*_x} Y_x) = Z_{**_{X_x}} Y_{*_x} X_x = \IX^1$, 
\item[-] $\kappa(\IY^2) = (m_{-1}) \circ (m_{-1})_{**} \kappa(j^1_x{\mathfrak s}) \cdot \IX = \widetilde{\IX^1}$, 
\item[-] $\kappa(\widetilde{\IY^1}) = (m_{-1})_* \circ (m_{-1})_{**} \; j^1_x{\mathfrak s} \cdot \IX = \IX^2$, 
\item[-] $\kappa(\widetilde{\IY^2}) = (m_{-1})_* \circ (m_{-1}) \; j^1_x{\mathfrak s} \cdot \kappa(j^1_x{\mathfrak s}) \cdot \IX \neq \widetilde{\IX^2}$. 
\end{enumerate}
Therefore, 
$$\begin{array}{ccl}
R(X_x, Y_x) Z_x & = & \nabla_{X_x} \nabla_Y Z - \nabla_{Y_x} \nabla_X Z \\
& = & \frac{1}{4} \Pi \Bigl( \IX^1 -_{**} \IX^2 -_{**} \widetilde{\IX^1} +_{**} \widetilde{\IX^2} \Bigr) - \\
& & \frac{1}{4} \Pi \Bigl( \IX^1 -_{**} \widetilde{\IX^1} -_{**} \IX^2 +_{**} \kappa(\widetilde{\IY^2}) \Bigr) \\
& = &  \frac{1}{4} \Pi \Bigl( \bigl(\IX^1 -_{**} \IX^2 -_{**} \widetilde{\IX^1} +_{**} \widetilde{\IX^2}\bigr), \\
& & \;\;\;\;\;\;\; \bigl(\IX^1 -_{**} \widetilde{\IX^1} -_{**} \IX^2 +_{**} \kappa(\widetilde{\IY^2}) \bigr) \Bigr) \\
& = & \frac{1}{4} \Pi \Bigl(\widetilde{\IX^2}, \kappa(\widetilde{\IY^2}) \Bigr)
\end{array}$$
In other words,
$$\begin{array}{ccl}
R(X_x, Y_x) Z_x & = & \frac{1}{4} \Pi \Bigl((m_{-1}) \circ (m_{-1})_{*} \; \kappa(j^1_x{\mathfrak s}) \cdot j^1_x{\mathfrak s} \cdot \IX, \\
& & \;\;\;\;\;\;\; (m_{-1})_* \circ (m_{-1}) \; j^1_x{\mathfrak s} \cdot \kappa(j^1_x{\mathfrak s}) \cdot \IX \Bigr) \\
& = & \frac{1}{4} \Pi \Bigl(\kappa(j^1_x{\mathfrak s}) \cdot j^1_x{\mathfrak s} \cdot \IX, j^1_x{\mathfrak s} \cdot \kappa(j^1_x{\mathfrak s}) \cdot \IX \Bigr)
\end{array}$$
\cqfd

\section{Third order affine extension}\label{third-order-affine-extension}

As for the order $2$, on can prove that affine jets of order $3$ extending a given $1$-jet always exist, provided $(1,1,1)$-jets, that is elements of the groupoid $\b^{(1,1,1)}_{nh}(\PP(M))$ are allowed. 
Recall from \aref{b111} that the later are of the type
$$\xi = j^1_xb,$$ 
where 
$$b : U_x \to \b^{(1,1)}_{nh}(\PP(M)) : x' \mapsto j^1_{x'}b_{x'}$$ 
is some local bisection of $\b^{(1,1)}_{nh}(\PP(M))$ and the various
$$b_{x'} : U_{x'} \to \b^{(1)}(\PP(M)),$$
for $x' \in U_x$, form a smooth family of local bisections of $\b^{(1)}(\PP(M))$. Recall also that when $\xi$ lies in $\b^{(1,1,1)}(\PP(M))$, we may assume that $b_{x'}(x') = b_x(x')$ and that $b_{x'}$ is tangent to $\e$ for all $x'$ or equivalently that $(\beta \circ b_{x'})_{*_{x'}} = b_{x'}(x')$ (cf.~observation following \dref{b111h}).

\begin{dfn}\label{affine111jet} A $(1,1,1)$-jet $\xi$ is affine if it belongs to $\b^{(1,1,1)}(\PP(M))$, if its $(1,1)$-part $p(\xi) = p_*(\xi) = p_{**}(\xi)$ is affine and if for any vector fields $X$, $Y$, $Z$ in ${\mathfrak X}(M)$, we have 
\begin{equation}\label{affine111}
\xi \Bigl(\nabla_{X_x}\nabla_Y Z \Bigr) = \nabla_{\xi X_x} \nabla_{b_xY} b_{\centerdot} Z.
\end{equation}
\end{dfn}
Let us say a few words about the right hand side of (\ref{affine111}). The vector field $b_xY$ is defined on a neighborhood $V_y$ of $y$ by $(b_xY)_{y'} = b_x(x') Y_{x'}$ with $b_x^0(x') = y'$. The notation $b_{\centerdot} Z$ stands for the family $T_{y'}$ of vector fields parameterized by $y' = b_x^0(x') \in V_y$~:
$$T_{y'} : V_{y'} \to TM : y'' = b_{x'}^0(x'') \mapsto (b_{x'}Z)_{y''} = b_{x'}(x'')Z_{x''}.$$ 
It is differentiated covariantly in the direction of the vector field $y' \mapsto (b_x Y)_{y'}$ and the result, that depends twofold on the variable $y'$  is being covariantly differentiated in the direction of $\xi X_x \in T_yM$. 

\begin{rmk} It is also important to notice that a $(1,1,1)$-jet that satisfies (\ref{affine111}) alone does not necessarily have an affine $(1,1)$-part. Indeed, let $b$ denote a local bisection of $\b^{(1)}(\PP(M))$ such that $j^1_xb = S(\xi)$. Then 
$$\xi \Bigl(\nabla_{X_x}\nabla_Y Z \Bigr) = \nabla_{\xi X_x}b \bigl(\nabla_Y Z \bigr).$$
So $j^1_x j^1_\centerdot b_\centerdot$ is affine if and only if 
$$\nabla_{\xi X_x}\Bigl(b \bigl(\nabla_Y Z \bigr) - \nabla_{b_x Y} b_{\centerdot} Z \Bigr) = 0$$
for all $X$, $Y$, $Z$ in $\IX(M)$. The latter relation only means that, for any vector fields $Y$ and $Z$, the two local vector fields 
$$U = b \bigl(\nabla_Y Z \bigr) \quad \mbox{and} \quad V = \nabla_{b_x Y} b_{\centerdot} Z$$ induce the same map 
$p^v \circ U_{*_y} = p^v \circ V_{*_y} : T_yM \to T_yM$, where $p^v$ denotes the projection $p^v : T^2M \to TM$ induced from the horizontal distribution on $T^2M$ associated to $\nabla^{\mathfrak s}$. Still, the vectors $U_y$ and $V_y$ might not agree in general. Equivalently, $j^1_x b_x$ might not coincide with $j^1_x b = S(\xi)$. 
\end{rmk}

The following statement follows directly from \pref{uae}.
\begin{prop}\label{uaebis} Given a symmetry jet ${\mathfrak s} : M \to \b^{(1,1)}(\PP(M))$ and the corresponding distribution $\d$ on $\b^{(1)}(\PP(M))$, the (tautological) distribution 
$${\mathfrak D}_{S(\xi)} = S_{*_\xi}(\d_\xi)$$ 
along $\im S \subset \b^{(1,1)}(\PP(M))$ corresponds to a groupoid morphism 
$${\mathbb S} : \b^{(1)}(\PP(M)) \to \b^{(1,1,1)}(\PP(M))$$ 
whose image consists of affine $(1,1,1)$-jets and such that $p \circ {\mathbb S}$ coincides with $S$.
\end{prop} 

\Pf For any $\xi \in \b^{(1)}(\PP(M))$ let $a_\xi$ denote some local bisection $x_1 \mapsto a_\xi(x_1)$ such that $j^1_xa_\xi = S(\xi)$, for $x = \alpha(\xi)$. Then, for each point $a_\xi(x_1)$, the expression $a_{a_\xi(x_1)}$ denotes a local bisection, that depends smoothly on $x_1$, whose first jet at $x_1$ coincides with $S(a_\xi(x_1))$. It is tautological that the $(1,1,1)$-jet  
$$j^1_x (j^1_{x_1} a_{a_\xi(x_1)}) = j^1_x (S \circ a_\xi)$$ 
is affine. Indeed,
$$\xi \Bigl( \nabla_{X_x} \nabla_Y Z \Bigr) = \nabla_{\xi X_x} a_\xi \Bigl( \nabla_{Y}Z \Bigr) = \nabla_{\xi X_x} \nabla_{a_\xi Y} a_{a_\xi(\centerdot)} Z.$$
Moreover, $j^1_x (S \circ a_\xi)$ belongs to $\b^{(1,1,1)}(\PP(M))$~:
\begin{enumerate}
\item[-] $p(j^1_x (S \circ a_\xi)) = S \circ a_\xi(x) = S(\xi)$,
\item[-] $p_*(j^1_x (S \circ a_\xi)) = j^1_x(p \circ S \circ a_\xi) = j^1_x a_\xi = S(\xi)$,
\item[-] $p_{**}(j^1_x (S \circ a_\xi)) = j^1_x(p_* \circ S \circ a_\xi) = j^1_x a_\xi = S(\xi)$.
\end{enumerate}
\ \\
Observe now that $D(j^1_x (S \circ a_\xi)) = (S \circ a_\xi)_{*_x}(T_xM) = S_{*_\xi}(\d_\xi)$. Furthermore, the section 
$${\mathbb S} : \b^{(1)}(\PP(M)) \to \b^{(1,1,1)}(\PP(M)) : \xi \mapsto {\mathbb S}(\xi) = j^1_x(S \circ a_\xi)$$ 
is a groupoid morphism. Indeed, if $(\xi_1, \xi_2) \in \b^{(1)} \times_{(\alpha, \beta)} \b^{(1)}$, then
$$D\Bigl({\mathbb S}(\xi_1 \cdot \xi_2)\Bigr) = S_{*_{\xi_1 \cdot \xi_2}}(\d_{\xi_1 \cdot \xi_2}) = S_{*_{\xi_1 \cdot \xi_2}}(\d_{\xi_1} \cdot \d_{\xi_2}) = S_{*_{\xi_1}}(\d_{\xi_1} ) \cdot S_{*_{\xi_2}} (\d_{\xi_2})$$
(cf.~\rref{1jets-as-planes}) implies that ${\mathbb S}(\xi_1 \cdot \xi_2) = {\mathbb S}(\xi_1) \cdot {\mathbb S}(\xi_2)$. 
\cqfd

\begin{prop}\label{uae-order3} The affine $(1,1,1)$-jet ${\mathbb S}(\xi)$ is the unique affine $(1,1,1)$-jet whose first order is $\xi$.
\end{prop}

\Pf The idea of the proof is to compute 
\begin{equation}\label{xi-nabla}
\xi \Bigl(\nabla_{X_x} \nabla_Y Z \Bigr) - \nabla_{\xi X_x} \nabla_{b_x Y} b_{\centerdot} Z,
\end{equation}
for a $(1,1,1)$-jet $\zeta = j^1_xj^1_\centerdot b_\centerdot$ in $\b^{(1,1,1)}(\PP(M))$ whose second order part is $S(\xi)$, in terms of ${\mathfrak s}$ and ${\mathbb S}$ so as to make the condition that (\ref{xi-nabla}) vanishes equivalent to $\zeta = {\mathbb S}(\xi)$. Using the formula (\ref{nablacarre}) with $\xi_o = {\mathbb S}(-I_x)$ and taking into account the facts that $p(\xi)(\pi(\X^1, \X^2)) = \pi(\xi \cdot \X^1, \xi \cdot \X^2)$ for any $(1,1)$-jet $\xi$ and $p_*(\xi) \cdot \pi_*(\IX^1, \IX^2) = \pi_* (\xi \cdot \IX^1, \xi \cdot \IX^2)$ for any $(1,1,1)$-jet $\xi$ (cf.~(\ref{pi*})), the first term $\xi(\nabla_{X_x} \nabla_Y Z)$ can be rewritten~:
$$\begin{array}{l}
\quad \displaystyle{\frac{1}{4} \xi \Bigl\{ \pi \Bigl[ \pi_{*} \Bigl(\IX, (m_{-1})_* \circ (m_{-1})_{**} j^1_x{\mathfrak s} \cdot \IX \Bigr),} \\
\quad \displaystyle{\pi_{*} \Bigl( (m_{-1}) \circ (m_{-1})_{**} \; {\mathbb S}(-I_x) \cdot \IX, (m_{-1}) \circ (m_{-1})_{*} \; {\mathbb S}(-I_x) \cdot j^1_x{\mathfrak s} \cdot \IX \Bigr) \Bigr]\Bigr\}} \\  
= \displaystyle{\frac{1}{4} \pi \Bigl[ \pi_{*} \Bigl({\mathbb S}(\xi) \cdot \IX, (m_{-1})_* \circ (m_{-1})_{**} {\mathbb S}(\xi) \cdot j^1_x{\mathfrak s} \cdot \IX \Bigr),} \\
\quad \displaystyle{\pi_{*} \Bigl((m_{-1}) \circ (m_{-1})_{**} {\mathbb S}(\xi) \cdot {\mathbb S}(-I_x) \cdot \IX, (m_{-1}) \circ (m_{-1})_{*} {\mathbb S}(\xi) \cdot  {\mathbb S}(-I_x) \cdot j^1_x{\mathfrak s} \cdot \IX \Bigr) \Bigr]} \\
= \displaystyle{\frac{1}{4} \pi \Bigl[ \pi_{*} \Bigl({\mathbb S}(\xi) \cdot \IX, (m_{-1})_* \circ (m_{-1})_{**} j^1_y{\mathfrak s} \cdot {\mathbb S}(\xi) \cdot \IX \Bigr),} \\
\quad \displaystyle{\pi_{*} \Bigl( (m_{-1}) \circ (m_{-1})_{**} {\mathbb S}(-I_y) \cdot {\mathbb S}(\xi) \cdot \IX, (m_{-1}) \circ (m_{-1})_{*} {\mathbb S}(-I_y) \cdot j^1_y{\mathfrak s} \cdot {\mathbb S}(\xi) \cdot \IX \Bigr) \Bigr],} 
\end{array}$$
where $\IX = Z_{**_{Y_x}Y_{*_x}X_x}$. The third equality follows from the fact that ${\mathbb S}(\xi)$ commutes with ${\mathbb S}(-I_x)$ and with $j^1_x{\mathfrak s}$. This is really the key point here and the main property of ${\mathbb S}(\xi)$ that distinguishes it from other $(1,1,1)$-jets. Indeed, ${\mathbb S}$ being a groupoid morphism, we have 
$${\mathbb S}(\xi) \cdot {\mathbb S}(-I_x) = {\mathbb S}(\xi \cdot -I_x) = {\mathbb S}(-I_y \cdot \xi) = {\mathbb S}(-I_y) \cdot {\mathbb S}(\xi)$$ 
and
$$\begin{array}{lll}
{\mathbb S}(\xi) \cdot j^1_x\fs & = & j^1_x (S \circ a_\xi) \cdot j^1_x(S \circ -I) = j^1_x \Bigl((S \circ a_\xi) \cdot (S \circ -I)\Bigr) \\
& = & j^1_x\Bigl(S \circ (a_\xi \cdot -I)\Bigr) = j^1_x\Bigl(S \circ (-I \cdot a_\xi)\Bigr) = j^1_x\Bigl((S \circ -I) \cdot (S \circ a_\xi)\Bigr) \\
& = & \Bigl.j^1_y\fs \cdot {\mathbb S}(\xi).
\end{array}$$

For the second term of (\ref{xi-nabla}), notice first that since $\zeta$ belongs to $\b^{(1,1,1)}(\PP(M))$, we may assume that $b_{x'}(x') = b_{x}(x')$ and $(\beta \circ b_{x'})_{*_{x'}} = b_{x'}(x')$. Let $\gamma : (-\eps, \eps) \to M$ be a path tangent to $X_x$ at $0$. Then $\tilde{\gamma} = \beta \circ b_x \circ \gamma$ is tangent to $\xi X_x$ at $y$ and observe that 
$$\nabla_{\xi X_x} \nabla_{b_x Y} b_{\centerdot} Z = \displaystyle{\frac{1}{2} \pi \Bigl( \frac{d}{dt} \bigl(\nabla_{b_x Y} b_{\centerdot} Z\bigr)_{\tilde{\gamma}(t)} \Bigr|_{t=0}, \pmb{-} {\mathfrak s}(x) \cdot \frac{d}{dt} \bigl(\nabla_{b_x Y} b_{\centerdot} Z\bigr)_{\tilde{\gamma}(t)} \Bigr|_{t=0} \Bigr)}.$$
Now, for each $t \in (- \eps, \eps)$ let $\tau_t : (-\eps, \eps) \to M$ be a path tangent to $Y_{\gamma(t)}$ at $0$. Again, the path $\tilde{\tau}(t) = \beta \circ b_{\gamma(t)} \circ \tau_t$ is tangent to $(b_xY)_{\tilde{\gamma}(t)}$.
Then
$$\begin{array}{cll}
\Bigl(\nabla_{b_x Y} b_{\centerdot} Z\Bigr)_{\tilde{\gamma}(t)} & = & \displaystyle{ \nabla_{(b_x Y)_{\tilde{\gamma}(t)}} \bigl(b_{\gamma(t)} Z\bigr)} \\
& = & \displaystyle{\frac{1}{2} \pi \Bigl( \frac{d}{ds} (b_{\gamma(t)}Z)_{\tilde{\tau}_t(s)} \Bigr|_{s=0}, \pmb{-} {\mathfrak s}(\tilde{\gamma}(t)) \cdot \frac{d}{ds} (b_{\gamma(t)}Z)_{\tilde{\tau}_t(s)} \Bigr|_{s=0} \Bigr).} 
\end{array}$$
Besides,
$$\begin{array}{cll}
\displaystyle{\frac{d}{ds} \Bigl(b_{\gamma(t)}Z \Bigr)_{\tilde{\tau}_t(s)} \Bigr|_{s=0}} & = & \displaystyle{\frac{d}{ds} b_{\gamma(t)} \bigl(\tau_t(s)\bigr) Z\bigl(\tau_t(s)\bigr) \Bigr|_{s=0}} \\
& = & j^1_{\gamma(t)} b_{\gamma(t)} \cdot Z_{*_{\gamma(t)}} Y_{\gamma(t)}
\end{array}$$
and
$$\frac{d}{dt}  j^1_{\gamma(t)} b_{\gamma(t)} \cdot Z_{*_{\gamma(t)}} Y_{\gamma(t)} \Bigr|_{t=0} = j^1_x j^1_\centerdot b_\centerdot \cdot Z_{**_{Y_x}} Y_{*_x} X_x.$$
Therefore,
$$\begin{array}{cll}
\displaystyle{\frac{d}{dt} \Bigl(\nabla_{b_x Y} b_{\centerdot} Z \Bigr)_{\tilde{\gamma}(t)} \Bigr|_{t=0}} & = & (m_{\frac{1}{2}})_* \pi_* \Bigl( \zeta \cdot Z_{**_{Y_z}} Y_{*_x} X_x, \\
& & \; \qquad \qquad (m_{-1})_* \circ (m_{-1})_{**} j^1_y{\mathfrak s} \cdot \zeta \cdot Z_{**_{Y_z}} Y_{*_x} X_x \Bigr).
\end{array}$$
We pursue our computation of $\nabla_{\xi X_x} \nabla_{b_x Y} b_{\centerdot} Z$ as in the proof of \tref{curvthm}, and obtain
$$\begin{array}{cll}
\nabla_{\xi X_x} \nabla_{b_x Y} b_{\centerdot} Z & = & \displaystyle{\frac{1}{4} \pi \Bigl[ \pi_* \Bigl( \zeta \cdot \IX, (m_{-1})_* \circ (m_{-1})_{**} \;j^1_y{\mathfrak s} \cdot \zeta \cdot \IX \Bigr),} \\
& & \quad \;\;\; \pi_* \Bigl((m_{-1}) \circ (m_{-1})_{**} \; {\mathbb S}(-I_y) \cdot \zeta \cdot \IX, \\
& & \qquad \quad \;\; (m_{-1}) \circ (m_{-1})_{*} \; {\mathbb S}(-I_y) \cdot j^1_y{\mathfrak s} \cdot \zeta \cdot \IX\Bigr) \Bigr].
\end{array}$$
Thus $\nabla_{\xi X_x} \nabla_{b_x Y} b_{\centerdot} Z$ admits the same expression as $\xi(\nabla_{X_x} \nabla_Y Z)$ with $\zeta$ instead of ${\mathbb S}(\xi)$. We claim that $\zeta$ is affine if and only it coincides with ${\mathbb S}(\xi)$. Indeed, since $\zeta \cdot \IX$ and ${\mathbb S}(\xi) \cdot \IX$ are in the same $(p \times p_* \times p_{**})$-fiber, there exists a $W \in T_yM$ such that 
$${\mathbb S}(\xi) \cdot \IX = A^{\zeta \cdot \IX}_{\p} (W) \stackrel{\rm not}{=} \zeta \cdot \IX \pmb{+} W.$$
(cf.~(\ref{param-P-fiber}) in \aref{third-der}). Then the relation (\ref{action-sur-TM}) implies that 
$$\begin{array}{rcl}
\displaystyle{j^1_y{\mathfrak s} \cdot {\mathbb S}(\xi) \cdot \IX} & = & j^1_y{\mathfrak s} \cdot \zeta \cdot \IX \pmb{+} -W, \\
\displaystyle{{\mathbb S}(-I_y) \cdot {\mathbb S}(\xi) \cdot \IX} & = & {\mathbb S}(-I_y) \cdot \zeta \cdot \IX \pmb{+} - W,\\
\displaystyle{{\mathbb S}(-I_y) \cdot j^1_y{\mathfrak s} \cdot {\mathbb S}(\xi) \cdot \IX} & = & {\mathbb S}(-I_y) \cdot j^1_y{\mathfrak s} \cdot \zeta \cdot \IX \pmb{+} W.
\end{array}$$
Each one of the previous relations remain valid if the two elements of $T^3M$ appearing in each side of the equality are both multiplied by an even number of negative signs. Hence defining $\ell_i \in \EL(T^3M)$, $i =1, 2, 3$ by  
\begin{enumerate}
\item[-] $\ell_1 \cdot \IY = (m_{-1})_* \circ (m_{-1})_{**} \;j^1_y{\mathfrak s} \cdot \IY$,
\item[-] $\ell_2 \cdot \IY = (m_{-1}) \circ (m_{-1})_{**} \; {\mathbb S}(-I_y) \cdot \IY$,
\item[-] $\ell_3 \cdot \IY = (m_{-1}) \circ (m_{-1})_{*} {\mathbb S}(-I_y) \cdot j^1_y{\mathfrak s} \cdot \IY$,
\end{enumerate}
we see that 
$$\begin{array}{l}
\ell_i \cdot {\mathbb S}(\xi) \cdot \IX = \ell_i \cdot \zeta \cdot \IX \pmb{+} - W \qquad i = 1, 2 \\
\ell_3 \cdot {\mathbb S}(\xi) \cdot \IX = \ell_3 \cdot \zeta \cdot \IX \pmb{+} W.
\end{array}$$
This implies that
$$\begin{array}{l}
\pi_* \Bigl( {\mathbb S}(\xi) \cdot \IX, \ell_1 \cdot {\mathbb S}(\xi) \cdot \IX\Bigr) = \pi_* \Bigl(\zeta \cdot \IX, \ell_1 \cdot \zeta \cdot \IX \Bigr) \pmb{+} 2W \\
\pi_* \Bigl(\ell_2 \cdot {\mathbb S}(\xi) \cdot \IX, \ell_3\cdot {\mathbb S}(\xi) \cdot \IX \Bigr) = \pi_* \Bigl(\ell_2 \cdot \zeta \cdot \IX, \ell_3 \cdot \zeta \cdot \IX \Bigr) \pmb{+} -2W.
\end{array}$$
Hence
$$\xi \Bigl(\nabla_{X_x} \nabla_Y Z \Bigr) - \nabla_{\xi X_x} \nabla_{b_x Y} b_{\centerdot} Z,$$
which coincides with 
$$\begin{array}{l}
\displaystyle{\frac{1}{4} \pi \Bigl[ \pi_{*} \Bigl({\mathbb S}(\xi) \cdot \IX, \ell_1 \cdot {\mathbb S}(\xi) \cdot \IX \Bigr), \pi_{*} \Bigl( \ell_2 \cdot {\mathbb S}(\xi) \cdot \IX, \ell_3 \cdot {\mathbb S}(\xi) \cdot \IX \Bigr) \Bigr] - }\\
\displaystyle{\frac{1}{4} \pi \Bigl[ \pi_* \Bigl( \zeta \cdot \IX, \ell_1 \cdot \zeta \cdot \IX \Bigr), \pi_* \Bigl(\ell_2 \cdot \zeta \cdot \IX, \ell_3 \cdot \zeta \cdot \IX\Bigr) \Bigr], }
\end{array}$$
is equal to $W$, i.e. we have shown that
\begin{equation}\label{nabla-nabla}
\xi \Bigl(\nabla_{X_x} \nabla_Y Z \Bigr) - \nabla_{\xi X_x} \nabla_{b_x Y} b_{\centerdot} Z = \Pi \Bigl( {\mathbb S}(\xi) \cdot \IX, \zeta \cdot \IX\Bigr).
\end{equation}
Hence $\zeta$ is affine is and only if $W = 0$, that is if and only $\zeta \cdot \IX = {\mathbb S}(\xi) \cdot \IX$ for all $\IX \in T_x^3M$ which in turn implies that $\zeta = {\mathbb S}(\xi)$.
\cqfd

\section{The curvature as a measure of the integrability of affine jets}

This section is devoted to proving that the affine $(1,1,1)$-jet ${\mathbb S}(\xi)$ extending $\xi$, which exists and is unique, as proven in the previous section, is a genuine $3$-jet depending on whether the $1$-jet $\xi$ does preserve both the covariant derivative of the torsion and the curvature. More precisely, a $(1,1,1)$-jet in $\b^{(1,1,1)}(\PP(M))$ is holonomic if and only if it is symmetric, that is, preserved by the involutions $\kappa$ and $\kappa_*$. Moreover, a $(1,1,1)$-jet is invariant under $\kappa$ (\rp $\kappa_*$) if and only if its first order preserves the covariant derivative of the torsion (\rp the curvature). The first statement, which is the content of the next proposition, is obtained from the relation (\ref{torsion}) by differentiating both sides. The second one, \pref{prop-curvature}, involves computing the curvature tensor, evaluated on vectors $X_x, Y_x, Z_x \in T_xM$, in terms of the second derivatives of $X, Y, Z$. 

\begin{prop}\label{kappa*} Let $\xi \in \b^{(1)}(\PP(M))$ and suppose $S(\xi)$ is holonomic. If $x = \alpha(\xi)$, let $X_x$, $Y_x$ and $Z_x$ be three vectors in $T_xM$ that extend to vector fields $X$, $Y$ and $Z$. Let also ${\mathfrak X} = Z_{**_{Y_x}} Y_{*_x} X_x$ in $T^3M$, then 
\begin{equation}\label{kappa*/nablaT}
\Pi \Bigl({\mathbb S}(\xi) \cdot {\mathfrak X}, \kappa_*({\mathbb S}(\xi)) \cdot {\mathfrak X}\Bigr) = \xi \Bigl( (\nabla_{Z_x}T^\nabla) (Y_x, X_x) \Bigr) - \Bigl(\nabla_{\xi Z_z}T^\nabla\Bigr)(\xi Y_x, \xi X_x),
\end{equation}
\end{prop}
Thus, when $S(\xi)$ is $\kappa$-invariant, the affine extension ${\mathbb S}(\xi)$ is $\kappa_*$-invariant if and only if $\xi$ preserves $\nabla T^\nabla$. In particular, ${\mathbb S}(\xi)$ is automatically $\kappa_*$-invariant when the connection $\nabla$ is torsionless.\\

\Pf Let $t \in (-\varepsilon, \varepsilon) \mapsto (\xi_t, \ZE_t) \in \b^{(1)}(\PP(M)) \times_{(\alpha, p^2)} T^2M$ be a smooth path whose first component $\xi_t$ is tangent to $\d_{\xi_0}$ at $t = 0$ with non-vanishing velocity vector $\dot{\xi}_0$. Set ${\mathfrak X} = \frac{d\ZE_t}{dt}|_0$, $X_t = p(\ZE_t)$, $Y_t = p_*(\ZE_t)$, $X_x = X_0$, $Y_x = Y_0$, $Z_x = \alpha_{*_{\xi_0}}\dot{\xi}_0 = p_* \circ p_* ({\mathfrak X})$, $\Y = \frac{dX_t}{dt}|_0 = p_*({\mathfrak X})$ and $\X = \frac{dY_t}{dt}|_0 = p_{**}({\mathfrak X})$. Then, equation (\ref{torsion}) holds for any $t \in (-\varepsilon, \varepsilon)$~: 
\begin{equation}\label{torsion-t}
\pi \Bigl(S(\xi_t) \cdot \ZE_t, \kappa(S(\xi_t)) \cdot \ZE_t \Bigr) = \xi_t(T^\nabla(Y_t, X_t)) - T^\nabla (\xi_t(Y_t), \xi_t(X_t)).
\end{equation}
Notice that 
$$\frac{d}{dt}S(\xi_t) \cdot \ZE_t \Bigr|_{t=0} = \rho^{(1,1)}_*\Bigl(\frac{d}{dt}S(\xi_t)\Bigr|_{t=0}, \frac{d}{dt}\ZE_t \Bigr|_{t=0} \Bigr) = \rho^{(1,1)}_*\Bigl(S_{*_{\xi_0}}(\dot{\xi}_0), \IX \Bigr) = {\mathbb S}(\xi_0) \cdot \IX,$$
thanks to the hypothesis that $\dot{\xi}_0 \in \d_{\xi_0}$ and the fact that $D({\mathbb S}(\xi)) = S_{*_\xi} (\d_{\xi})$. Similarly,
$$\frac{d}{dt}\kappa(S(\xi_t)) \cdot \ZE_t \Big|_{t=0} = \kappa_*({\mathbb S}(\xi_0)) \cdot \IX,$$
thanks to \rref{kappa*-diff} which implies that $\rho^{(1,1)}_{*}(\kappa^M_*(S_{*_{\xi}}(\dot{\xi}_0)), \IX) = \kappa_*({\mathbb S}(\xi)) \cdot \IX$.
Thus, the derivative with respect to $t$ of (\ref{torsion-t}), evaluated at $t=0$, yields 
\begin{multline}\label{torsion-der}
\pi_*\Bigl({\mathbb S}(\xi) \cdot {\mathfrak X}, \kappa_*({\mathbb S}(\xi)) \cdot {\mathfrak X}\Bigr) = \\
S(\xi) \cdot T^\nabla_{*_{(Y_x, X_x)}}(\X, \Y) -_* T^\nabla_{*_{(\xi Y_x, \xi X_x)}}(S(\xi) \cdot \X, S(\xi) \cdot \Y),
\end{multline}

Observe that, thanks to the hypothesis that $\xi$ preserves the torsion, both sides of (\ref{torsion-der}) are vectors in $T_{0_M}TM = TM \oplus TM$. Whence their vertical component are well-defined and agree, that is,
\begin{multline}
p^v \biggl( \pi_*\Bigl({\mathbb S}(\xi) \cdot {\mathfrak X}, \kappa_*({\mathbb S}(\xi)) \cdot {\mathfrak X}\Bigr) \biggr) = \\
p^v \biggl(S(\xi) \cdot T^\nabla_{*_{(Y_x, X_x)}}(\X, \Y) -_* T^\nabla_{*_{(\xi Y_x, \xi X_x)}}(S(\xi) \cdot \X, S(\xi) \cdot \Y) \biggr).
\end{multline}

Notice that $p^v$ coincides with $\widetilde{\nabla} : T^2M \to TM$, the vertical projection introduced in \lref{tilde-nabla}. In particular, the right-hand side of the previous relation coincides with the difference of the vertical projections of each term. Suppose that $\X = Y_{*_x} Z_x$ and $\Y = X_{*_x}Z_x$, where $X$ and $Y$ are local vector fields on $M$ and observe that 
$$\widetilde{\nabla} \circ T^\nabla_{*_{(Y_x, X_x)}} (\X, \Y) = \nabla_{Z_x} \bigl(T^\nabla(Y, X)\bigr),$$
Indeed, $\widetilde{\nabla} (\W) = \nabla_{V_x} U$ when $\W = U_{*_x} V_x$. In particular, if $\X$ and $\Y$ are both tangent to the horizontal distribution $\h^\nabla = \widetilde{\nabla}^{-1}(0_{TM})$, then 
$$\widetilde{\nabla} \circ T^\nabla_{*_{(Y_x, X_x)}} (\X, \Y) =  \bigl(\nabla_{Z_x}T^\nabla\bigr) (Y_x, X_x).$$

Moreover, the action of an affine jet $S(\xi)$ on $T^2M$ commutes with the vertical projection $\widetilde{\nabla}$, that is
$$\widetilde{\nabla}(S(\xi) \cdot \X) = \xi \, \widetilde{\nabla}(\X).$$ 
Altogether, under the hypothesis that ${\mathfrak X}$ is such that $p_*({\mathfrak X})$ and $p_{**}({\mathfrak X})$ both belong to $\h^\nabla$, (\ref{torsion-der}) becomes
\begin{equation}\label{torsion-der-bis}
\Pi \Bigl({\mathbb S}(\xi) \cdot {\mathfrak X}, \kappa_*({\mathbb S}(\xi)) \cdot {\mathfrak X}\Bigr) = \xi \bigl((\nabla_{Z_x} T^\nabla) (Y_x, X_x) \bigr) - \bigl(\nabla_{\xi Z_x} T^\nabla\bigr) (\xi Y_x, \xi X_x),
\end{equation}
where we have used the fact that $p^v \circ \pi_*$ coincides with $\Pi$. \\

Finally, the assumption that $\X$ and $\Y$ belong to the horizontal distribution $\h^\nabla$ is not restrictive due to the fact that the left-hand side of (\ref{torsion-der-bis}) does only depend on $X_x$, $Y_x$ and $Z_x$ (cf.~\rref{vert-part}). 

\cqfd

\begin{rmk} Removing the hypothesis that $\xi$ preserves the torsion, we can still prove that
\begin{multline}\label{kappa*-gen}
\widetilde{\nabla} \Bigl( \pi_* \bigl({\mathbb S}(\xi) \cdot {\mathfrak X}, \kappa_*({\mathbb S}(\xi)) \cdot {\mathfrak X}\bigr) \Bigr) = \\
\xi \bigl((\nabla_{Z_x} T^\nabla) (Y_x, X_x) \bigr) - \bigl(\nabla_{\xi Z_x} T^\nabla\bigr) (\xi Y_x, \xi X_x).
\end{multline}
\end{rmk}

Now we will see that ${\mathbb S}(\xi)$ is $\kappa$-invariant if and only if $\xi$ preserves the curvature of $\nabla$. This will thus show that ${\mathbb S}(\xi)$ is a $3$-jet if and only if $\xi$ preserves the following three tensors~: the torsion of $\nabla$, its covariant derivative and the curvature of $\nabla$. 

\begin{prop}\label{prop-curvature} Let $\xi \in \b^{(1)}(\PP(M))$ be a $1$-jet that preserves the torsion of $\nabla$, let $x = \alpha(\xi)$ and let $X_x$, $Y_x$ and $Z_x$ be three vectors in $T_xM$ that extend to vector fields $X$, $Y$ and $Z$. Let also ${\mathfrak X} = Z_{**_{Y_x}} Y_{*_x} X_x$ in $T^3M$. Then
\begin{equation}\label{curvature}
\Pi \Bigl({\mathbb S}(\xi) \cdot {\mathfrak X}, \kappa ({\mathbb S}(\xi)) \cdot {\mathfrak X} \Bigr) = \xi \Bigl(R^\nabla(X_x, Y_x)Z_x\Bigr) - R^\nabla(\xi X_x, \xi Y_x) \xi Z_x.
\end{equation}
In particular, if $\xi$ preserves the torsion tensor $T^\nabla$ and the curvature tensor $R^\nabla$, the affine extension ${\mathbb S}(\xi)$ is $\kappa$-invariant. 
\end{prop}

\Pf Write ${\mathbb S}(\xi) = j^1_x j^1_\centerdot b_\centerdot$, with $b_x$ (\rp $b_{x'}$) a local bisection of $\b^{(1)}(\PP(M))$ tangent to $\d_\xi$ (\rp $\d_{b_x(x')}$). Observe that the vector fields $b_x Y$ and $b_x Z$ extend the vectors $\xi Y_x$ and $\xi Z_x$ respectively and can therefore be used to compute the curvature, as is done below. 
$$\begin{array}{lll}
&& \xi \Bigl(R^\nabla(X_x, Y_x)\; Z_x\Bigr) - R^\nabla(\xi X_x, \xi Y_x) \; \xi Z_x \\
& = & \xi \Bigl( \nabla_{X_x} \nabla_Y Z - \nabla_{Y_x} \nabla_X Z - \nabla_{[X, Y]_x} Z\Bigr) \\
&& - \Bigl( \nabla_{\xi X_x} \nabla_{b_x Y} b_x Z -  \nabla_{\xi Y_x} \nabla_{b_x X} b_x Z - \nabla_{[b_x X, b_x Y]_y} b_x Z \Bigr). 
\end{array}$$
The vector fields $X$ and $Y$ may be chosen so that their bracket at $x$ vanishes, or equivalently that $\kappa(X_{*_x} Y_x) = Y_{*_x} X_x$. Moreover, the assumption that $\xi$ preserves the torsion $T^\nabla$ implies that the bracket $[b_x X, b_x Y]_y$ vanishes as well, as shown below~: 
$$\begin{array}{lll}
[b_x X, b_x Y]_y & = & \pi \Bigl( \bigl(b_x Y\bigr)_{*_y} (\xi X_x), \kappa \bigl((b_x X)_{*_y} (\xi Y_x) \bigr)\Bigr) \\
& = &  \pi \Bigl(j^1_xb_x \cdot Y_{*_x} X_x, \kappa \bigl(j^1_x b_x \cdot X_{*_x} Y_x \bigr) \Bigr)\\
& = &  \pi \Bigl(S(\xi) \cdot Y_{*_x} X_x, \kappa (S(\xi)) \cdot \kappa(X_{*_x} Y_x) \Bigr)\\
& = &  \pi \Bigl(S(\xi) \cdot Y_{*_x} X_x, \kappa (S(\xi)) \cdot Y_{*_x} X_x\Bigr)\\
& = & \xi \Bigl(T^\nabla(X_x, Y_x) \Bigr) - T^\nabla \Bigr(\xi X_x, \xi Y_x \Bigr)  \quad(\mbox{\pref{prop-torsion}}) \\
& = & 0\quad(\mbox{if $\xi$ preserves the torsion}).
\end{array}$$
Thus 
$$\begin{array}{lll}
&& \xi \Bigl(R^\nabla(X_x, Y_x)\; Z_x\Bigr) - R^\nabla(\xi X_x, \xi Y_x) \; \xi Z_x \\
& = & \xi \Bigl( \nabla_{X_x} \nabla_Y Z \Bigr) - \Bigl( \nabla_{\xi X_x} \nabla_{b_x Y} b_x Z \Bigr) \\
&& - \xi \Bigl( \nabla_{Y_x} \nabla_X Z \Bigr) + \Bigl( \nabla_{\xi Y_x} \nabla_{b_x X} b_x Z \Bigr).
\end{array}$$
The relation (\ref{nabla-nabla}) established in the proof of \pref{uae-order3} and applied to $\zeta = j^2_x b_x$ which is also en extension of $S(\xi)$ under the hypothesis that $S(\xi)$ is holonomic yields
$$\begin{array}{lll}
&& \xi \Bigl(R^\nabla(X_x, Y_x)\; Z_x\Bigr) - R^\nabla(\xi X_x, \xi Y_x) \; \xi Z_x \\
& = & \Pi \Bigl( {\mathbb S}(\xi) \cdot \IX, \zeta \cdot \IX \Bigr) - \Pi \Bigl( {\mathbb S}(\xi) \cdot \kappa(\IX), \zeta \cdot \kappa(\IX) \Bigr) \\
& = & \Pi \Bigl( {\mathbb S}(\xi) \cdot \IX, \zeta \cdot \IX \Bigr) - \Pi \Bigl( \kappa ( {\mathbb S}(\xi)) \cdot \IX, \kappa(\zeta) \cdot \IX \Bigr) \\
& = & \Pi \Bigl( {\mathbb S}(\xi) \cdot \IX, \zeta \cdot \IX \Bigr) - \Pi \Bigl( \kappa ( {\mathbb S}(\xi)) \cdot \IX, \zeta \cdot \IX \Bigr) \\
& = & \Pi \Bigl( {\mathbb S}(\xi) \cdot \IX, \kappa ( {\mathbb S}(\xi)) \cdot \IX \Bigr) 
\end{array}$$
For the second equality, we have used the fact that $\kappa(Z_{**_{Y_x}}Y_{*_x} X_x) = Z_{**_{X_x}}X_{*_x} Y_x$, for the third the relation (\ref{i-kappa-bis}) in \aref{third-der} and for the fourth, the fact that $\kappa(j^2_xb_x) = j^2_x b_x$.
\cqfd

\begin{rmk}
It is possible to write a more general formula for (\ref{curvature}) when the assumption that $\xi$ preserves the torsion is not fulfilled. 
\end{rmk}

\begin{cor} Suppose the torsion of $\nabla^{\mathfrak s}$ vanishes identically. Then the curvature $R$ of $\nabla^{\mathfrak s}$ may be recovered from (\ref{curvature}) by particularizing $\xi$. Indeed, considering some linear homothety $m_a : TM \to TM : X_x \mapsto a X_x$ with $a \neq \pm 1$, we obtain the following expression for the curvature~:
$$R (X_x, Y_x)Z_x = \frac{1}{a (1 - a^2)} \Pi \Bigl({\mathbb S}(m_a) \cdot {\mathfrak X}, \kappa ({\mathbb S}(m_a)) \cdot {\mathfrak X}\Bigr),$$
where $X$, $Y$, $Z \in \IX(M)$ and $\IX = Z_{**_{Y_x}} Y_{*_x} X_x$.
\end{cor}

The subgroupoid of $1$-jets that preserve the torsion and the curvature can be characterized in terms of $\d^\fs$ as follows.

\begin{dfn}\label{int-locus} The integrability locus of a distribution $\d$ on a manifold $W$ is the set of points $w \in W$ such that the bracket of any pair of local vector fields near $w$ tangent to $\d$ belongs to $\d$ at $w$, or
$$\In(\d) = \Bigl\{w \in W; A, B \in \Gamma \d \Longrightarrow [A, B]_w \in \d\Bigr\}.$$
When $w \in \In(\d)$, we also say that $\d$ is flat at $w$.
\end{dfn}
The following is a classical result of differential geometry~:

\begin{prop}\label{int/oscul} Let $\d$ be a distribution on a manifold $W$. If $w \in \In(\d)$, there exists an embedded submanifold $F \subset W$ which is osculatory to $\d$ at $w$. The second order jet of $F$ at $w$ is unique. 
\end{prop}

\begin{lem}\label{Int-kappa} Let ${\mathfrak s}$ be a symmetry jet on $M$. Then $\xi \in \In(\d^{\mathfrak s})$ if and only if $\kappa({\mathbb S}(\xi)) = {\mathbb S}(\xi)$. In particular $\In(\d^\fs)$ is the set of $1$-jets that preserve the torsion and the curvature, that is $\In(\d^\fs) = \b(T^\fs, R^\fs)$ with the notation introduced in \dref{B(Q)}.
\end{lem}

\Pf Observe that $\xi \in \Int(\d^\fs)$ if and only if there exists a local bisection $b$ of $\b^{(1)}(\PP(M))$ that is osculatory to $\d_b^\fs$ at $x = \alpha(\xi)$. In other terms, $Tb = D(j^1b)$ is tangent to $\d_b^\fs = D(S\circ b)$ at $x$ or $j^1b$ is tangent to $S \circ b$ at $x$. The latter statement is equivalent to $(j^1b)_*(T_xM) = (S \circ b)_*(T_xM) = \d_b^\fs$ that is to $j^2_xb = {\mathbb S}(b(x))$. The latter equality says that the affine jet ${\mathbb S}(b(x))$ is $\kappa$-invariant.
\cqfd

\begin{rmk}
Observe that the integrability locus of $\d^{\mathfrak s}$ is a subgroupoid of $\b^{(1)}(\PP(M))$ that contains $I \cup -I$. It is all of $\b^{(1)}(\PP(M))$ if and only if $R^{\mathfrak s}$ vanishes identically. It would be interesting to know whether $\In(\d^{\mathfrak s})$ determines $R^{\mathfrak s}$ in general.
\end{rmk}

\section{Why there are no other tensors than the torsion and the curvature~?}\label{why-there-are-no}

The purpose of the section is to explain why, if we pursue this procedure no new tensor appear. First of all, (semi-holonomic) affine extensions of all order exist. At order four, the affine extension $S^{(1,1,1,1)}(\xi)$ of $\xi \in \b^{(1)}(\PP(M))$ is defined through 
$$D(S^{(1,1,1,1)}(\xi)) = {\mathbb S}_{*_{\xi}}(\d^{\mathfrak s}_\xi).$$ 
Set ${\mathfrak D}^{\fs}_{S(\xi)} \stackrel{\rm def}{=} S_{*_\xi}(\d^{\mathfrak s}_\xi)$. The relation
\begin{equation}\label{iterated-integrability}
[{\mathfrak D}^{\fs}, {\mathfrak D}^{\fs}]_{S(\xi)} = [S_{*} \d^{\mathfrak s}, S_* \d^{\mathfrak s}]_{S(\xi)} = S_{*_\xi} [\d^{\mathfrak s}, \d^{\mathfrak s}]_\xi
\end{equation}
implies that if $\xi$ belongs to the integrability locus of $\d^{\mathfrak s}$ (cf.~\dref{int-locus}), then $S(\xi)$ automatically belongs to that of ${\mathfrak D}^{\fs}$. As a consequence, if $\xi$ preserves $T$, $\nabla T$ and $R$, then $S^{(1,1,1,1)}(\xi)$ is automatically $\kappa$-invariant (\pref{kappa-inv}). This says that the $\kappa$-invariance of $S^{(1,1,1,1)}(\xi)$ is not anymore obstructed by some tensor. In addition, differentiating the relations (\ref{kappa*/nablaT}) and (\ref{curvature}) imply that the $\kappa_*$ and $\kappa_{**}$-invariance of $S^{(1,1,1,1)}(\xi)$ depends on the $\xi$-invariance of $\nabla R$ and $\nabla \nabla T$ (cf.~\pref{next-order-inv}). This process can be iterated, showing that if $\xi$ preserves the various covariant derivatives of the torsion and the curvature, then the affine extension of $\xi$ is holonomic at any order. \\

In order to write proofs of the results announced previously, the following result, due to Tapia, is particularly useful. It generalizes to Lie groupoids the well-known theorem of E.~Cartan according to which any closed subgroup of a Lie group is an embedded Lie subgroup~:

\begin{thm}\label{closed-subgr}\cite{Ta} Let $G \rightrightarrows M$ be a locally trivial Lie groupoid. Then any closed (algebraic) subgroupoid of $G$ is an embedded Lie subgroupoid.
\end{thm} 
Locally trivial Lie groupoids are those Lie groupoids for which the map $\alpha \times \beta : G \to M \times M$ is a surjective submersion. This condition is certainly satisfied by $\b^{(1)}(\PP(M))$ (cf.~\lref{loc-triv}).

\begin{dfn}\label{B(Q)} Given a family $\{Q_1, ..., Q_k\}$ of $(1, p)$-tensors on $M$, the subgroupoid of $\b^{(1)}(\PP(M))$ consisting of those $1$-jets $\xi$ that preserve all $Q_i$'s in the sense that 
$$Q_i (\xi X_1, ..., \xi X_{p}) = \xi (Q_i(X_1, ..., X_{p}))$$ 
for all $X_1,..., X_p \in \IX(M)$ and all $i = 1, ..., k$ is denoted by $\b(Q_1, ..., Q_k)$.
\end{dfn}

\begin{cor} For any choice of tensors $Q_1, ..., Q_k$, the (algebraic) subgroupoid $\b(Q_1, ..., Q_k)$ is an embedded Lie subgroupoid. In particular, given a symmetry jet ${\mathfrak s}$, this yields a set of embedded Lie subgroupoids ${\mathcal B}(T)$, ${\mathcal B}(R)$, ${\mathcal B}(T, \nabla T, R)$, etc ..., for $T$ (\rp $R$) the torsion (\rp curvature) of the associated connection $\nabla^{\mathfrak s}$. 
\end{cor} 

\begin{lem}\label{d-tangent-to-b(Q)} Given a tensor $Q$ on $M$, the distribution $\d^{\mathfrak s}$ is tangent to ${\mathcal B}(Q)$ at $\xi$ if and only if $\xi$ preserves the first covariant derivative of $Q$. 
\end{lem}

\Pf A vector $X_\xi = \frac{d \xi_t}{dt}\bigr|_{t=0}$ in $\d^{\mathfrak s}_\xi$ belongs to $T_\xi\b(Q)$ if and only if 
$$\frac{d}{dt} \xi_t Q \Bigl(X_1^t, ..., X_p^t \Bigr) \Bigl|_{t=0} = \frac{d}{dt} Q \Bigl(\xi_t X_1^t, ..., \xi_t X_p^t \Bigr) \Bigr|_{t=0},$$
for any paths $t \mapsto X_i^t$ of vectors in $T_{\alpha(\xi_t)}M$. Equivalently,
$$S(\xi) \cdot Q_{*_{(X_1^0, ..., X_p^0)}} \Bigl( \X_1, ..., \X_p \Bigr) = Q_{*_{(\xi X_1^0, ..., \xi X_p^0)}} \Bigl( S(\xi) \cdot \X_1, ..., S(\xi) \cdot \X_p\Bigr),$$
where $\X_i = \frac{d}{dt}X_i^t \bigl|_{t=0}$. Projecting horizontally with respect to $\h^\nabla$ and assuming, without loss of generality, that $\X_i$ is tangent to the horizontal distribution $\h^\nabla$, as is done in the proof of \pref{kappa*}, we obtain 
$$\xi \Bigl( \nabla_{Z_x} Q (X_1^0, ...,X_p^0)\Bigr) = \nabla_{\xi Z_x} Q \Bigl( \xi X_1^0, ..., \xi X_p^0 \Bigr),$$
where $Z_x = \alpha_{*_\xi}(X_\xi)$.
\cqfd

Now we prove existence of affine $(1,1,1,1)$-jets. Let $\xi \in \b^{(1)}(\PP(M))$, and let $a_\xi$ denote some local bisection $U \ni x_1 \mapsto a_\xi(x_1) \in \b^{(1)}(\PP(M))$ tangent to $\d^{\mathfrak s}$ at $\xi$, or equivalently such that $j^1_xa_\xi = S(\xi)$, for $x = \alpha(\xi)$. Similarly, for each $1$-jet $a_\xi(x_1)$, $x_1\in U$,  let $a^2_\xi(x_1)$ (instead of $a_{a_\xi(x_1)}$) denote a local bisection tangent to $\d^{\mathfrak s}$ at the point ${a_\xi(x_1)}$. Iterating this procedure, we obtain a family of local bisections $x_k \mapsto a^k_\xi(x_1, ..., x_{k-1})(x_k)$, $k = 1, 2, ...$ of $\b^{(1)}(\PP(M))$ such that 
\begin{enumerate}
\item[-] $a^k_\xi(x_1, ..., x_{k-1})(x_{k-1}) = a^{k-1}_\xi(x_1, ..., x_{k-2})(x_{k-1})$,
\item[-] $j^1_{x_{k-1}} a^k_\xi(x_1, ..., x_{k-1}) = S \bigl(a^{k-1}_\xi(x_1, ..., x_{k-2})(x_{k-1}) \bigr).$
\end{enumerate}
\begin{prop}\label{affine-ext-order-4} Let $\xi \in \b^{(1)}(\PP(M))$, with $\alpha (\xi) = x$. Then the $(1,1,1,1)$-jet 
$$S^{(1,1,1,1)}(\xi) = j^1_xj^1_{x_1}a^2_\xi(x_1) = j^1_x({\mathbb S} \circ a_\xi).$$ 
is affine in the sense that for any $x \in M$ and $X, Y, Z, T \in \IX(M)$, we have
$$\xi \Bigl( \nabla_{X_x} \nabla_Y \nabla_Z T \Bigr) = \nabla_{(\xi X_x)} \nabla_{(a_{\xi}(x_1) Y_{x_1})} \nabla_{(a^2_{\xi}(x_1, x_2) Z_{x_2})} \bigl(a^3_\xi(x_1, x_2, x_3) T_{x_3} \bigr).$$  
\end{prop} 
\Pf \pref{saffine} implies the following sequence of equalities~:
$$\begin{array}{lll}
\xi \Bigl( \nabla_{X_x} \nabla_Y \nabla_Z T \Bigr) & = & \nabla_{\xi X_x} \Bigl( a_\xi \bigl(\nabla_Y \nabla_Z T \bigr) \Bigr), \\
a_\xi(x_1) \Bigl(\nabla_{Y_{x_1}} \nabla_Z T \Bigr) & = & \nabla_{a_\xi(x_1) Y_{x_1}} \Bigl ( a^2_\xi(x_1) \bigl(\nabla_Z T \bigr) \Bigr), \\
a^2_\xi(x_1)(x_2) \Bigl(\nabla_{Z_{x_2} }T \Bigr) & = & \nabla_{a^2_\xi(x_1)(x_2) Z_{x_2}} \Bigl( a^3_\xi(x_1, x_2) Z \Bigr).
\end{array}$$
\cqfd

\begin{rmk}
Observe that 
$$D(S^{(1,1,1,1)}(\xi)) = {\mathbb S}_{*_\xi}(\d_\xi).$$
\end{rmk}

\begin{rmk} This argument applies to any order. Let us denote $k\cdot (1)$ a sequence $(1, ..., 1)$ with $k$ times the number $1$. Given $\xi \in \b^{(1)}(\PP(M))$ with $x = \alpha(\xi)$, the $k \cdot (1)$-jet $S^{k \cdot (1)} (\xi) = j^1_xj^1_{x_1} ... j^1_{x_{k-1}} a^{k}(x_1, ..., x_{k-1})$ is affine in the sense that its $k$ parts of order $k-1$ are affine and agree and it preserves the $k$-th power of $\nabla$. 
\end{rmk}

\begin{prop}\label{next-order-inv} Let $\xi \in \b^{(1)}(\PP(M))$ be such that ${\mathbb S}(\xi)$ is holonomic and let $S^{(1,1,1,1)}(\xi)$ be the affine $(1,1,1,1)$-jet extending $\xi$. Let ${\mathbb X} \in T^4_xM$ with $x = \alpha(\xi)$ and its four projections on $T_xM$ denoted by $X_x, Y_x, Z_x, T_x$. Then
\begin{multline}\label{nabla-nabla-T}
\Pi \Bigl(S^{(1,1,1,1)}(\xi) \cdot {\mathbb X}, \kappa_{**}(S^{(1,1,1,1)}(\xi)) \cdot {\mathbb X}\Bigr) \\ 
= \xi \Bigl( \nabla_{T_x}\nabla T^\nabla (Z_x, Y_x, X_x) \Bigr) - \Bigl(\nabla_{\xi T_x} \nabla T^\nabla \Bigr)(\xi Z_x, \xi Y_x, \xi X_x)
\end{multline}
\begin{multline}\label{nabla-R}
\Pi \Bigl(S^{(1,1,1,1)}(\xi) \cdot {\mathbb X}, \kappa_{*}(S^{(1,1,1,1)}(\xi)) \cdot {\mathbb X}\Bigr) \\ 
= \xi \Bigl( \nabla_{T_x} R^\nabla (X_x, Y_x, Z_x) \Bigr) - \Bigl( \nabla_{\xi T_x} R^\nabla \Bigr)(\xi X_x, \xi Y_x, \xi Z_x),
\end{multline}
where $\Pi : T^4M \times_{(P,P)}T^4M \to TM$ with $P = p \times p_* \times p_{**} \times p_{***}$ is defined in a similar fashion as for the case of $T^3M$ (cf.~(\ref{Pi}) in \aref{third-der}).
\end{prop}

\Pf The idea is, of course, to differentiate (\ref{kappa*/nablaT}) and (\ref{curvature}) with respect to $\xi$. Strictly speaking, this would not be possible because these two relations are only valid for $\xi$'s preserving the torsion. However, thanks to Tapia's \tref{closed-subgr}, we know that ${\mathcal B}(T)$ is a submanifold and we may thus differentiate the relations (\ref{nabla-nabla-T}) and (\ref{nabla-R}) in the direction of its tangent space. \\

As a first step, observe that the hypothesis that $\xi$ preserves the torsion and its first covariant derivative implies (cf.~\lref{d-tangent-to-b(Q)}) that 
$$\d_\xi \subset T_\xi {\mathcal B}(T).$$
In particular $\alpha_{*_\xi}$ restricted to $T_\xi \b(T)$ is a submersion. So given any ${\mathbb X} \in T^4_xM$, there exists a $X_\xi \in T_\xi {\mathcal B}(T)$ such that 
$$(X_\xi, {\mathbb X}) \in T_{(\xi, \IX)} \bigl({\mathcal B}(T) \times_{(\alpha, p^3)} T^3M\bigr).$$
Let $(- \varepsilon, \varepsilon) \to {\mathcal B}(T) \times_{(\alpha, p^3)} T^3M : t \mapsto (\xi_t, \IX_t)$ be a path tangent to $(X_\xi, {\mathbb X})$. Then the relation
\begin{equation}\label{abcd}
\Pi \bigl({\mathbb S}(\xi_t) \cdot {\mathfrak X}_t, \kappa_o({\mathbb S}(\xi_t)) \cdot {\mathfrak X}_t \Bigr) = 
\xi_t \bigl(Q (Z_t, Y_t, X_t) \bigr) - Q (\xi_t Z_t, \xi_t Y_t, \xi_t X_t),
\end{equation}
for $X_t = p \circ p (\IX_t)$, $Y_t = p_* \circ p(\IX_t)$, $Z_t = p_* \circ p_*(\IX_t)$, holds true for any $t \in (- \varepsilon, \varepsilon)$, where $(\kappa_o, Q)$ is either $(\kappa_*, \nabla T)$ or $(\kappa, R)$. Differentiating both sides with respect to $t$ and evaluating at $t = 0$ yields~:
\begin{multline}\label{abdc}
\Pi_*\bigl( S^{(1,1,1,1)}(\xi) \cdot {\mathbb X}, \kappa_{o_*}(S^{(1,1,1,1)}(\xi)) \cdot {\mathbb X} \bigr) = S(\xi) \Bigl(Q_* (Z_{*_x}T_x, Y_{*_x}T_x, X_{*_x}T_x) \Bigr) -_*\\
Q_* \Bigl(S(\xi) \cdot Z_{*_x}T_x, S(\xi) \cdot Y_{*_x}T_x, S(\xi) \cdot Z_{*_x}T_x \Bigr),
\end{multline}
where $T_x = \frac{d}{dt}p^3(\IX_t)\bigl|_{t=0}$. Because the left-hand side of (\ref{abcd}) vanishes at $t = 0$, the left-hand side of (\ref{abdc}) is a vector in $T_{0_x}TM = T_xM \oplus T_xM$ whose vertical component coincides with 
$$\Pi \bigl( S^{(1,1,1,1)}(\xi) \cdot {\mathbb X}, \kappa_{o_*}(S^{(1,1,1,1)}(\xi)) \cdot {\mathbb X} \bigr).$$
As to the right-hand side of (\ref{abdc}), its treatment is essentially contained in the proof of \pref{kappa*}. Indeed, its vertical component coincides with the difference of the vertical projection with respect to $\h^\nabla$ of each terms. So assuming the local vector fields $X$, $Y$ and $Z$ to be tangent to the horizontal distribution $\h^\nabla$, we reach the desired equalities.
\cqfd

\begin{prop}\label{kappa-inv} Suppose $\xi \in \b^{(1)}(\PP(M))$ preserves $T$, $\nabla T$ and $R$. Then the affine $(1,1,1,1)$-jet $S^{(1,1,1,1)}(\xi)$ is $\kappa$-invariant. In particular if $\xi$ preserves also $\nabla \nabla T$ and $\nabla R$, then $S^{(1,1,1,1)}(\xi)$ is holonomic.
\end{prop}

\Pf Let $\xi \in \b(T, \nabla T, R)$. Then $\xi$ belongs to $\In (\d^{\mathfrak s})$ (cf.~\pref{Int-kappa}), which implies that $S(\xi)$ belongs to $\In ({\mathfrak D}^{\mathfrak s})$ (cf.~(\ref{iterated-integrability})). Therefore there exists a local bisection $b$ in $\im S \subset \b^{(1,1)}_{nh}(\PP(M))$, thus of type $b = S \circ b_o$ for a local bisection $b_o$ of $\b^{(1)}(\PP(M))$, which is osculatory to ${\mathfrak D}^{\mathfrak s}$ at $S(\xi)$. Equivalently, the distributions $Tb$ and ${\mathfrak D}^{\mathfrak s}_b$ are tangent at $S(\xi)$ and, therefore, the corresponding local bisections $j^1b$ and ${\mathbb S} \circ b_o$ of the groupoid $\b^{(1,1,1)}_{nh}(\PP(M))$ are tangent at $x = \alpha(\xi)$~:
\begin{equation}\label{transitory}
(j^1b)_{*_x} (T_xM) = ({\mathbb S}_{*_\xi} \circ b_{o*_x}) (T_xM).
\end{equation}
On the other hand, the fact that $T_{S(\xi)}b = {\mathfrak D}^{\mathfrak s}_{S(\xi)}$ implies that 
$$(S_{*_\xi} \circ b_{o*_x}) (T_xM) = S_{*_\xi} (\d^{\mathfrak s}_\xi)$$ and therefore that $b_{o*_x} (T_xM) = \d^{\mathfrak s}_\xi$. Thus (\ref{transitory}) says that $j^2_xb = S^{(1,1,1,1)}(\xi)$. Hence the affine $(1,1,1,1)$-jet $S^{(1,1,1,1)}(\xi)$ is in fact a $(2,1,1)$-jet, whence is fixed by $\kappa$.
\cqfd

\begin{rmk} The subgroupoid $\b(T, R, ..., \nabla^k T, \nabla^k R, ...)$ is the largest subset of $\b^{(1)}(\PP(M))$ entirely foliated by the leaves of $\d^\fs$ and containing all of them. 
\end{rmk}

\section{Distributions on associated bundles}\label{hor-ass-bdle}


The groupoid $\b^{(1)}(\PP(M))$ acts on a vector bundle $E \to M$ exactly when the latter is associated to the frame bundle (cf.~\rref{associated}). We describe the horizontal distribution induced by a symmetry jet $\fs$ on $M$ as the collection of $-1$-eigenspaces of a canonical bundle map $TE \to TE$ over the identity, in the same spirit as for the distribution $\d^\fs$. \\

Let $\rho : \b^{(1)}(\PP(M)) \times_{(\alpha, \pi)} E \to E : (\xi, e) \mapsto \xi \cdot e$ be a linear groupoid action of $\b^{(1)}(\PP(M))$ on a vector bundle $\pi : E \to M$. Then a symmetry jet ${\mathfrak s}$ induces a linear connection $\nabla^{\mathfrak s}$ on $E$ through the formula~:  
\begin{equation}\label{sym-jet-->conn-on-E}
\nabla^{\mathfrak s}_{X_x} e = \frac{1}{2} \pi \Bigl(e_{*_x} X_x, \pmb{-}{\mathfrak s}(x) \cdot e_{*_x}X_x \Bigr),
\end{equation}
for $X_x$ a vector in $T_xM$ and $e$ is a local section of $E$ defined near $x$. The bold minus sign $\pmb{-}$ stands for $m_{-1} \circ L_{-I*}$, where $L_{-I*}$ is the differential of the action of the bisection $-I$ on $E$, that is 
$$L_{-I*} : TE \to TE : \frac{de_t}{dt}\Bigl|_{t=0} \mapsto \frac{d \rho(-I_{\pi(e_t)}, e_t)}{dt}\Bigl|_{t=0}.$$ 
The dot in formula (\ref{sym-jet-->conn-on-E}) denotes the derived action 
$$\rho^{(1)} : \b^{(1,1)} (M) \times_{(\alpha, \pi_*)} TE \to TE$$ (cf.~\dref{derived-action}) and the map $\pi$ is defined similarly to (\ref{i}). Simply observe that 
$$P : TE \to E \times TM : V_e \mapsto (p(V_e) = e, \pi_{*_e}(V_e))$$ is an affine fibration modeled on $E_{\pi(e)}$, the fiber of $\pi$ through $e$. Whence the map 
$$\pi : TE \times_{(P, P)} TE \to E : (V_e, V'_e) \mapsto \pi(V_e, V'_e).$$ 
It is not difficult to adapt the first part of the proof of \pref{symm-jet<-->conn} and show that (\ref{sym-jet-->conn-on-E}) defines a linear connection on $E$. \\

Now given a bisection $b$ of $\b^{(1)}(\PP(M))$ and a horizontal distribution $\d$ along $b$, we have two bundle automorphisms over $\phi_b : M \to M : x \mapsto \beta \circ b(x)$~:
$$\Phi_b : E \to E : e \mapsto \rho(b(\pi(e)), e)$$
$$\psi_{(b, \d)} : TM \to TM : X_x \mapsto \beta_* \circ \bigr(\alpha_*\big|_{\d_{b(x)}}\bigl)^{-1} (X_x),$$ 
both over $\phi_b : M \to M : x \mapsto \beta \circ b(x)$. Consider the map $\Psi_{(b, \d)} : TE \to TE$ over $\Phi_b$ defined by 
\begin{equation}\label{Psi}
\Psi_{(b, \d)} (V_e) = \rho_*\Bigl(\overline{\pi_{*}(V_e)}^{(\d_{b(\pi(e))}, \alpha_*)}, V_e\Bigr),
\end{equation}
where $\overline{\pi_{*}(V_e)}^{(\d_{b(\pi(e))}, \alpha_*)}$ denotes the lift of $\pi_{*}(V_e)$ in $\d_{b(\pi(e))}$ with respect to $\alpha_*$. It is also the unique vector $X$ in $\d_{b(\pi(e))}$ for which the pair $(X, V_e)$ belongs to the fiber product $TG \times_{(\alpha_*, \pi_*)}TE$ which is the source  of the map $\rho_*$. \\

\begin{figure}[h!]
\begin{center}
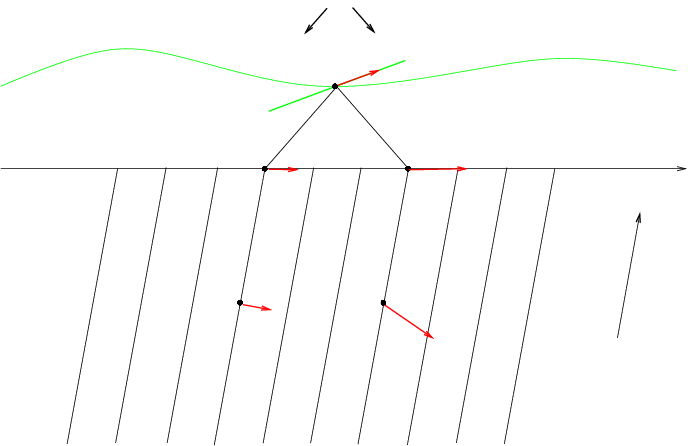 
\caption{The map $\Psi_{(b, \d)}$}
\end{center}
\end{figure}

It is easy to see that the map $\Psi_{(b, \d)}$ preserves the vertical tangent space $T^\pi E$, coincides with $(\Phi_b)_*$ vertically and with $\psi_{(b, \d)}$ horizontally, that is to say
$$\left\{
\begin{array}{lcl}\label{gen-second-deriv}
\Psi_{(b, \d)} \circ i^\pi = i^\pi \circ (\Phi_b)_*\\
\pi_* \circ \Psi_{(b, \d)} = \psi_{(b, \d)} \circ \pi_*,
\end{array}\right.$$
where $i^\pi : T^\pi E \to TE$ denotes the canonical inclusion. In particular, an invariant horizontal distribution on a groupoid $G$ that acts on a fibration $\pi : E \to M$ induces a representation of the group of bisections of $\b^{(1)}(\PP(M))$ into the group $GL(TE)$ of fiberwise linear diffeomorphisms of $TE$. \\

Notice also that if $\d$ (\rp $\d'$) is a horizontal distribution along a bisection $b$ (\rp $b'$) respectively, then 
\begin{equation}\label{multiplicativity}
\Psi_{(b, \d)} \circ \Psi_{(b', \d')} = \Psi_{(b \cdot b', \d \cdot \d')}.
\end{equation}

Now suppose $b = -I$ and $\d$ coincides with the distribution $\d^\fs$ on $\b^{(1)}(\PP(M))$ induced from a symmetry jet $\fs$ on $M$. Consider the bundle map
$$\Theta_{(b,\d)} \stackrel{\Def}{=} \Psi_{(b, Tb)} \circ \Psi_{(b, \d)} = m_{-1*} \circ \Psi_{(b, \d)} = \Psi_{(I, Tb \cdot \d)}.$$ 
It is a bundle map over the identity on $E$ which coincides with the identity on $T^\pi E$ and with the map $\psi_{(b, Tb)} \circ \psi_{(b, \d)} = \psi_{(b, \d)}$ horizontally.
Since $b \cdot b = I$, $Tb \cdot Tb = TI$, $\d \cdot \d = TI$ (cf.~\pref{e-iota}) and $Tb \cdot \d = \d \cdot Tb$, the map $\Theta_{(b,\d)}$ is involutive~:
$$\begin{array}{ccl}
\Theta_{(b,\d)} \circ \Theta_{(b,\d)} & = &  \Psi_{(b, Tb)} \circ \Psi_{(b, \d)} \circ  \Psi_{(b, Tb)} \circ \Psi_{(b, \d)} \\
& = & \Psi_{(I, Tb \cdot \d \cdot Tb \cdot \d)} \\
& = & \Psi_{(I, TI)} = \id.
\end{array}$$ 
This implies that each tangent space $T_eE$ splits into a direct sum of eigenspaces for the eigenvalues $+1$ and $-1$. Of course $T^\pi E$ is contained in the $+1$-eigenspace of $\Theta_{(b, \d)}$. Moreover, since $\psi_{(b, Tb)} = \id$ and $\psi_{(b, \d)} = -I$, we have $\pi_* \circ \Theta_{(b, \d)} = -I \circ \pi_*$, which implies that $T^\pi E$ coincides with the $+1$-eigenspace and the $-1$-eigenspace is thus a horizontal distribution denoted by $\h = \h^{\mathfrak s}$ on $E$~: 
$$\h = \Ker(\Theta_{(b,\d)} + I) = \im (- \Theta_{(b,\d)} + I)$$

\begin{rmk}\label{associated} The groupoid $\b^{(1)}(\PP(M))$ acts linearly on a vector bundle $\pi : E \to M$ if and only if the latter is associated to the principal bundle of frames of $TM$. 
\end{rmk}

\begin{lem} Let ${\mathfrak s}$ be a symmetry jet on $M$ and let $\d$ be the corresponding distribution on $\b^{(1)}(\PP(M))$. Then the horizontal distribution $\h^{\mathfrak s}$ on $E$ is the canonical horizontal distribution associated to the connection $\nabla^{\mathfrak s}$ (cf.~\rref{hor-dist}). 
\end{lem} 

\Pf Firstly, let us show that the map $\Psi_{(b, \d)}$ coincides with the $\rho^{(1,1)}$-action of ${\mathfrak s}$ on $TE$, that is, for any $V_e \in T_eE$~:
$$\Psi_{(b, \d)} (V_e) = \rho^{(1,1)} \Bigl({\mathfrak s}(\pi(e)), E \Bigr),$$
Suppose $\pi_*(V_e) = X \in T_xM$. Then,
$$\begin{array}{ccl}
\vspace{.2cm}
\Psi_{(b, \d)} (V_e) & = & \rho_{*_{(b(x), e)}} \Bigl(\overline{X}^{(\d_{b(x)}, \alpha_*)}, V_e\Bigr) \\
\vspace{.1cm}
& = & \displaystyle{\rho^{(1,1)}\Bigl({\mathfrak s}(x), V_e \Bigr)},\\
\end{array}$$
(cf.~\rref{derived-action-of-pair-groupoid}). It is proven now that a vector $V_e$ in $T_eE$ is in the eigenspace of $\Psi_{(b,\d)}$ for the eigenvalue $-1$ if and only if $\nabla^{\mathfrak s}_{X_x} e = 0$, for $X_x = \pi_*(V_e)$ and $e$ a local section of $\pi : E \to M$ such that $e_{*_x}X_x = V_e$~: 

$$\nabla^{\mathfrak s}_{X_x} e = \frac{1}{2} \pi \Bigl(e_{*_x} X_x, (\pmb{-} {\mathfrak s}(x) \cdot e_{*_x}X_x) \Bigr) = 0$$
if and only if
$$m_{-1*} \Bigl({\mathfrak s}(x) \cdot e_{*_x} X_x \Bigr) = - e_{*_x} X_x,$$
or equivalently, if and only if 
$$\Theta_{(b,\d)} (e_{*_x} X_x) = -e_{*_x} X_x.$$
Since $\nabla^{\mathfrak s}_{X_x} e = 0$ characterizes vectors $V_e$ in $\h^{\mathfrak s}$, the claim is proven.
\cqfd

\section{Levi-Civita connection}

Given a pseudo-Riemannian metric $g$ on a manifold $M$, it is well-known that there exists a unique torsionless affine connection $\nabla$ on $M$ that preserves the metric $g$ in the sense that $\nabla g = 0$. Let us reprove this fact in terms of symmetry jets. \\

Consider the vector bundle $\pi : S^2(M) \to M$, consisting of covariant symmetric $2$-tensors. The groupoid $\b^{(1)}(\PP(M))$ naturally acts on $S^2(M)$~:
$$\rho : \b^{(1)}(\PP(M)) \times_{(\alpha, \pi)} S^2(M) \to S^2(M) : (\xi, h) \mapsto \xi \cdot h,$$
with $(\xi \cdot h) (X_y, Y_y) = h (\xi^{-1} X_y, \xi^{-1} Y_y)$ for $X_y, Y_y \in T_yM$, $y = \beta(\xi)$. \\

Now, any (holonomic or not) symmetry jet ${\mathfrak s} : M \to \b^{(1,1)}(\PP(M))$ induces a horizontal distribution $\d = \d^\fs$ along the bisection $b = -I$ and an involutive vector bundle map $\Psi_{(b, \d)}$ (cf.~(\ref{Psi}))~:

$$\Psi_{(b, \d)} : TS^2(M) \to TS^2(M) : X_h \mapsto \rho_{*_{(-I_x, h)}}\Bigl(\ol{\pi_*(X_h)}^{(\d_{-I_x}, \alpha_*)}, X_h \Bigr).$$
covering the identity map $S^2(M) \to S^2(M)$. Notice that in this case $\Psi_{(b, Tb)} = \id$, hence $\Psi_{(b, \d)} = \Theta_{(b, \d)}$. Any leaf $L$ of the horizontal distribution 
$$(E_\Psi^{-1})_{h_x} = \im (- {\mathfrak s}(x) + I)$$
by $(-1)$-eigenspaces of $\Psi_{(b, \d)}$ such that $\pi|_L : L \to M$ is a diffeomorphism is a parallel symmetric $2$-tensor. 

\begin{rmk}
Notice that a pseudo-Riemannian metric $h$ on $M$ --- or for that matter any covariant tensor on $M$ --- determines the subgroupoid $\OO^{(1)}$ of $\b^{(1)}(\PP(M))$ consisting of $1$-jets that preserve $h$. It contains $-I$ as subgroupoid (provided $h$ is a $2p$-tensor). Since $\OO^{(1)}$ is closed, it is necessarily a Lie subgroupoid (cf.~\tref{closed-subgr}). Denotes by $\OO^{(1,1)}(M)$ (\rp $\OO^{(1,1)}_h(M)$) the groupoid of $1$-jets of local bisections of $\OO^{(1)}$ (\rp $\OO^{(1,1)}(M) \cap \b^{(1,1)}(\PP(M))$). Now if a symmetry jet $\fs$ on $M$ takes its values in $\OO^{(1,1)}_h(M)$, then $h$ is parallel for the connection $\nabla^\fs$. Indeed, the map $\Psi_{b, \d}$ obviously preserves the tangent spaces to the section $h$ of $S^2(M) \to M$. Conversely, if a symmetry jet $\fs$ is such that $h$ is parallel for $\nabla^\fs$, then $\fs(M) \subset \OO^{(1)}(M)$. The problem is to understand whether $\OO^{(1)}(M)$ supports a unique distribution along $-I$ which lies in $\e \cap T\OO^{(1)}(M)$. The proof of the next proposition goes along a different path. 
\end{rmk}

\begin{prop}\label{LC}
Given a pseudo-Riemannian metric $h$, there is a unique holonomic symmetry jet ${\mathfrak s}$ such that the $(-1)$-eigenspace of $\Psi_{(-I, \d^\fs)}$ in $T_{h_x}S^2(M)$ coincides with $h_{*_x}(T_xM)$ at all points $x$ in $M$.
\end{prop}

\begin{rmk} A relatively convincing argument, though not a complete proof, in favor of the previous statement is that the dimension of the manifold of $2$-jets  over $-I_x$ and the dimension of the space of $1$-jets of symmetric $2$-tensors over $h_x$ coincide, a fact that is specific to pseudo-Riemannian metrics (as opposed to symplectic structures for instance, see \rref{symplectic} below). 
\end{rmk}

\Pf Fix a point $x$ in $M$. Let ${\mathfrak s}^0(x) = j^2_x s^0_x, {\mathfrak s}(x) = j^2_xs_x$ be two $2$-jets extending $-I_x$. Then ${\mathfrak s}(x) - {\mathfrak s}^0(x)$ may be thought of as a bilinear map $T_xM \times T_xM \to T_xM$, or equivalently a linear map $A : T_xM \to \End(T_xM, T_xM)$ which is symmetric since both ${\mathfrak s}_0(x)$ and ${\mathfrak s}(x)$ belong to $\b^{(2)}(\PP(M))$ (cf.~\rref{affine-str}). \\

Of course a $2$-jet $\xi$ extending $-I_x$ induces an involution $\Psi_\xi$ of $T_{h_x}S^2(M)$ whose $(-1)$-eigenspace is a horizontal $n$-plane $(E_\Psi^{-1})_{h_x}$. The latter corresponds to the $1$-jet $\T(\xi) = j^1_xh$ of some local section $h$ of $S^2(M)$ extending $h_x$ through $(E_\Psi^{-1})_{h_x} = h_{*_x} (T_xM)$. Our aim is to show that the map $\T$ is a bijective correspondence as this implies that if $h$ is a pseudo-Riemannian metric on $M$, then $\fs(x) = \T^{-1}(j^1_xh)$ defines a holonomic symmetry jet for which $h$ is parallel. \\

Let $j^1_xh^0 = \T(\fs^0(x))$ and $j^1_x h = \T(\fs(x))$. Then for any $X_x \in T_xM$, the vector $h_{*_x} (X_x) - h^0_{*_x} (X_x)$ is canonically identified to an element, denoted $h(X_x)$, of the fiber of $\pi$ over $x$, that is to a symmetric tensor. More precisely, for any lift $Y_{h_x}$ of $X_x$ in $T_{h_x}S^2(M)$, we have
$$\begin{array}{ccl}
h(X_x) & = & \frac{1}{2}\Bigl(- {\mathfrak s}(x) + I \Bigr) \cdot Y_{h_x} - \frac{1}{2}\Bigl(- {\mathfrak s}^0(x) + I \Bigr) \cdot Y_{h_x} \\
& = & - \frac{1}{2}\Bigl({\mathfrak s}(x) \cdot Y_{h_x} - {\mathfrak s}^0(x) \cdot Y_{h_x} \Bigr) \\
& = & - \frac{1}{2} \Bigl(\rho_{*_{(-I_x, h_x)}} \bigl((s_x)_{*_{x}}(X_x), Y_{h_x} \bigr) - \rho_{*_{(-I_x, h_x)}}\bigl((s^0_x)_{*_{x}}(X_x),  Y_{h_x}\bigr)\Bigr)\\
& = & - \frac{1}{2}\rho_{*_{(-I_x, h_x)}} \Bigl((s_x)_{*_{x}}(X_x) - (s^0_x)_{*_{x}}(X_x),  Y_{h_x} - Y_{h_x}\Bigr)\\
& = & - \frac{1}{2}\rho_{*_{(-I_x, h_x)}} \Bigl(A(X_x), 0_{h_x}\Bigr) \\
& = & - \frac{1}{2} \Bigl(h_x (A(X_x) \cdot , \cdot ) + h_x ( \cdot, A(X_x) \cdot) \Bigr).
\end{array}$$
This shows in particular that for $\fs^0(x)$ and $X_x$ fixed, the correspondence $A(X_x) \mapsto h(X_x)$ is a linear map between vector spaces of the same dimension. When $h_x$ is non-degenerate and $A$ is symmetric it is also injective as we prove now. Suppose that for all $Y_x, Z_x \in T_xM$, we have 
$$h_x (A(X_x) Y_x, Z_x) + h_x (Y_x, A(X_x) Z_x) = 0.$$
In other word the $h_x$-adjoint of $A(X_x)$ is $-A(X_x)$. Now, let us suppose that $A$ is symmetric. Then 
$$\begin{array}{cl}
h_x (A(X_x)Y_x, Z_x) & = - h_x(Y_x, A(X_x) Z_x) = - h_x(Y_x, A(Z_x) X_x) \\
& = h_x(A(Z_x) Y_x, X_x) = h_x(A(Y_x) Z_x, X_x) \\
& = - h_x(Z_x, A(Y_x) X_x) = - h_x(Z_x, A(X_x) Y_x),
\end{array}$$
which implies, by symmetry of $h_x$ that $A$ vanishes. 
\cqfd

\begin{rmk}\label{symplectic} Let $\omega$ be a symplectic structure on $M$. Then a similar construction yields a correspondence between $2$-jets extending $-I_x$ and $1$-jets of closed $2$-forms extending a fixed one $\omega_x$. On the one hand, the dimension of the space of $2$-jets in $p^{-1}(-I_x)$ is $n^2(n+1)/2$. On the other hand, the dimension $d$ the space of $1$-jets of closed forms, that is $j^1_x\omega$ such that $(d\omega)_x = 0$ coincides with the difference between the dimension of the space of all jets at $x$ of $2$-forms, that is $n^2(n-1)/2$, and the dimension of the space of $3$-forms at $x$, that is $n(n-1)(n-2)/6$. Hence $d = n(n+1)(n+2)/6$, which is precisely the dimension of the space of totally symmetric $3$-tensors at $x$, whose we know parameterize the space of symplectic torsionless affine connections \cite{T}. 
\end{rmk}

\section{Relation with Lie algebroid connections}\label{groupoid} 

It is known that affine connections and Lie algebroid connections on the Lie algebroid of the general linear groupoid $GL(TM) = \b^{(1)}(\PP(M))$ are the same thing. The point of this section is to show a short path between a symmetry jet $\fs$ and a Lie algebroid connection on the Lie algebroid of $\b^{(1)}(\PP(M))$. The former may be thought of as the rank $n = \dim M$ distribution $\d^\fs$ along the bisection $-I$ in $\b^{(1)}(\PP(M))$, the latter is in fact a rank $n$ distribution on $\b^{(1)}(\PP(M))$ tangent to the $\beta$-fibers. Going from one to the other is done by elementary manipulations on distributions.\\

We briefly recall Atiyah's sequence which is another way to describe a linear connection and which yields the motivations for the definition of a Lie algebroid connection. We refer to \cite{dSW} for an elegant introduction to this topic. A lie algebroid connection is then compared in the case of $\b^{(1)}(\PP(M))$ to a symmetry jet.  \\

To a principal fiber bundle $\pi : P \to M$ with structure group $G$ is naturally associated a Lie groupoid, called the {\sl gauge groupoid}, as described hereafter. Its set of objects is defined to be $G_P = (P \times P)/G$, the quotient of $P\times P$ through the diagonal action of $G$. Equivalently, $G_P$ may be described as the set of equivariant maps between any two fibers of $P$. The image by the source $\alpha$ (respectively target $\beta$) of an equivariant map $\xi : P_x \to P_y$ is $x$ (respectively $y$). Finally, the product is the composition of maps. Notice that for the principal bundle $\pi : \F(M) \to M$ of frames of $TM$, the gauge groupoid is precisely $\b^{(1)}(\PP(M))$. \\

The Lie algebroid (cf.~\cite{Mck-05}) of the Lie groupoid $G = \b^{(1)}(\PP(M))$ coincides with the set $T\F(M)/G$ of $G$-invariant vector fields on $\F(M)$ with anchor induced by the projection $\pi_* : T\F(M) \to TM$. 
In this case the anchor is surjective. Lie algebroids with surjective anchor are called {\sl transitive Lie algebroids}. \\

Now a principal connection on $\pi : \F(M) \to M$ is an invariant horizontal distribution on $\F(M)$. Equivalently, it is a section of the anchor $\rho : L \to TM$. 

\begin{dfn}
A {\sl connection on a transitive Lie algebroid} is defined to be a section of its anchor or a splitting of the exact sequence of vector bundles

$$0 \to \Ker \rho \to L \to TM \to 0.$$
The induced map $L \to \Ker \rho$ is called a {\sl connection form}. 
\end{dfn}

To rely the Lie algebroid point of view to the the symmetry jet point of view on a linear connection, we will describe below how a Lie algebroid connection on the Lie algebroid of $\b^{(1)}(\PP(M))$ directly induces a symmetry jet. \\

Let $\sigma^\beta : TM \to T_{\varepsilon(M)}G^\beta$ be a Lie algebroid connection. The right-invariant distribution on $G$ associated to the image of $\sigma^\beta$ is denoted hereafter by $\Sigma^\beta$. Define a horizontal distribution on $G$ along $M$ by
$$\EN^\sigma = \{\sigma (X) + F(\sigma(X)); X \in TM\},$$
where $F : TG^\beta \to TG^\alpha$ is the flip automorphism, defined by
\begin{equation}\label{flip}
F : T_xG^\beta \to T_xG^\alpha : Y \to Y - \alpha_{*}(Y)
\end{equation} that exchanges the tangent space to the $\beta$-fibers over a point in $M$ with that of the corresponding $\alpha$-fiber. It is also induced by the canonical isomorphisms $T_xG^\beta \simeq T_xG/T_xM \simeq T_xG^\alpha$ with the normal bundle to $TM$. Notice that $\fb(\EN^\sigma_{x}) = -I_x$ (the \emph{bouncing} map $\fb$ is defined in \rref{bouncing}). Let $b$ denote the bisection $-I$ of $\b^{(1)}(\PP(M))$ and let $L_{b}$ denote the diffeomorphism of $G$ induced by the left action of $b$~:
$$L_{b} : G \to G : \xi \mapsto b(\beta(\xi)) \cdot \xi.$$
Then 
$$\d^\sigma \stackrel{\Def}{=} (L_{b})_*(\EN^\sigma)$$ 
is a horizontal distribution along $b$ that is contained in the holonomic distribution $\e^{(1)}$ (cf.~\dref{def-e} and \nref{nota-e}) because $\fb(\d^\sigma_{-I_x}) = \fb(\EN^\sigma_{\varepsilon(x)}) = -I_x$ or equivalently, a symmetry jet. \\ 

Conversely, let $\d \subset \e$ be a horizontal distribution on $G = \b^{(1)}(\PP(M))$ along the bisection $b = -I$, then $\d$ induces a Lie algebroid connection $\sigma$ as follows. First translate $\d$ to a horizontal distribution $\EN^\d$ along $\varepsilon(M)$ via $L_{b}$~:
$$\EN^\d = (L_{b})_*(\d).$$
Notice that $\fb(\EN^\d_{\varepsilon(x)}) = -I_x$. Then define
$$\sigma^\d : TM \to TG^\beta : X \mapsto \frac{1}{2}\Bigl((\alpha_{*_x}\Big|_{\EN^\d})^{-1}(X) + X \Bigr).$$
One verifies easily that $\sigma$ is indeed a Lie algebroid connection~:
$$\begin{array}{l}
(\alpha_* \circ \sigma^\d)(X) =  \frac{1}{2}\alpha_* \Bigl((\alpha_{*_x}\Big|_{\EN^\d})^{-1}(X) + X\Bigr) = \frac{1}{2} (X + \alpha_{*_x}(X)) = X \\
(\beta_* \circ \sigma^\d)(X) =  \frac{1}{2}\beta_* \Bigl((\alpha_{*_x}\Big|_{\EN^\d})^{-1}(X) + X\Bigr) =  \frac{1}{2}(-X + \beta_{*_x}(X)) = 0
\end{array}$$

Finally, one can verify that a Lie algebroid connection on $\b^{(1)}(\PP(M))$ and its corresponding horizontal distribution along $-I$ are two realizations of a same linear connection. \\

To go directly from a Lie algebroid connection $\sigma^\beta$ on $G = \b^{(1)}(\PP(M))$ to the corresponding global distribution $\d$ on $\b^{(1)}(\PP(M))$, one can proceed as follows. Consider the map 
$$\varphi_{\sigma^\beta} : \alpha^*TM \to T\b^{(1)}(\PP(M))$$
\begin{equation}\label{s-->d}
\varphi_{\sigma^\beta} (\xi, X) = (L_{\xi})_{*_x} \circ \sigma^\beta (X) - (R_{\xi})_{*_y} \circ F \circ \sigma^\beta \circ \xi (X),
\end{equation}
The image of $\varphi_{\sigma^\beta}$ is an invariant horizontal distribution (see Figure 5). Let us check that the image of $\varphi_{\sigma^\beta}$ indeed consists of eigenvectors for the eigenvalue $-1$ of the map $\psi_\xi$, $\xi \in \b^{(1)}(\PP(M))$ defined in terms of $\d$ along $b = -I$ by (\ref{psi}) in \sref{UAE}. Let $X \in T_xM$, then $\sigma(X) = 1/2 (X + \ol{X})$, where $\ol{X}$ stands for the lift with respect to $\alpha_*$ of $X$ in $\EN_{\varepsilon(x)}$. Then $F(\sigma(X)) = 1/2(-X + \ol{X})$. Thus

$$\begin{array}{cll}
\psi_\xi \Bigl(\varphi_{\sigma^\beta} (\xi, X) \Bigr) & = & \ol{\xi X} \cdot \Bigl[(L_{\xi})_{*_x} \bigl(\sigma(X)\bigr) - (R_{\xi})_{*_y} \bigl(F(\sigma(\xi X))\bigr) \Bigr] \cdot -\ol{X} \\
& = & \Bigl( (L_{\xi})_{*_x} \bigl(\sigma(X)\bigr) \Bigr) \cdot -\ol{X} + \ol{\xi X} \cdot \Bigl( - (R_{\xi})_{*_y} \bigl(F(\sigma(\xi X))\bigr) \Bigr)\\
& = & 0_\xi \cdot \Bigl(\frac{1}{2} (X + \ol{X})\Bigr) \cdot \Bigl(\frac{1}{2} ( -\ol{X} - \ol{X})\Bigr) \\
&& + \Bigl(\frac{1}{2}(\ol{\xi X} + \ol{\xi X}) \Bigr) \cdot \Bigl(\frac{1}{2} (\xi X - \ol{\xi X})\Bigr) \cdot 0_\xi \\
& = & 0_\xi \cdot \Bigl(\frac{1}{2} (-\ol{X} - X)\Bigr) + \Bigl(\frac{1}{2}(\ol{\xi X} - \ol{\xi X}) \Bigr) \cdot 0_\xi \\
& = & \Bigl.- \varphi_{\sigma^\beta} (\xi, X)
\end{array}$$

\begin{figure}[h!]\label{sigmad}
\begin{center}
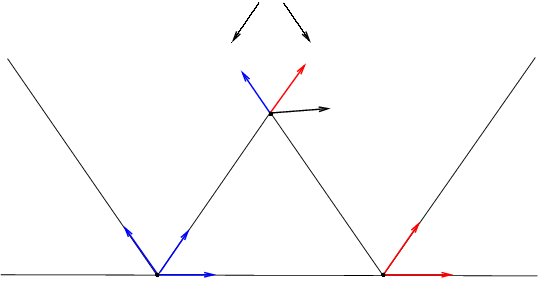 
\caption{The distribution $\d$ on $\b^{(1)}$ from the map $\sigma$}
\end{center}
\end{figure}

\begin{rmk} For a general groupoid, the data of a Lie algebroid connection is not equivalent to that of an invariant horizontal distribution. Indeed, consider for instance the pair groupoid $\PP(M) = M \times M \rightrightarrows M$ on the manifold $M$. The anchor $\rho : T_xG^\beta \to T_xM$ of the associated Lie algebroid being an isomorphism, its inverse yields a canonical Lie algebroid connection. In contrast, the data of an invariant horizontal distribution $\d$ on $M \times M$ is equivalent to that of a groupoid morphism $M \times M \to \b^{(1)}(\PP(M))$ over the identity or, in other words, a global trivialization of $TM$. This implies that the manifold $M$ is parallelizable, which is a noticeable restriction on the type of pair groupoids that admit an invariant horizontal distribution. \\

In general, the groupoid morphism $G \to \b^{(1)}(\PP(M))$ over the identity is the data that is missing in a Lie algebroid connection. More precisely, a Lie algebroid connection endowed with a Lie groupoid morphism $\psi : G \to \b^{(1)}(\PP(M))$ over the identity $M \to M$ is an invariant horizontal distribution on $G$. Indeed, the formula (\ref{s-->d}) makes now sense if $\xi$ acting on $T_xM$ is replaced by $\psi (\xi)$.
\end{rmk}

\section{Geodesics and parallel transport} 

Let ${\mathfrak s}$ be a symmetry jet, let $S$ be the corresponding affine extension section and $\d^{\mathfrak s}$ the associated invariant horizontal distribution on $\b^{(1)}(\PP(M))$. Let also $\sigma^\beta : TM \to T\b^{(1)}(\PP(M))^\beta$ denote the corresponding Lie algebroid connection and $\sigma^\alpha = \iota \circ \sigma^\beta$ its dual Lie algebroid connection. The right-invariant distribution on $\b^{(1)}(\PP(M))$ associated to the image of $\sigma^\alpha$ (\rp $\sigma^\beta$) is denoted hereafter by $\Sigma^\alpha$ (\rp $\Sigma^\beta$). \\

We consider smooth paths $\gamma : I \to M$ defined on some interval $I$ containing $0$ that are in addition regular, that is $d \gamma /dt$ does not vanish. Let $\T$ be some distribution defined on $\b^{(1)}(\PP(M))$ for which $\beta_*$ restricts to a fiberwise linear isomorphism. Then an $\beta$-lift of $\gamma$ through $\xi \in \b^{(1)}(\PP(M))$ is a path 
$$L(\gamma) = L_\xi^\T(\gamma) : I' \to \b^{(1)}(\PP(M))$$ defined on an eventually smaller interval $I' \subset I$, still containing $0$, such that

\begin{enumerate}
\item[-] $\alpha \circ L(\gamma) = \gamma$,
\item[-] $\displaystyle{\frac{dL(\gamma)}{dt}(t) \in \T_{L(\gamma)(t)}}$ for all $t \in I'$,
\item[-] $L(\gamma)(0) = \xi$.
\end{enumerate}
Such a lift is the flow line through $\xi$ of the vector field obtained by lifting the velocity vector field of $\gamma$ (restricted to a subinterval where $\gamma$ is injective) to a vector field tangent to $\T$ and defined on $\alpha^{-1}(\gamma(I))$, so it certainly does exist, although not in general on the entire interval $I$. \\

Given an affine connection $\nabla$, a common way to think about the induced parallel transport is as follows. Given a path $\gamma : 0 \in I \to M$, denote by $\ol{\gamma}_X(t)$ the lift of $\gamma$ through the vector $X$ in $T_xM$ tangent to the horizontal distribution $\h$ on $TM$. The parallel transport along $\gamma$ is then defined to be 
$$\tau^\gamma(t) : T_{\gamma(0)}M \to T_{\gamma(t)}M : X \mapsto \ol{\gamma}_X(t).$$

\begin{lem}\label{parallel-transport} Let $\gamma : I \to M$ be regular path in $M$ with $\gamma(0) = x$. Then the lift $L^{\Sigma^\alpha}_{\varepsilon(x)} (\gamma)$ is the parallel transport along $\gamma$. In other words,
$$L^{\Sigma^\alpha}_{\varepsilon(x)}(\gamma)(t) = \tau^\gamma(t).$$
\end{lem}

\Pf Recall from \sref{hor-ass-bdle} the description of $\h$ as eigenspace~:
$$\h = \Ker(\Theta_{(b,\d)} + I) = \im (- \Theta_{(b,\d)} + I)$$
and notice that 
$$\Theta_{(b,\d)} = \Psi_{(b, Tb)} \circ \Psi_{(b, \d)} = \Psi_{(I, \EN)},$$
since $Tb \cdot \d = \EN$. In other terms, the lift of a vector $Y \in T_xM$ to a vector $\ol{Y}$ tangent to $\h$ at $X$ admits the following expression~:
$$\begin{array}{ccl}
\ol{Y} & = & \frac{1}{2}\Bigl(- \Theta_{(b, \d)}(\ol{Y}) + \ol{Y} \Bigr) \\
& = & \frac{1}{2}\Bigl(- \rho^{(1)}_{*_{}} \bigl(\ol{Y}^{\EN_{\varepsilon(x)}, \alpha_*}, \ol{Y} \bigr) + \ol{Y} \Bigr) \\
& = &\frac{1}{2}\Bigl(\rho^{(1)}_{*_{}} \bigl( - \ol{Y}^{\EN_{\varepsilon(x)}, \alpha_*}, - \ol{Y} \bigr) + \rho^{(1)}_{*_{}} \bigl( I_* Y, \ol{Y} \bigr) \Bigr) \\
& = & \frac{1}{2} \rho^{(1)}_{*_{}} \bigl( - \ol{Y}^{\EN_{\varepsilon(x)}, \alpha_*} +  I_* Y, 0_X \bigr) \Bigl.\\
& = & \rho^{(1)}_*\bigl( \sigma^\alpha(Y), 0_X \bigr). \Bigl.
\end{array}$$
This implies that the tangent vector to the path $L^{\Sigma^\alpha}_{\varepsilon(x)}(\gamma)(t) \cdot X$ belongs to~$\h$. Indeed, set 
$$\tau(t) = L^{\Sigma^\alpha}_{\varepsilon(x)}(\gamma)(t).$$ 
Then
$$\begin{array}{ccl}
\displaystyle{\frac{d}{dt} \rho^{(1)}\Bigl(\tau(t), X\Bigr)\Bigr|_{t}} & = & \displaystyle{\rho^{(1)}_*\Bigl( \frac{d}{dt}\tau(t)\Bigr|_{t}, 0_X \Bigr)} \\
& = & \displaystyle{\rho^{(1)}_*\Bigl( (R_{\tau(t)})_{*_{\gamma(t)}} \sigma^\alpha \bigl(\frac{d \gamma}{dt}\bigr), 0_X \Bigr)} \\
& = & \displaystyle{\rho^{(1)}_*\Bigl( m_*\bigl\{\sigma^\alpha (\frac{d \gamma}{dt}), 0_{\tau(t)}\bigr\}, 0_X \Bigr)} \\
& = & \displaystyle{\rho^{(1)}_*\Bigl( \sigma^\alpha (\frac{d \gamma}{dt}), \rho^{(1)}_* \bigl(0_{\tau(t)}, 0_X \bigr) \Bigr)} \quad \mbox{(since $\rho$ is an action)} \\
& = & \displaystyle{\rho^{(1)}_*\Bigl( \sigma^\alpha \bigl(\dot{\gamma}(t)\bigr), 0_{\tau(t) \cdot X} \Bigr)} \in \h. 
\end{array}$$
Observe that in the case of the distributions $\Sigma^\alpha$ (and $\Sigma^\beta$), the lift of a path is defined on the entire interval $I$. This is due to the right-invariance of $\Sigma^\alpha$
\cqfd

\begin{figure}[h!]\label{geodesics-fig}
\begin{center}
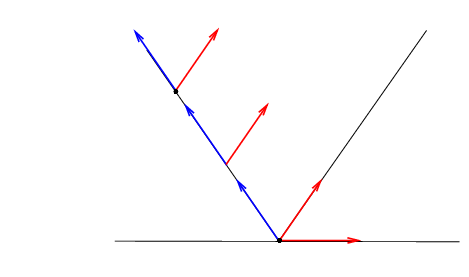 
\caption{The distribution $\d$ on $\b^{(1)}$ from the map $\sigma^\beta$}
\end{center}
\end{figure}

Consider now a path $\gamma : I \to M$, with $\gamma(0) = x$. Then 
$$v(t) = \Bigr(L^{\Sigma^\beta}_{\varepsilon(x)}(\gamma)(t)\Bigr) \cdot \frac{d\gamma}{dt}\Bigr|_{t}$$
is a path in $T_xM$, the inverse parallel transport along $\gamma$ of $\frac{d \gamma}{dt}$. Conversely, given a path $v : I \to T_xM$, it is associated to a unique path $\gamma = {\mathcal I}(v)$ (${\mathcal I}$ for ``integration'') in $M$ such that the inverse parallel transport of $\frac{d \gamma}{dt}$ along $\gamma$ yields $v$. 
The path $\gamma$ is constructed as follows. First consider the path
$$\sigma^\beta \circ v :  I \to \Sigma^\beta_{\varepsilon(x)}.$$ 
It can be left-translated along the $\alpha$-fiber of $x$ into a path of vector fields $X_\xi = (L_\xi)_*(\sigma^\beta(v(t)))$ in $\Sigma^\beta$ along $\alpha^{-1}(x)$ which can be flipped with respect to $\d^\fs$ and become a path of vector fields $F(X_\xi)$ in $\Sigma^\alpha$ tangent to $\alpha^{-1}(x)$ (see \fgref{geodesics-fig}). Indeed, it follows from (\ref{s-->d}) that the flip of $\Sigma^\beta$ with respect to $\d^\fs$ is $\Sigma^\alpha$. The flow line of the latter through $x$, denoted $\ol{\gamma}(t)$, $t \in I' \subset I$ is the parallel transport along the path ${\mathcal I}(v)(t) = \gamma(t) = \beta \circ \ol{\gamma}(t)$. \\ 

These observations allow us to characterize geodesics in terms of $\sigma$, hence, indirectly, in terms of ${\mathfrak s}$.

\begin{prop}\label{geodesics} A path $\gamma : I \to M$ is a geodesic if and only if its associated path $\ul{\gamma} : I \to T_{\gamma(0)}M$ is constant. In particular, for any $X \in TM$, the geodesic $\gamma_X$ tangent to $X$ at $0$ is the path associated to the constant path $\ul{\gamma}_X(t) = X$.
\end{prop}
Notice that the homogeneity property $\gamma_{cX} (t)= \gamma_X(ct)$ appears clearly in the construction. \\

\Pf This follows readily from \pref{parallel-transport} and the characterization of geodesics as paths with constant velocity.
\cqfd

\lref{parallel-transport} also yields an interpretation of the lift with respect to $\d^\fs$ of a path in $M$.
 
\begin{lem}\label{d-lift}
Let $\gamma : I \to M$ be a path in $M$ and $\xi : T_{\gamma(0)} M \to T_yM$ a linear map. Then the $\alpha$-lift of $\gamma$ tangent to $\d$ and passing through $\xi$ is the path~:
\begin{equation}\label{d-lift&parallel-t}
L_\xi^\d(\gamma)(t) : T_{\gamma(t)}M \to T_{\gamma'(t)}M : X \mapsto \tau^{\gamma'}(t) \circ \xi \circ (\tau^{\gamma}(t))^{-1}(X),
\end{equation} 
where 
$$\gamma' = {\mathcal I}\bigl(\xi \circ (\tau^\gamma(t))^{-1} (\frac{d\gamma}{dt})\bigr).$$ 
\end{lem}

\Pf Setting $\widetilde{\gamma} = L_\xi^\d(\gamma)$, the relation to be proven is 
$$\widetilde{\gamma}(t) = \tau^{\gamma'}(t) \circ \xi \circ (\tau^{\gamma}(t))^{-1},$$
or, in other terms,
\begin{equation}\label{d-lift-equ}
\widetilde{\gamma}(t) \circ \tau^\gamma (t) = \tau^{\gamma'}(t) \circ \xi.
\end{equation}
Equivalently, the path $\widetilde{\gamma} \circ \tau^\gamma$ is everywhere tangent to the distribution $\Sigma^\alpha$. Fix $t \in I$ and set
\begin{enumerate}
\item[-] $X = \frac{d \gamma}{dt}\bigl|_t$,
\item[-] $\ol{X}_1 = \frac{d \widetilde{\gamma}}{dt}\bigl|_t$,
\item[-] $\ol{X}_2 = \frac{d \tau^\gamma}{dt}\bigl|_t$.
\end{enumerate}
The relation (\ref{d-lift-equ}) becomes
$$\ol{X}_1 \cdot \ol{X}_2 \in \Sigma^\alpha_{\widetilde{\gamma}(t) \circ \tau^\gamma (t)},$$
for all $t \in I$, or equivalently,
$$\bigl(R_{\widetilde{\gamma}(t) \cdot \tau^\gamma(t)}^{-1}\bigr)_{*} \bigl(\ol{X}_1 \cdot \ol{X}_2\bigr) \in \im \sigma^\alpha.$$
Remembering the characterization of $\d^\fs$ as the distribution consisting of the traces in $\e$ of the $(-1)$-eigenspaces of the map $\psi$ (cf.~(\ref{psi})), we know that
$$L \cdot \ol{X}_1 \cdot R = -\ol{X}_1,$$
where $L$ (\rp $R$) is the $\alpha$-lift (\rp $\beta$-lift) of $\beta_*\ol{X}_1 = \widetilde{\gamma}(t)(X)$ (\rp $\alpha_*\ol{X}_1 = X$) in $\d_{\widetilde{\gamma}(t)}$. Equivalently, 
$$\ol{X}_1 \cdot R = \iota_*(L) \cdot -\ol{X}_1,$$
and since $\iota_*(L) = -L$,
$$\ol{X}_1 \cdot R = - \bigl(L \cdot \ol{X}_1\bigr).$$
Moreover, because $\ol{X}_2$ belongs to the right invariant distribution $\Sigma^\alpha$, 
$$\ol{X}_2 = (R_{\tau^\gamma(t)})_{*}\sigma^\alpha(X) = \bigl(\frac{1}{2} X + \frac{1}{2} (m_{-1})_* R\bigr) \cdot 0_{\tau^\gamma(t)}.$$ 
Therefore, using the linearity of the differential of the multiplication $m$ in $\b^{(1)}(\PP(M))$, we compute
$$\begin{array}{lll}
\bigl(R_{\widetilde{\gamma}(t) \cdot \tau^\gamma(t)}^{-1}\bigr)_{*} \bigl(\ol{X}_1 \cdot \ol{X}_2\bigr) & = & \ol{X}_1 \cdot \ol{X}_2 \cdot 0_{(\widetilde{\gamma}(t) \cdot \tau^\gamma(t))^{-1}} \\
& = & \ol{X}_1 \cdot \bigl(\frac{1}{2} X + \frac{1}{2} (m_{-1})_* R\bigr) \cdot 0_{\tau^\gamma(t)} \cdot 0_{(\widetilde{\gamma}(t) \cdot \tau^\gamma(t))^{-1}} \Bigl. \\
& = & \bigl( \frac{1}{2}\ol{X}_1 + \frac{1}{2} \ol{X}_1 \bigr) \cdot \bigl(\frac{1}{2} X + \frac{1}{2} (m_{-1})_* R\bigr) \cdot 0_{(\widetilde{\gamma}(t))^{-1}} \Bigl.\\
& = & \Bigl( \frac{1}{2}\bigl(\ol{X}_1\cdot X \bigr) + \frac{1}{2} \bigl( \ol{X}_1 \cdot (m_{-1})_* R\bigr) \Bigr)  \cdot 0_{(\widetilde{\gamma}(t))^{-1}} \Bigl.\\
& = &  \frac{1}{2}\Bigl(\ol{X}_1 + (m_{-1})_* (\ol{X}_1\cdot R) \Bigr)  \cdot 0_{(\widetilde{\gamma}(t))^{-1}}\\
& = &  \frac{1}{2}\Bigl(\ol{X}_1 - (m_{-1})_* (L \cdot \ol{X}_1) \Bigr)  \cdot \frac{1}{2} \Bigl( \iota_*\ol{X}_1 - \iota_*\ol{X}_1 \Bigr)\\
& = &  \frac{1}{2}\Bigl(\ol{X}_1  \cdot \iota_*\ol{X}_1 \Bigr) - \frac{1}{2} \Bigl((m_{-1})_* (L \cdot \ol{X}_1) \cdot \iota_*\ol{X}_1 \Bigr)\\
& = & \frac{1}{2} \widetilde{\gamma}(t)(X) - \frac{1}{2} (m_{-1})_*L \Bigl.\\
& = & \sigma^\alpha\bigl(\widetilde{\gamma}(t)(X)\bigr)\Bigl.
\end{array}$$
\cqfd

\section{Geodesic symmetries and locally symmetric spaces}\label{geod-sym&loc-sym-sp}

Let $\fs$ be a symmetry jet on $M$, we would like to describe the geodesic symmetries of the affine connexion $\nabla^\fs$ directly in terms of $\d^\fs$. More generally, for each $\xi \in \b^{(1)}(\PP(M))$, we construct hereafter the extension of $\xi$ through geodesics~:
$$\varphi_\xi : U_{\alpha(\xi)} \to U_{\beta(\xi)} : \exp_{\beta(\xi)}  \circ \; \xi \circ \exp_{\alpha(\xi)}^{-1},$$ where $U_{\alpha(\xi)}$ and $U_{\beta(\xi)}$ are some neighborhoods of $\alpha(\xi)$ and $\beta(\xi)$ respectively. Of course, the geodesics symmetries are the maps $s_x = \varphi_{-I_x}$, $x \in M$. We would like to show that the maps $\varphi_\xi$ are the best candidates affine transformations, that is, when $S^{k \cdot (1)}(\xi)$ is affine, it coincides with $j^k_{\alpha(\xi)}\varphi_\xi$. \\

Let $\xi$ be a fixed element in $\b^{(1)}(\PP(M))$. If $\gamma : I \to M$ is a regular path passing through $\alpha(\xi)$, let $L(\gamma) = L_\xi^{\d} (\gamma)$ be the $\alpha$-lift of $\gamma$ tangent to $\d = \d^\fs$ and passing through $\xi$. It is defined on some subinterval $I' \subset I$. Given $X \in TM$, let $\gamma_X$ denote the maximal geodesic through $X$. The homogeneity property for geodesics implies that the image of  $L(\gamma_X)$ is independent on $X$ in ${\mathbb R}X$. Define $b_\xi$ to be the union over all vectors $X$ in the unit sphere $ST_xM$ of $T_xM$ (relative to some Riemannian metric on $M$) of the images of the paths $L(\gamma_X)$. Near $\xi$, the set $b_\xi$ is an embedded submanifold tangent to $\d^{\mathfrak s}$ at $\xi$, and, in fact, a local bisection. Indeed, the sphere $ST_xM$ being compact, a common domain $I'$ may be chosen for the various lifts $L(\gamma_X)$. The submanifold $b_\xi$ is a best candidate leaf of $\d^\fs$, that is, if there exists a submanifold of $\b^{(1)}(\PP(M))$ that is tangent to $\d^\fs$ up to order $k$, then $b_\xi$ is such a manifold. It follows easily from \pref{geodesics} and \lref{d-lift} that $b_\xi$ is the set of linear maps~:

$$T_{\gamma_X(t)}M \to T_{\varphi_\xi(\gamma_X(t))}M : X' \mapsto \tau^{\gamma_{\xi(X)}}(t) \circ \xi \circ (\tau^{\gamma_X}(t))^{-1}(X'),$$
where $t$ varies in $I'$ and $X$ in $ST_xM$. This is almost $\varphi_\xi$, but not quite yet. More precisely, the base map is $\varphi_\xi$. So to obtain $\varphi_\xi$, we apply the bouncing map 
$\fb$ to $Tb_\xi$ (cf.~\rref{bouncing}), that is, 
$$j^1\varphi_\xi = \fb(Tb_\xi).$$ 
In particular, the geodesic symmetries $s_x$, $x \in M$ are realized as local bisections in $\b^{(1)}(\PP(M))$ as follows~:
$$j^1 s_x = \fb(Tb_{-I_x}).$$
Notice that if $b_\xi$ is tangent to $\d^\fs$ up to order $k$, then $j^1\varphi_\xi$ is tangent to $\d^\fs$ up to order $k-1$. To summarize, we have proven the following proposition.

\begin{prop}\label{geodesic-symmetries} To produce $j^1 \varphi_\xi$, $\alpha$-lift each geodesic through $x$ to a path tangent to $\d^\fs$ passing through $\xi$ then make it a holonomic bisection via $\fb$. In other terms
$$j^1\varphi_\xi = \fb(Tb_\xi).$$
\end{prop}

This procedure is relatively simple. The distribution $\d^\fs$ does not admit $n$-dimensional leaves in general. Nevertheless, it always has $1$-dimensional leaves and gathering the ones passing through a given $\xi$ and that project via $\alpha$ onto the geodesics through $\alpha(x)$ (a natural family of paths filling a neighborhood of $x$) yields a bisection whose projection onto $M \times M$ is the graph of $\varphi_\xi$. In case $\d^\fs$ does admit a leaf $D$ through $\xi$, it will coincide near $\xi$ with $b_\xi$, which becomes automatically a holonomic bisection since $\d^\fs \subset \e$.

\begin{rmk} When $\xi \in \In{\d^\fs} = \b(T, R)$ (cf.~\dref{int-locus}), then $b_\xi$ is osculatory to $\d^\fs$ and when $\xi \in \b(T, \nabla T, R)$, then $j^1\varphi_\xi$ is osculatory to $\d^\fs$. It seems important to make the following distinction. If a bisection $b$ is tangent to $\d^\fs$ up to order $k$, that is $j^{k-1}_{\alpha(\xi)} b = S^{k \cdot (1)}(b(x))$ then the holonomic bisection $\fb(Tb)$ is tangent to $\d^\fs$ only up to order $k-1$, unless $b$ is, in addition, tangent to the holonomic distribution $\e$ up to order $k+1$, in which case the bisection $\fb(Tb)$ is tangent to $\d^\fs$ up to order $k$ as well, that is $S^{k \cdot (1)}(b(x))$ is holonomic. 
\end{rmk}

It is worthwhile noticing that Emmanuel Giroux has shown in \cite{G} that a plane field admits near each point a canonical $2$-jet of {\sl path-osculatory surfaces} (our terminology, to distinguish from osculatory). Path-osculatory, means that each path through the point in the surface is osculatory to the distribution. In the context of the present paper, we obtain for each $\xi \in \b^{(1)}(\PP(M))$, a canonical surface whose second order jet at $x$ coincides with Giroux's osculatory surfaces. \\

Now it is easy to show the well-known result that for a locally symmetric space, that is a space endowed with a torsionless affine connection for which the curvature tensor is parallel, the local diffeomorphisms $\varphi_\xi$ is affine if and only if $\xi \in \b(R)$ (cf.~\cite{H}, Lemma 1.2.~p.~200). In particular, for those spaces --- and only them --- each geodesic symmetry is affine. 

\begin{prop}\label{lss} Let ${\mathfrak s}$ be a symmetry jet whose associated connection $\nabla$ is locally symmetric, that is satisfies $T^\nabla = 0$ and $\nabla R^\nabla = 0$. Then through any $\xi \in \b^{(1)}(\PP(M))$ that preserves the curvature passes a $n$-dimensional leaf of $\d^{\mathfrak s}$.
\end{prop}

\Pf Consider the Lie subgroupoid $\b(R)$. Since $\nabla R = 0$, \lref{d-tangent-to-b(Q)}  implies that $\d^{\mathfrak s}$ is tangent to $\b(R)$ along $\b(R)$. Moreover, \lref{Int-kappa} implies that $\d^{\mathfrak s}$ in involutive along $\b(R)$. Hence the Lie subgroupoid $\b(R)$ is foliated by leaves of~$\d^\fs$.
\cqfd

%
%

\section{Relation with Kobayashi's admissible sections}\label{Kobayashi}

A theorem due to Kobayashi \cite{K} asserts that there is a bijective correspondence between torsionless affine connections and admissible sections of the bundle of $2$-frames $\F^{(2)}(M) \to \F^{(1)}(M)$. Since jets occur in our construction as well, it seems relevant to compare the two approaches. This section contains a brief description of Kobayashi's theorem, interpreted in terms of affine extensions as in \cite{H}, and a direct correspondence between Kobayashi's admissible section and our symmetry jet. \\

Given a manifold $M$ of dimension $n$ and a non-negative integer $k$, {\em the bundle of $k$-frames of $M$}, denoted $\pi^k : \F^{(k)}(M) \to M$ is the set of $k$-jets at $0$ of local diffeomorphisms $\R^n \to M$ defined near $0$, endowed with the canonical projection that sends a jet $j^k_0f$ to $f(0) \in M$. Observe that the bundle of $1$-frames is the usual bundle of frames of $M$. There are obvious maps $\pi^{l \to k} : \F^{(l)}(M) \to \F^{(k)}(M)$ for each pair $l > k$. The group $GL(n, \R)$ acts on the right on each $\F^{(k)}(M)$ through $j^k_0f \cdot A =  j^k_0 (f \circ A)$. A section $\F^{(1)}(M) \to \F^{(2)}(M)$ is said to be admissible if it is equivariant with respect to these $GL(n, \R)$-actions. \\

The main ingredient of this correspondence is the property of uniqueness of affine extension applied to $\R^n$ endowed with the trivial connection $\nabla_o$ and $M$ endowed with a torsionless affine connection $\nabla$ (it is a consequence of \pref{uae} applied to the disconnected affine manifold $(M \cup \R^n, \nabla \cup \nabla_o)$). It ensures existence of an admissible section $s^{\nabla} : \F^{(1)}(M) \to \F^{(2)}(M)$  that maps $\xi : \R^n \to T_{f(0)}M$ to the unique affine $2$-jet that extends it. It is indeed equivariant with respect to the $GL(n, \R)$-actions on $\F^{(1)}(M)$ and $\F^{(2)}(M)$ since $s^{\nabla}(\xi) \cdot A$ is an affine $2$-jet that extends $\xi \circ A$, implying the relation $s^{\nabla}(\xi) \cdot A = s^{\nabla}(\xi \circ A)$.\\

Conversely, given an admissible section $s$, it directly induces a linear horizontal distribution $\h^s$ on $TM$, hence an affine connection $\nabla^s$. Indeed, let $\xi \in \F^{(1)}(M)$. Since $s(\xi) = j^2_0f$ for some local diffeomorphism $f : \R^n \to M$ defined near $0$, one may define for $X \in T_{f(0)}M$
$$\h_X = f_{**_0} (H_{f_{*_0}^{-1}}(X)),$$
where $H_Z$ denotes the natural horizontal plane in $T_ZT\R^n$. The $GL(n, \R)$-invariance of $s^\nabla$ guarantees the independence of $\h$ on the initial choice of a basis $\xi$ in $T_xM$ as well as the linearity of $\h$, in the sense that $\h_{a X + b Y} = m_{a*} (\h_{X}) +_* m_{b*} (\h_Y)$. \\

One can alternatively observe, as is done in \cite{K} that the pullback of the  ${\mathcal Gl}(n, \R)$-component of the canonical form $\theta^{(2)}$ on $\F^{(2)}(M)$ via $s$ yields a connection form on $\F^{(1)}(M)$. \\

Now we have two natural ways to think about torsionless affine connections, one in terms of admissible sections and another one in terms of holonomic symmetry jets. It is tempting to close the triangle and show how to induce a symmetry jet naturally from an admissible section and vice-versa. So doing, we are going to enlarge the Kobayashi's correspondence to connection with torsion. As can be expected it simply consists in allowing admissible section to take values in $\F^{(1,1)}(M)$ the set of $1$-jets of sections of $\F^{(1)}(M) \to M$. Observe that $\F^{(1,1)}(M)$ is a bundle over $M$ for the canonical projection $\pi^{(1,1)} (j^1_x e) = \pi^1(e(x))$ and that its elements are also horizontal planes tangent to $\F^{(1)}(M)$, the correspondence being $j^1_xe \mapsto e_{*_x}(T_xM)$. \\ 

Consider the groupoid action $\rho^{(1)} : \b^{(1)}(\PP(M)) \times_{(\alpha, \pi^1)} \F^{(1)}(M) \to \F^{(1)}(M)$, its~derivative 
$$\rho^{(1)}_{*} : T\b^{(1)}(\PP(M)) \times_{(\alpha_*, \pi^1_*)} T\F^{(1)}(M) \to T\F^{(1)}(M) : (X_\xi, B_e) \mapsto X_\xi \cdot B_e,$$ 
and the induced groupoid action
$$\rho^{(1,1)} : \b^{(1,1)}_{nh}(\PP(M)) \times_{(\alpha, \pi^{(1,1)})} \F^{(1,1)}(M) \to \F^{(1,1)}(M) : (j^1_xb, j^1_xe) \mapsto j^1_xb \cdot j^1_xe = j^1_x (b \cdot e).$$
Observe that  
$$D \Bigl(j^1_xb \cdot j^1_xe \Bigr) = \rho^{(1)}_*(b_{*_x}(T_xM), e_{*_x}(T_xM))$$
and that $\b^{(1)}(\PP(M))$ and $\b^{(1,1)}_{nh}(\PP(M))$ act on $\F^{(1)}(M)$ and $\F^{(1,1)}(M)$ respectively through $GL(n, \R)$-equivariant maps. More is true~:

\begin{lem} The action $\rho^{(1,1)}$ of $\b^{(1,1)}(\PP(M))$ on $\F^{(1,1)}(M)$ is simply transitive, that is, given two elements $j^1_{x_1} e_1$ and $j^1_{x_2} e_2$ in $\F^{(1,1)}(M)$, there exists a unique element in $\b^{(1,1)}(\PP(M))$ mapping $j^1_{x_1} e_1$ on $j^1_{x_2} e_2$. It is denoted by $m(j^1_{x_1} e_1, j^1_{x_2} e_2)$. 
\end{lem}

\Pf Define $\xi$ to be the linear isomorphism $T_{x_1}M \to T_{x_2}M$ that maps the basis $e_1(x_1)$ to $e_2(x_2)$. Let also $\varphi : \Op\{x_1\} \to \Op\{x_2\}$ be a local diffeomorphism of $M$ such that $j^1_{x_1} \varphi = \xi$. The relation $j^1_{x_1}b \cdot j^1_{x_1}{e}_1 = j^1_{x_2}{e}_2$ is satisfied by $b$ defined by
$$b(x_1) \bigl(e_1(x_1)\bigr) = e_2(\varphi(x_1)).$$
Moreover, from its construction $j^1_xb$ belongs to $\b^{(1,1)}(\PP(M))$. Suppose $j^1_xb_o \in \b^{(1,1)}(\PP(M))$ also satisfies $j^1_{x_1}b_o \cdot j^1_{x_1}e_1 = j^1_{x_2}e_2$. Then $j^1_{x_1}(b_o^{-1} \cdot b)$ fixes $j^1_{x_1}e_1$. In other terms, 
\begin{equation}\label{DD}
D(j^1_{x_1}(b_o^{-1} \cdot b)) \cdot D(j^1_{x_1}e_1) = D(j^1_{x_1}e_1).
\end{equation} 
Let $X\in T_{x_1}M$ and let $Y$ denotes the $\alpha_*$-lift of $X$ in $D(j^1_{x_1}(b_o^{-1} \cdot b)) \subset T_{\varepsilon (x_1)} \b^{(1,1)}(\PP(M))$, then $Y = X + A$ with $A = \frac{da_t}{dt}\bigl|_{t=o} \in \Ker(\alpha_* \times \beta_*) \simeq \End(T_{x_1}M)$. Let also $Z$ denote the lift of $X$ in $D(j^1_{x_1}e_1)$. Then the relation (\ref{DD}) implies that
$$\begin{array}{cll}
Z & = & \rho_*^{(1)}\Bigl(X + A, Z \Bigr) \\
& = & \rho_*^{(1)}\Bigl(X, Z \Bigr) + \rho_*^{(1)}\Bigl(A, 0_{e_1(x_1)} \Bigr) \\
& = & \displaystyle{Z + \frac{d}{dt} a_t(e_1(x_1))\Bigl|_{t=0},}
\end{array}$$
which implies that $A$ vanishes. This holds for any $X \in T_{x_1}M$ and thus 
$$D(j^1_{x_1}(b_o^{-1} \cdot b)) = T_{x_1}M.$$
\cqfd

Now, an admissible section $s$ induces a symmetry jet ${\mathfrak s}$ through the relation~:
$${\mathfrak s} (x) = m(s(e_x), s(-e_x)),$$
where $e_x$ is some element in $\F^{(1)}(M)$, whose choice does not affect the value of $m(s(e_x), s(-e_x))$ because $m$ is $GL(n, \R)$-invariant. Moreover, its first order part $p({\mathfrak s} (x))$ is $- I_x$. As ${\mathfrak s}(x)$ consists of affine $(1,1)$-jets, the induced connection $\nabla^{\mathfrak s}$ coincides with~$\nabla^s$. \\

Conversely, given a symmetry jet ${\mathfrak s}$, we obtain an admissible section $s$ as follows. The distribution $\d^\fs$ along $-I$ associated to ${\mathfrak s}$ induces the distribution $E_{-1}^\psi$ on $\F^{(1)}(M)$ consisting of eigenspaces for the eigenvalue $-1$ of the induced involution 
$$\psi : T\F^{(1)}(M) \to T\F^{(1)}(M): B_{e_x} \mapsto \rho^{(1)}_{*} \Bigl(\ol{X}^{\d_{-I_x}, \alpha_*}, B_{e_x} \cdot (-I) \Bigr),$$ 
where $X = \pi^1_{*_{e_x}}(B_{e_x})$. The vertical tangent space $\Ker(\pi^1)_*$ consists of fixed vectors by $\psi$. Besides, 
$$\pi^1_*\circ \psi = - \psi,$$ 
hence $E_{-1}^\psi$ is $n$-dimensional and transverse to $\pi^1$. Define 
$$s : \F^{(1)}(M) \to \F^{(1,1)}(M) : e \mapsto s(e),$$
such that $s(e) = j^1_x\tilde{e}$ if and only if $D(j^1_x\tilde{e}) = (E_{-1}^\psi)_e.$
When ${\mathfrak s}$ is holonomic, $s$ takes its values in $\F^{(2)}(M) \subset \F^{(1,1)}(M)$ thanks to the following construction. First observe that if $s(j^1_0 f)$ is the unique affine $2$-jet extension (with respect to the connection $\nabla^{\mathfrak s}$ associated to ${\mathfrak s}$) of some $1$-jet $j^1_0 f \in \F^{(1)}(M)$, then ${\mathfrak s}(f(0)) \cdot s(j^1_0 f)$ is still an affine extension, namely the affine extension of $- j^1_0 f$. Now because $s$ is admissible, $s(- j^1_0 f) = s(j^1_0 f) \cdot (-I)$, where $s(j^1_0 f) \cdot (-I)$ refers to the action of the element $-I \in GL(n, \R)$ on $\F^{(2)}(M)$.  Whence $s(j^1_0 f)$ has to satisfy  the following implicit relation~:

\begin{equation}\label{implicit}
{\mathfrak s}(f(0)) \cdot s(j^1_0 f) = s(j^1_0 f) \cdot (-I),
\end{equation}
Since the $GL(n, \R)$-action commutes with the action of $\b^{(2)}(\PP(M))$, the previous relation is equivalent to 
\begin{equation}\label{implicit-bis}
s(j^1_0 f)^{-1} \cdot {\mathfrak s}(f(0)) \cdot s(j^1_0 f) = -I,
\end{equation}
where if $s(j^1_0 f) = j^2_0 f$, then $s(j^1_0 f) = j^2_{f(0)}f^{-1}$, which can be solved by means of \lref{sec-der-commutation}. Indeed, let $\theta_1$ be any element in $\F^{(2)}(M)$ that extends $j^1_0f$. Then
$$\theta_1^{-1} \cdot {\mathfrak s}(f(0)) \cdot \theta_1 = \eta_x,$$
for some $2$-jet $\eta_x$ of local diffeomorphism of $\R^n$ extending $-I$. Solutions to (\ref{implicit-bis}) are in bijective correspondence with $2$-jets $\theta$ of local diffeomorphisms of $\R^n$ that extend $I$ and satisfy
$$\theta \cdot \eta_x \cdot \theta^{-1} = -I,$$
or equivalently 
$$-I \cdot \theta \cdot \eta_x = \theta.$$
Supposing $\theta = j^2_xf$ and $\eta_x = j^2_xh$, \lref{sec-der-commutation} implies that 
$$d^2f(X_x) = d^2(-I \cdot f \cdot h)(X_x) = - d^2f(X_x) + d^2(-h)(X_x).$$
Hence $d^2f = \frac{1}{2} d^2(-h)$. 
\cqfd

\appendix

\section{Notations}\label{notation}

The following standard pieces of notation are frequently used in the text~:
\begin{enumerate}
\item For a manifold $N$, the canonical projection $TN \to N$ is denoted most often by the letter $p$ and sometimes by the symbol $p^N$. In particular, $p : T^2M \to TM$ is the projection of $TN$ onto $N$ for $N = TM$.
\item The canonical inclusion of a manifold into its tangent bundle as the zero section is systematically denoted by the letter $i$.
\item When $\pi : M \to N$ is a submersion, $i^\pi : T^\pi M \hookrightarrow TM$ denotes the canonical inclusion of the sub-bundle $T^\pi M$ consisting of vectors tangent to the fibers of $\pi$ in $TM$. 
\item If $\pi : M \to N$ is a submersion and $\varphi : N_0 \hookrightarrow N$ is an immersion, the notation $T_{N_0}^\pi M$ stands for the set of vectors tangent to the fibers of $\pi$ and belonging to $\pi^{-1}(\varphi(N_0))$. 
\item Given a fibration $\pi : E \to B$, the fiber $\pi^{-1}(b)$ is denoted by $E^\pi_b$ or $E_b$ when unambiguous and the natural inclusion of $E^\pi_b$ in the total space $E$ is denoted by $i^\pi_b : E^\pi_b \to E$.
\item For a vector bundle $\pi : E \to B$, we have a natural inclusion 
$$i^\pi_e : E_{\pi(e)} \hookrightarrow T_e^\pi E : e' \mapsto \frac{d(e + te')}{dt}\Bigr|_{0}$$ for each element $e$ in $E$. Notice that $i^\pi_e = (i^\pi_{\pi(e)})_{*_e}$.
\item If $\pi_1 : M_1 \to N$ and $\pi_2 : M_2 \to N$ are two submersions. Then $M_1 \times_{(\pi_1, \pi_2)} M_2$ denotes the fiber-product manifold $\{(m_1, m_2) \in M_1 \times M_2; \pi_1(m_1) = \pi_2(m_2)\}$. It is naturally endowed with a projection onto N. Observe that 
$$T(M_1 \times_{(\pi_1, \pi_2)} M_2) = TM_1 \times_{(\pi_{1*}, \pi_{2*})} TM_2.$$
\item Given two vector bundles $\pi_1 : E_1 \to B_1$ and $\pi_2 : E_2 \to B_2$ and a map $\varphi : E_1 \to E_2$, the expression $\varphi : (E_1, \pi_1) \to (E_2, \pi_2)$ means that $\varphi$  is a morphism of vector bundles.
\end{enumerate}

\section{Lie groupoids of jets of bisections}\label{jet-of-bis}

We refer to the books \cite{dSW}, \cite{Mck-87} and \cite{Mck-05} for an introduction to the basics on groupoids. Notice though that,  like Mackenzie but unlike da Silva and Weinstein, we call the source map $\alpha$ and the target map $\beta$. The following notions will be needed in the text~: Lie groupoid, locally trivial Lie groupoid, groupoid morphism, Lie subgroupoid, bisection, local bisection, Lie groupoid action, Lie algebroid. \\
 
For a Lie groupoid $G \rightrightarrows M$, the symbol $\alpha$ denotes the source map, $\beta$ the target map, $m : G \times_{(\alpha, \beta)}G \to G$ the multiplication, $\varepsilon : M \to G$ the natural inclusion and $\iota$ the inversion. The collection of smooth local bisections of $G$ is denoted by $\b_{\ell}(G)$ and the subcollection of local bisections defined on a neighborhood of $x \in M$ by $\b_{\ell, x}(G)$. Let us now introduce the jet extensions groupoids (cf.~\cite{Ehresmann}). Given a Lie groupoid $G \rightrightarrows M$, a point $x \in M$ and any integer $k \geq 0$, there is an equivalence relation $\sim^k_x$ on the set $\b_{\ell, x}(G)$ of local bisections of $G \rightrightarrows M$ defined near $x$, namely $b_1\sim^k_x b_2$ if $b_1(x) = g = b_2(x)$ and $b_1$ and $b_2$ have the same Taylor expansion of order $k$ at $x$ with respect to some (hence any) pair of  local coordinate systems around $x$ and $g$. The equivalence class of $b$ with respect to $\sim^k_x$ is commonly denoted by $j^k_xb$ and called the $k$-jet of $b$ at $x$. 

\begin{dfn} The set of all $k$-jets of local bisections of $G$ is denoted by $\b^{(k)}(G)$ and called the $k$-jet extension of the Lie groupoid $G$.
\end{dfn}

The set $\b^{(k)}(G)$ is another Lie groupoid over $M$ when it is endowed with the obvious structure maps, namely
\begin{enumerate}
\item[-] $\alpha : \b^{(k)}(G) \to M : j^k_xb \mapsto x$, 
\item[-] $\beta : \b^{(k)}(G) \to M : j^k_xb \mapsto \beta \circ b (x)$, 
\item[-] $m(j^k_yb' , j^k_xb) = j^k_x (b' \cdot b)$ when $y = \beta \circ b(x)$,
\item[-] $\varepsilon : M \to \b^{(k)}(G) : x \mapsto j^k_x\varepsilon$,
\item[-] $\iota :  \b^{(k)}(G) \to  \b^{(k)}(G) : j^k_xb \mapsto j^k_yb^{-1}$, where $y = \beta \circ b(x)$. 
\end{enumerate}
Let us observe that for $k > l$, the natural projection $p^{k \to l} : \b^{(k)}(G) \to \b^{(l)}(G) : j^k_xb \mapsto j^l_xb$ is a Lie groupoid morphism. 

\begin{rmk}\label{1jets-as-planes}
Let $G \rightrightarrows M$ be a Lie groupoid and set $n = \dim M$. Consider the Grassmann bundle $\pi :\Gr_n(G) \to G$ of $n$-planes tangent to $G$. Notice that the data of an element $j^1_xb$ in $\b^{(1)}(G)$ is equivalent to the data of the map $b_{*_x} : T_xM \to T_{b(x)}G$, itself equivalent to the data of the plane $\im(b_{*_x})$. Moreover, the map 
\begin{equation}\label{D}
D : \b^{(1)}(G) \to \Gr_n(G) : j^1_xb \mapsto D(j^1_xb) = b_{*_x}(T_xM)
\end{equation} 
is a diffeomorphism onto the open subset of $\Gr_n(G)$ consisting of ``horizontal" planes, that is, planes that are transverse to the $\alpha$-fibers and the $\beta$-fibers. The latter subset is denoted by $\Gr^h_n(G)$ and supports thus a groupoid structure whose source and target map are $\alpha \circ \pi$ and $\beta \circ \pi$ respectively and whose multiplication is induced from the differential of the multiplication 
$$m_{*_{(g_1, g_2)}} : T_{g_1}G \times_{(\alpha_{*_{g_1}}, \beta_{*_{g_2}})} T_{g_2}G \to T_{g_1 \cdot g_2}G.$$ 
in $G$. The identity at $x$ is $T_xM$ and the inverse of $D \subset T_gG$ is $\iota_{*_g}(D)$.
\end{rmk}

\begin{nota}\label{b1} Consider the pair groupoid $\PP(M)$, that is the set $M \times M$, endowed with the two projections $\alpha = p_2$ and $\beta = p_1$ onto $M$ and the product $(x, y) \cdot (y, z) = (x,z)$. A local bisection $b$ of $\PP(M)$ is a partial diffeomorphism $U_b = \alpha(b) \to V_b = \beta(b)$. Therefore the groupoid $\b^{(k)}(\PP(M))$ is the proper subset of $J^k(M,M)$ consisting of $k$-jets of partial diffeomorphisms of $M$. In particular, the groupoid $\b^{(1)}(\PP(M))$ is the set of linear maps between any pair of tangent spaces to $M$. It is called in the literature the \emph{general linear groupoid} of the vector bundle $TM$ or the \emph{gauge groupoid} of the principal bundle of frames of $M$ and also denoted by ${\rm GL}(TM)$ or ${\rm Aut}({\mathcal F}(M))$. 
\end{nota}

The extension procedure can be iterated and the groupoid $\b^{(k_1)}(\b^{(k_2)}(G))$ is denoted hereafter by $\b^{(k_1,k_2)}(G)$. It contains $\b^{(k_1 + k_2)}(G)$ as an embedded subgroupoid. Observe that, in addition to the natural projections 
$$p^{k_1 \to l_1} : \b^{(k_1,k_2)}(G) \to \b^{(l_1, k_2)}(G) : j^{k_1}_xb \mapsto j^{l_1}_xb,$$
for $l_1 = 0, ..., k_1-1$, there are projections 
$$p^{k_2\to l_2}_* : \b^{(k_1,k_2)}(G) \to \b^{(k_1, l_2)}(G) : j^{k_1}_xb \mapsto j^{k_1}_x(p^{k_2 \to l_2} \circ b)$$
that are groupoid morphisms as well. On $\b^{(k_1+k_2)}(G) \subset \b^{(k_1,k_2)}(G)$, the map $p^{k_2\to 0}_*$ coincides with 
$p^{k_1+k_2 \to k_1}$. Similarly, given a sequence $(k_1, ..., k_I)$ of natural numbers, the groupoid $\b^{(k_1)}(...(\b^{(k_I)}(G))...)$ is denoted by $\b^{(k_1, ..., 
k_I)}(G)$ and supports a series of projections 
$$p^{k_i \to l_i}_{\underbrace{* ... *}_{i-1}} : \b^{(k_1, ..., k_i, ..., k_I)}(G) \to \b^{(k_1, ..., l_i, ..., k_I)}(G).$$ 
Notice, for instance, that any groupoid $\b^{(k_1, ..., k_I)}(G)$ is a subgroupoid of the groupoid $\b^{(1,...,1)}(G)$ with $k_1 + ... + k_I$ copies of~$1$. 

\begin{nota}\label{bkl}
Provided it does not generate any ambiguity, we will use the following abbreviations~:
\begin{enumerate}
\item[-] The projection $p^{k \to 0}$ on $\b^{(k)}(G)$, that extracts the $0$-th order part of a jet, is denoted by $p^k$ and coincides with $p \circ ... \circ p$ ($k$ factors) on $\b^{(1, ..., 1)}(G) \supset \b^{(k)}(G)$.
\item[-] Similarly, the projection $p^{k_i \to 0}_{* ... *}$ is denoted by $p^{k_i}_{* ... *}$ and coincides with $p_{* ... *} \circ ... \circ p_{* ... *}$ on $\b^{(k_1, ..., k_{i-1}, 1, ..., 1, k_{i+1}, ..., k_I)} \supset \b^{(k_1, ..., k_i, ..., k_I)}$. 
\item[-] We remove the superscripts from the projections $p^1$, $p_*^1$, ..., $p^1_{* ...*}$ and denote them by $p$, $p_*$, $p_{*... *}$. Observe that
\begin{equation}\label{p*=p*}
D\Bigl( p_{\underbrace{* ... *}_{i}}(j^1_xb)\Bigr) = (p_{\underbrace{* ... *}_{i-1}})_{*_{b(x)}} \Bigl( D(j^1_xb)\Bigr).
\end{equation}
\end{enumerate} 
\end{nota}

A local bisection $b$ of $G$, induces a local so-called {\bf holonomic} bisection  
$$j^kb : U \to \b^{(k)}(G) : x \mapsto j^k_xb$$
of the groupoid $\b^{(k)}(G)$. When $G \rightrightarrows M$ is locally trivial, there is a nice characterization of local holonomic bisections of $\b^{(1)}(G)$ as local bisections tangent to a certain distribution that we introduce hereafter. Consider the map $p : \b^{(1)}(G) \to G$ and its differential $p_{*_\xi} : T_\xi \b^{(1)}(G) \to T_{p(\xi)}G$ at $\xi \in\b^{(1)}(G)$. Observe that $D(\xi)$ is a subspace in $T_{p(\xi)}G$.
\begin{dfn}\label{def-e}
The holonomic distribution $\e$ or $\e^G$ on $\b^{(1)}(G)$ is defined by
$$\e_\xi = p_{*_\xi}^{-1} \bigl(D(\xi)\bigr).$$
\end{dfn}

\begin{prop}\label{e}
The distribution $\e$ has rank $(n + n(k - n))$, where $n = \dim M$ and $k = \dim G$, contains the distribution $\Ker p_*$ and is transverse to the $\alpha$ and $\beta$-fibers of $\b^{(1)}(G)$. It has the property that a local bisection $b : U \to \b^{(1)}(G)$ is holonomic if and only if it is tangent to $\e$. 
\end{prop}

\Pf Observe that $\Ker p_{*_\xi}$ is $(n(k - n))$-dimensional. Therefore, the plane $\e_\xi$ has rank $\dim D(\xi) + \dim \Ker p_{*_\xi} = n + n(k-n)$. \\

Because $D(\xi)$ is transverse to the $\alpha$ and $\beta$-fibers of $G$, its lift $p_*^{-1} (D(\xi))$ enjoys the same property in $\b^{(1)}(G)$. \\

Let $b$ be a local bisection in $\b^{(1)}(G)$. The condition $T_{b(x)}b \subset \e$ is equivalent to $p_*(T_{b(x)}b) = D(b(x))$, that is $p_*(j^1_xb) = b(x)$ (cf.~(\ref{p*=p*})) or $j^1_x(p \circ b) = b(x)$. So a bisection $b$ is tangent to $\e$ everywhere if and only if it coincides with $j^1b^0$.
\cqfd

\begin{rmk} In general, a bisection $b$ of $\b^{(1,...,1)}(M)$ is holonomic if $b$ is tangent to $\e$ and its projection $p \circ b$ is a holonomic bisection. 
\end{rmk}

\begin{rmk}\label{sh} The common terminology in the literature (cf.~e.g.~\cite{Lib}) is as follows. Elements in $\b^{(k)}(G)$ are called {\em holonomic} jets, elements in $\b^{(1,..., 1)}(G)$ ($k$ copies of $1$) that do not belong to $\b^{(k)}(G)$ are called {\em non-holonomic} jets. Amongst non-holonomic jets are {\em semi-holonomic} jets who are those non-holonomic jets for which the values of the projections $p, p_*, ..., p_{*...*}$ all agree.
\end{rmk}

\begin{dfn}\label{derived-action} A left action $\rho : G \times_{(\alpha, \pi)} E \to E$ of a groupoid $G \rightrightarrows M$ on a fiber bundle $\pi : E \to M$ (\cite{Mck-05}) induces an action of the groupoid $\b^{(1)}(G)\rightrightarrows M$ onto $\pi \circ p : TE \to M$ as follows~:
$$\rho^{(1)} : \b^{(1)}(G) \times_{(\alpha, \pi \circ p)}TE \to TE : \Bigl(j^1_x b, X_e\Bigr) \mapsto j^1_xb \cdot X_e = \rho_{*} \Bigl(b_{*_x}\bigl(\pi_{*_e} (X_e)\bigr), X_e\Bigr),$$
where $\rho_* : TG \times_{(\alpha_*, \pi_*)} TE \to TE$ is the differential of $\rho$. 
Iterating this procedure, we obtain actions
$$\rho^{(1,...,1)} : \b^{(1,...,1)} \times_{(\alpha, \pi \circ p^k)} T^kE \to T^k E$$
for any $k = 1, 2, ...$. 
\end{dfn}

\begin{rmk}\label{derived-action-of-pair-groupoid} In particular, starting from the trivial action of the pair groupoid $\PP(M)$ on $M$
$$\rho : (M\times M) \times_{(\alpha, \id)} M  \to M : ((y, x), x) \mapsto y$$
we obtain actions $\rho^{(1)}$, $\rho^{(1,1)}$, $\rho^{(1,1,1)}$, ... of $\b^{(1)}, \b^{(1,1)}, \b^{(1,1,1)}, ...$ on $TM$, $T^2M$, $T^3M, ...$ respectively. A groupoid action $\rho : G \times_{(\alpha, \pi)} E \to E$ is said to be {\sl effective} if $\rho(g_1, e) = \rho(g_2, e)$ for all $e \in E$ with $\pi(e) = \alpha(g_1) = \alpha(g_2)$ implies that $g_1 = g_2$. The various actions $\rho^{(1)}$, $\rho^{(1,1)}$, $\rho^{(1,1,1)}$, ... are effective. 
\end{rmk}

\begin{lem}\label{loc-triv} The $k$-jet extension of a locally trivial groupoid is locally trivial as well. In particular, all the groupoids $\b^{(... ,l,k)}(M) \rightrightarrows M$ are locally trivial.
\end{lem}

\section{Second tangent bundle}\label{second-t-b}

Given a manifold $M$, its second tangent bundle $T^2M$ is defined to be the tangent bundle $p : T(TM) \to TM$ of the total space of the tangent bundle to $M$. Its elements are denoted by calligraphed letters $\X, \Y, \ZE, ....$. It is endowed with several pieces of structure that we describe in the present section. Note that $T^2M$ is an example of a double vector bundle as treated in Chapter 9 of \cite{Mck-05}. The main specificity of $T^2M$ is that its \emph{flip} is canonically isomorphic to itself.\\

\noindent 
{\bf Vector bundle structures~:} $T^2M$ admits two distinct structures of vector bundle over the manifold $TM$~: \\

$\bullet$ $p : T^2M \to TM$ is the canonical projection, that is, if $t \mapsto X_t$ is a path in $TM$, then 
$$p\Bigl(\frac{dX_t}{dt}\Bigl|_{t=0}\Bigr) = X_0.$$ 
Fiberwise addition and scalar multiplication are denoted, as usual, by $+ : TM \times_{(p,p)}TM : (\X_1, \X_2) \mapsto \X_1 + \X_2$ and $m_a : \R \times TM : (a, \X) \mapsto a\X = m_a(\X)$ respectively. The fiber over a vector $X_x \in TM$ is denoted by $T_{X_x}TM$. \\

$\bullet$ $p_* : T^2M \to TM$ is the differential of the canonical projection $p : TM \to M$, that is 
$$p_*\Bigl(\frac{dX_t}{dt}\Bigl|_{t=0}\Bigr) = \frac{dp(X_t)}{dt}\Bigl|_{t=0}.$$
Fiberwise addition is the differential of the corresponding map on $TM$, that is 
$$+_* : T^2M \times_{(p_*, p_*)} T^2M : \Bigl(\X = \frac{dX_t}{dt}\Bigl|_0, \Y = \frac{dY_t}{dt}\Bigl|_0\Bigr) \mapsto \frac{dX_t + Y_t}{dt}\Bigl|_0,$$
where the path $t \mapsto X_t$ and $t \mapsto Y_t$ have been chosen to satisfy $p(X_t) = p(Y_t)$. This is not restrictive since $p_*(\X) = p_*(\Y)$. Similarly, scalar multiplication by a real $a$ is the differential of $m_a$ on $TM$~: 
$$m_{a*} : T^2M \to T^2M : \X = \frac{dX_t}{dt}\Bigl|_{t=0} \mapsto m_{a*}(\X) = \frac{d (a X_t)}{dt}\Bigl|_{t = 0}.$$
The $p_*$-fiber over a vector $X_x \in TM$ is denoted by $T^{X_x}TM$.

\begin{figure}[h!]  
\begin{center}
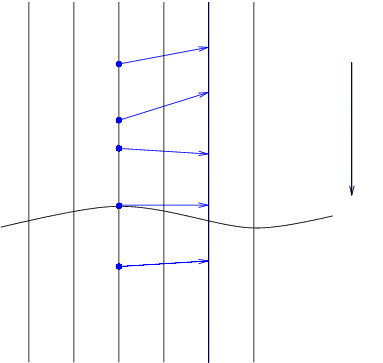 
\caption{The addition and scalar multiplication in $(T^2M, p_*)$}
\end{center}
\end{figure}

\begin{lem}
The map $p : T^2M \to TM$ (\rp $p_* : T^2M \to TM$) is a vector bundle morphism when $T^2M$ is endowed with the vector bundle structure induced by $p_*$ (\rp $p$). 
\end{lem}
The to maps $p \circ p$ and $p \circ p_*$ of $T^2M$ onto $M$ agree and are denoted by $p^2$. In other words the following diagram is commutative. The fiber of $p^2$ over a point $x \in M$ is denoted by $T^2_xM$. 

\vspace{.7cm}
\hspace{3cm}\xymatrix{
& T^2M\ar@{->}[ld]_p \ar@{->}[rd]^{p_*} & \\
TM \ar@{->}[rd]_p & & TM\ar@{->}[ld]^p \\
& M &}
\vspace{.7cm}\\

\noindent
{\bf Horizontal inclusions and projections~:} The two vector bundle structures $p, p_* : T^2M \to TM$ induce two natural inclusions of $TM$ into $T^2M$~:
$$i, i_* : TM \to T^2M : X_x \mapsto i(X_x) = 0_{X_x}, i_*(X_x) = {0_{*X_x}},$$ 
parameterizing the respective zero-sections denoted  by $0_{TM}, 0_{*TM}$. The map $i_*$ is the differential of the inclusion $i : M \to TM$ of $M$ as the zero-section of $TM$. The inclusion $i$ (\rp $i_*$) is a vector bundle morphism between $(TM, p)$ and $(T^2M, p_*)$ (\rp $(T^2M, p)$). The associated projections of $T^2M$ onto its two zero-sections are denoted by $e = i \circ p$ and $e_* = i_* \circ p_* $. Let ${\bf i}$ denote the inclusion $i \circ i = i_* \circ i$ of $M$ into $T^2M$. It parameterizes the intersection of the two zero-sections $0_{TM} \cap {0_*}_{TM} = 0_{0_M}$. Similarly, let ${\bf e}$ denote the projection $e \circ e_* = e_* \circ e$ of $T^2M$ onto $0_{0_M}$. These different maps satisfy the relations~:
$$\begin{array}{llll}
p \circ i = id_{TM} \;\;& p \circ e = p \;\;& p_* \circ i = i \circ p \;\;&  p_* \circ e = e \circ p \\
p_* \circ i_* = id_{TM} \;\; & p_* \circ e_* = p_* \;\; & p \circ i_* = i \circ p \;\; & p \circ e_* = e \circ p \\
\end{array}$$
and $${\bf e} = i^2 \circ p^2.$$

\noindent 
{\bf Vertical inclusion of $TM$~:} There is a canonical ``vertical inclusion" from $TM$ into $T^2M$ parameterizing $T^p_{0_M}TM = p^{-1}(0_M) \cap \, p_*^{-1}(0_M)$ (cf.~\sref{notation})
$$i^p_{0_M} : TM \to T^2M : X_x \mapsto i^p_{0_M}(X_x) = \frac{dtX_x}{dt} \Bigl|_{t=0}.$$
The restriction of $i^p_{0_M}$ to $T_xM$ is sometimes denoted by $i^p_{0_x}$. Observe that the map $i^p_{0_M}$ is a vector bundle map for both structures on $T^2M$. In particular, a vector $\X \in T^2M$ belongs to $\im (i^p_{0_M})$ if and only if $m_a(\X) = m_{a*}(\X)$ for all $a \in \R$.\\

\noindent 
{\bf Vertical inclusions of $TM \oplus TM$~:} There are also two canonical inclusions of $TM \oplus TM$ into $T^2M$ parameterizing the two ``kernels" $p^{-1}(0_M) = T_{0_M}TM$ and $p_*^{-1}(0_M) = T^pTM$ as follows~:
$$I_{p} : TM \oplus TM \to p^{-1}(0_M) \subset T^2M : (X^1_x, X_x^2) \mapsto i_{*_x}(X_x^1) + i^p_{0_M}(X_x^2)$$
$$I_{p_*} : TM \oplus TM \to p_*^{-1}(0_M) \subset T^2M : (X^1_x, X_x^2) \mapsto i(X_x^1) +_* i^p_{0_M}(X_x^2).$$
The map $I_{p}$ (\rp $I_{p_*}$) is a vector bundle morphism between $TM \oplus TM$ and $(T^2M, p)$ (\rp $(T^2M, p_*)$).

\begin{rmk}\label{Leibniz} Given vector fields $X$, $Y$ on $M$ and a function $f \in C^\infty(M)$, the section $fY$ of $p : TM \to M$ satisfies
\begin{equation}\label{scalar-mult}
(fY)_{*_x} (X_x) = X_xf \; Y_x + m_{f(x)*}(Y_{*_x}X_x),
\end{equation}
where $X_x f \; Y_x$ really means $I_{p_*}(f(x)Y_x, X_xf Y_x) = i(f(x)Y_x) +_* i^p_{0_M}(X_xf \; Y_x)$. Indeed, setting $(f \circ p) \times \id : TM \to \R \times TM : Z_x \mapsto (f(x), Z_x)$ and $m : \R \times TM \to TM : (a, X_x) \mapsto a X_x$, the section $fY$ can be rewritten as 
$$m \circ ((f \circ p) \times \id) \circ Y.$$ 
Hence 
$$\begin{array}{ccl}
(fY)_{*_x} (X_x) & = & m_{*_{(f(x), Y_x)}} \Bigl(f_{*_x} (X_x), Y_{*_x}X_x\Bigr) \\
& = & m_{*_{(f(x), Y_x)}} \Bigl(f_{*_x} (X_x), 0_{Y_x} \Bigr) + m_{*_{(f(x), Y_x)}} \Bigl(0_{f(x)}, Y_{*_x}X_x\Bigr) \\
& = & \Bigl. X_xf \; Y_x + m_{f(x)*} (Y_{*_x}X_x),
\end{array}$$
\end{rmk}

\noindent
{\bf Affine structure~:} The product $P = p \times p_* : T^2M \to TM \times_{(p, p)}TM$ yields on $T^2M$ yet another structure, of affine bundle of rank $n$ this time, whose fiber over $(X_x, Y_x)$ is modeled on the vector space $T_xM$ and denoted by $T_{X_x}^{Y_x}TM$. Observe that for a fixed vector $\X \in T^2M$, with $p(\X) = X_x$ and $p_*(\X) = Y_x$ the two maps
\begin{equation}\label{A}
\begin{array}{l}
A_\X : T_xM \to T_{X_x}^{Y_x}TM : V_x \mapsto \X + \Bigl(e (\X) +_* i^p_{0_M}(V_x)\Bigr) \\
A_\X : T_xM \to T_{X_x}^{Y_x}TM : V_x \mapsto \X +_* \Bigl(e_*(\X) + i^p_{0_M}(V_x)\Bigr)
\end{array}
\end{equation}
coincide and parameterize the fiber of $P$ through $\X$. Moreover, there is a canonical map
\begin{equation}\label{i}
\pi : T^2M \times_{(P, P)}T^2M \to TM : (\X_1, \X_2) \mapsto \pi(\X_1, \X_2),
\end{equation}
defined by
\begin{equation}\label{def-i}
\pi(\X_1,  \X_2) = V_x  \quad \mbox{if} \quad \X_1 =  A_{\X_2}(V_x) 
\end{equation}
One could also write, with a slight abuse of notation, the expression $\pi(\X_1,  \X_2) = \X_1 - \X_2$. The map $\pi$ satisfies 
\begin{equation}\label{abc}
\pi(\X^1, \X^2) = \pi(\X^1 +_o \X, \X^2 +_o \X)
\end{equation}
when $+_o$ denotes either $+$ or $+_*$ and $p_o(\X) = p_o(\X^i) \in T^2M$ for the corresponding projection $p_o$. \\

\noindent
{\bf Canonical Involution~:} A very important piece of the structure of $T^2M$ is its canonical involution, defined below~:

\begin{dfn}\label{involution} The canonical involution on $T^2M$ is commonly defined by means of local coordinates $(x^i, X^i, Y^i, \X^i)$ induced by local coordinates $x^i$ on $M$ as being the map $\kappa = \kappa_M: T^2M \to T^2M$ that flips the two middle sets of coordinates 
$$\kappa(x^i, X^i, Y^i, \X^i) = (x^i, Y^i, X^i, \X^i).$$
\end{dfn}

\noindent
{\bf Properties of $\kappa$~:} The involution $\kappa$ is an isomorphism between the two distinct vector bundle structures on $T^2M$, and, as such, satisfies the relations~:
\begin{equation}\label{prop-kappa}
\begin{array}{ll}
p_* \circ \kappa = p & p \circ \kappa = p_*\\
\kappa \circ i = i_* & \kappa \circ i_* = i \\
\kappa \circ e = e_* & \kappa \circ e_* = e \\
\kappa \circ m_{a*} = m_a \circ \kappa & \kappa \circ m_a = m_{a*} \circ \kappa .\\
\kappa(\X + \Y) = \kappa(\X) +_* \kappa(\Y) & \kappa(\X +_* \Y) = \kappa(\X) + \kappa(\Y) 
\end{array}
\end{equation}
It is thus also an endomorphism of the affine bundle $P : T^2M \to TM \times_{(p,p)} TM$ over the reflection map $\kappa_o(X,Y) = (Y,X)$. Furthermore, $\kappa$ fixes pointwise the image of $i_{0_M}^p$~:
$$\kappa \circ i^p_{0_M} = i^p_{0_M}.$$ 
This and the relations (\ref{A}) imply that for two vectors $\X_1$ and $\X_2$ in a same $P$-fiber,
\begin{equation}\label{i-kappa}
\pi(\X_1, \X_2) = \pi(\kappa(\X_1), \kappa(\X_2)).
\end{equation}

The following properties will not be used in the text but provide additional insight on $\kappa$ and its relation with the Lie bracket and the flows of vector fields.
\begin{prop}\label{involution-prop}(\cite{L}, \cite{K-}) The involution $\kappa$, sometimes denoted by $\kappa^M$, satisfies the following properties~: let $x \in M$ and $X, Y \in {\mathfrak X}(M)$
\begin{enumerate}
\item[-] $[X, Y]_x = \pi \bigl(Y_{*_x}X_x, \kappa(X_{*_x} Y_x) \bigr)$,
\item[-] $\displaystyle{\kappa(Y_{*_x} X_x) = \frac{d(\varphi_Y^t)_{*_x}X_x}{dt}\Bigl|_{0}}$, where $\varphi_Y^t$ denotes the flow of $Y$ at time $t$.
\end{enumerate}
\end{prop}

\begin{rmk}
It would be nice to have an intrinsic description of $\kappa$. Each one of the previous properties could serve as a  definition of $\kappa$ on $T^2M - T^pTM$. On $T^pTM$, we want $\kappa$ to coincide with
$$\kappa \Bigl(i(X_x) + i_{0_M}(V_x) \Bigr) = i_*(X_x) +_* i_{0_M}(V_x).$$
So $\kappa$ can be defined intrinsically separately on $T^2M - T^pTM$ and $T^pTM$ but smoothness of $\kappa$ across $T^pTM$ has then to be established. So this approach does not seem completely satisfactory.
\end{rmk}

\begin{rmk}
The involution $\kappa$ allows for an alternative characterization of distributions on $M$ that are involutive~: a distribution $\d$ on $M$ is involutive if and only if the subset $T\d \cap p_*^{-1}(\d)$ of $T^2M$ is $\kappa$-invariant ($\d$ is thought of a subbundle of $TM$; therefore $T\d \subset T^2M$).
\end{rmk}

\section{Structure of $\b^{(1,1)}(\PP(M))$}\label{b11}

The structure of the groupoid $\b^{(1,1)}(\PP(M))$ (cf.~\aref{jet-of-bis} and \nref{b1}) follows closely that of $T^2M$. As already observed in \aref{jet-of-bis}, it is endowed with two natural projections $p$ and $p_*$ onto $\b^{(1)}(\PP(M))$, whose definition we briefly recall. An element of $\b^{(1,1)}(\PP(M))$ is of the type $j^1_xb$ for some local bisection $b : U_x \to \b^{(1)}(\PP(M))$ defined in a neighborhood $U_x$ of $x$. Then
$$\left\{ \begin{array}{lllllllll}
p &:& \b^{(1,1)}(\PP(M)) &\to& \b^{(1)}(\PP(M)) &:& j^1_xb &\mapsto& b(x) \\
p_* &&&&&&&& j^1_x(p \circ b).
\end{array}\right.$$
We thus obtain a commutative diagram~: 

\vspace{.2cm}
\hspace{1.3cm}
\xymatrix{
& \b^{(1,1)}(\PP(M)) \ar@{->}[ld]_p \ar@{->}[rd]^{p_*} & & \\
\b^{(1)}(\PP(M)) \ar@{->}[rd]_p & & \b^{(1)}(\PP(M))\ar@{->}[ld]^p & (\star)\\
& M & &}
\vspace{.5cm}
Notice that $b^0 = p \circ b$ is a local bisection of the pair groupoid $M\times M$, that is a section $x \mapsto (f(x), x)$ of $\alpha$ such that $\beta \circ b^0 : x \mapsto f(x)$ is a local diffeomorphism of $M$. In the sequel $b^0 = p \circ b$ will systematically be identified with the local diffeomorphism $f = \beta \circ b^0$.

\begin{rmk}\label{bouncing} Under the identification of $\b^{(1,1)}(\PP(M))\sim \Gr_n^h(\b^{(1)}(\PP(M)))$ (cf.~\rref{1jets-as-planes}), the map $p_* : \b^{(1,1)}(\PP(M)) \to \b^{(1)}(\PP(M))$ coincides with the ``bouncing map", defined as follows
$$\fb : \Gr_n^h(\b^{(1)}(\PP(M))) \to \b^{(1)}(\PP(M)) : P_\xi \mapsto \fb(P_\xi) = \beta_{*_\xi} \circ \Bigl(\alpha_{*_\xi}\bigr|_{P_\xi}\Bigr)^{-1}.$$ 
\begin{figure}[h!]  
\begin{center}
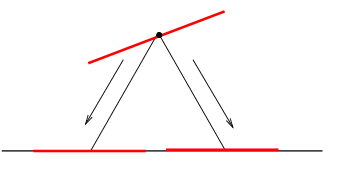 
\end{center}
\end{figure}
\end{rmk}

When $b$ is a local bisection of $\b^{(1)}(\PP(M))$, the map $\fb$ applied to $Tb$ yields a holonomic bisection~: 
$$x \mapsto \fb(T_{b(x)}b) = j^1_xb^0.$$
\begin{nota}\label{nota-e}
The holonomic distribution $\e^{M\times M}$ on $\b^{(1)}(\PP(M))$  (cf.~\dref{def-e}) is denoted by $\e^{(1)}$.
\end{nota}
Observe that if a $n$-plane $P_\xi$ is contained in the holonomic distribution $\e^{(1)}$ then $\fb(P_\xi) = \xi$. (One could describe $\e^{(1)}$ as the union of all horizontal $n$-planes having this property.) More generally, the linear maps $\fb(P_\xi)$ and $\xi$ coincide on $\alpha_*(P_\xi \cap \e^{(1)}_\xi)$. 

\begin{rmk} The distribution $\e^{(1)}$ is maximally non-integrable. 
It generalizes the canonical contact form $\alpha$ on the set of $1$-jets of local maps from $M$ to the real line $J^1(M \times \R) \simeq T^*M \times \R$, defined by $\alpha(X_\beta, V) = \beta(\pi_*(X_\alpha)) - V$ to the case of $1$-jets of maps from $M$ to $M$.
\end{rmk}

\begin{dfn}\label{bh}
The set of $(1,1)$-jets $\xi = j^1_xb$ for which $p(\xi) = p_*(\xi)$ is denoted hereafter $\b^{(1,1)}_{sh}(\PP(M))$ (``sh" for semi-holonomic, cf.~\rref{sh}). Notice that $\b^{(2)}(\PP(M)) \subsetneqq \b^{(1,1)}_{sh}(\PP(M))$. 
\end{dfn}

\begin{rmk}\label{h-e}
A $(1,1)$-jet $\xi$ belongs to $\b^{(1,1)}_{sh}(\PP(M))$ if and only if $D(\xi) \subset \e^{(1)}$.
\end{rmk}

\begin{prop}\label{e-iota}
Along the bisection $b_o = -I$, the distribution $\e^{(1)}$ coincides with the family of $(-1)$-eigenspaces of the involution $\iota_{*}$.
\end{prop}

\Pf Let $x\in M$ and set $\xi = -I_x \in \b^{(1)}(\PP(M))$. Then 
$$\iota_{*_\xi} : T_\xi \b^{(1)}(\PP(M)) \to T_\xi \b^{(1)}(\PP(M))$$ 
is an involutive linear isomorphism. Hence the tangent space to $\b^{(1)}(\PP(M))$ at $\xi$ decomposes into a direct sum of $+1$ and $-1$ eigenspaces for $\iota_{*_\xi}$~:
$$T_\xi \b^{(1)}(\PP(M)) =  E_{+1} \oplus E_{-1}.$$ 
We claim that $\e^{(1)}_\xi = E_{-1}$. Let $X_\xi \in T_\xi\b^{(1)}(\PP(M))$. Then, because $p$ is a groupoid morphism and $\iota(\xi) = \xi$, we have
$$p_{*_{\xi}} \circ \iota_{*_\xi}(X_\xi) = \iota_{*_{p(\xi)}} \circ p_{*_\xi}(X_\xi).$$
Moreover $\iota_{*_{p(\xi)}} : T_xM \times T_xM : (X^1, X^2) \mapsto (X^2, X^1)$. Thus $X_\xi \in E_{-1}$ implies that $p_{*}(X_\xi) \in D(\xi)$. Conversely, $p_{*}(X_\xi) \in D(\xi)$ implies that $X_\xi \in E_{-1} + \Ker p_{*_\xi}$. But $\Ker p_{*_\xi} \subset E_{-1}$. Indeed, if $\exp(tA)$ is a one-parameter subgroup in the Lie group $p^{-1}(x,x) \subset \b^{(1,1)}(\PP(M))$, then $\iota (\exp (tA)) = \exp(-tA)$, whence $\iota (\xi \cdot \exp(tA)) = \xi \cdot \exp (-tA)$ ($\xi$ is central), which implies that $\iota_{*_\xi} (X) = -X$ for $X \in \Ker p_{*_\xi}$. 
\cqfd

\begin{cor}\label{invertibility} Any element of $\b^{(1,1)}(\PP(M))$ whose first order part belongs to the bisection $-I \subset \b^{(1)}(\PP(M))$ is its own inverse.
\end{cor}

\Pf Let $\xi \in \b^{(1,1)}_{sh}(\PP(M))$, with $p(\xi) = -I_x$ then $D(\xi) \subset \e^{(1)}_{-I_x}$. Hence, from \rref{1jets-as-planes}, we know that $D(\iota(\xi))$ coincides with $\iota_{*_{-I_x}}(D(\xi))$ which equals $D(\xi)$ by the previous result, implying that $\iota(\xi) = \xi$. 
\cqfd

\begin{rmk}\label{f**}
The $2$-jet $j^2_xf$ of a local diffeomorphism $f : U \subset M \to M$  of $M$ is equivalently described as the map 
$$f_{**_x} : T^2_xM \to T^2_yM : \X_{X_x} \mapsto (f_{*})_{*_{X_x}}(\X_{X_x}) \stackrel{\rm not}{=} f_{**_{X_x}}(\X_{X_x}).$$
Observe (with pieces of notation introduced in \sref{second-t-b}) that 
\begin{equation}\label{rel-f**}
\begin{array}{rll}
p \circ f_{**_x} &=& f_{*_x} \circ p \\
p_* \circ f_{**_x} &=& f_{*_x} \circ p_* \\
f_{**_x} \circ i^p_{0_x} &=& i^p_{0_y} \circ f_{*_x}.
\end{array}
\end{equation}
\end{rmk}

More generally, consider the natural left action $\rho^{(1,1)}$ of $\b^{(1,1)}(\PP(M))$ on $T^2M$ (cf.~\dref{derived-action} and \rref{derived-action-of-pair-groupoid})~:
\begin{equation}\label{(1,1)-action}
\begin{array}{cllll}
\rho^{(1,1)} & : & \b^{(1,1)}(\PP(M)) \times_{(\alpha, p^2)} T^2M & \to & T^2M \\
& & \displaystyle{\Bigl(j^1_xb, \X = \frac{dX_t}{dt}\Bigr|_{0}\Bigr) } & \mapsto & \displaystyle{ j^1_xb \cdot \X = \frac{d (b \cdot X_t)}{dt}\Bigr|_{0}. }
\end{array}
\end{equation}
In particular if $\X = Y_{*_x}X_x \in T^2M$ and $\beta(b(x)) = y$, then
\begin{equation}\label{coucou}
j^1_xb \cdot \X = (bY)_{*_y}(b(x) X_x).
\end{equation}

The relevance of the next lemma is essentially that it allows to carry the canonical involution $\kappa$ of $T^2M$ to a canonical involution on $\b^{(1,1)}_{sh}(\PP(M))$. Indeed, given a $(1,1)$-jet $\xi$, the following expression~:
$$\kappa(\xi) \cdot \X = \kappa \bigl( \xi \cdot \kappa(\X)\bigr).$$
define a map $T^2_xM \to T^2_yM$, which corresponds, thanks to the next lemma, to a new $(1,1)$-jet, defined to be $\kappa(\xi)$ (see~\dref{kappa-definition}). Recall that for a groupoid $G \rightrightarrows M$ and two point $x$ and $y$ in $M$, the notation $G_{x,y}$ stands for the subset $\alpha^{-1}(x) \cap \beta^{-1}(y)$ of $G$. 

\begin{lem}\label{11-jetsasmaps} Through the action $\rho^{(1,1)}$, the set of $(1,1)$-jets $(\b^{(1,1)}(\PP(M)))_{x,y}$ is in one-to-one correspondence with the set of maps $\ell : T^2_xM \to T^2_yM$ enjoying the following properties~:
\begin{enumerate}
\item[(a)] $\ell$ is a vector bundle morphism from $(T^2_xM, p)$ to $(T^2_yM, p)$ over a linear map $p(\ell) : T_xM \to T_yM$ and a vector bundle morphism from $(T^2_xM, p_*)$ to $(T^2_yM, p_*)$ over a linear map $p_*(\ell) : T_xM \to T_yM$~: \\
\vspace{.2cm}
\hspace{1.3cm}
\xymatrix{
T^2_xM \ar[r]^{\ell} \ar[d]_p & T^2_yM \ar[d]^p & & & T^2_xM \ar[r]^{\ell} \ar[d]_{p_*} & T^2_yM \ar[d]^{p_*} \\
T_xM \ar[r]^{p(\ell)} & T_yM & & & T_xM \ar[r]^{p_*(\ell)} & T_yM
}
\item[(b)]  $\ell$ coincides with $p(\ell)$ on each fiber of the vertical sub-bundle $T^p_{0_M}TM$~: \\
\vspace{.2cm}
\hspace{1.3cm}
\xymatrix{
T^p_{0_x}TM \ar[r]^{\ell} & T^p_{0_y}TM \\
T_xM \ar[u]^{i^p_{0_x}}\ar[r]^{p(\ell)} & T_yM \ar[u]_{i^p_{0_y}}
} 
\end{enumerate}
A $(1,1)$-jet $\xi$ corresponds to a map $\ell_\xi$ with $p(\ell_\xi) = p(\xi)$ and $p_*(\ell_\xi) = p_*(\xi)$. Genuine $2$-jets $j^2_xf \in \b^{(2)}(\PP(M))$ induce maps, also denoted by $f_{**_x}$, that, in addition, commute with $\kappa$, i.e.
\begin{equation}\label{f**-kappa}
f_{**_x} \circ \kappa = \kappa \circ f_{**_x}.
\end{equation}
\end{lem}

\begin{rmk}\label{rem-f**-kappa}
The relation (\ref{f**-kappa}) holds for any smooth map $f : M \to N$, where $\kappa$ stands either for the involution on $T^2M$ or on $T^2N$.
\end{rmk}

The next lemma is used in the proof of the previous one.

\begin{lem}\label{decomposition}
Let $U$ denote some open subset of $M$ and let $\{X^1, ..., X^n\}$ be a set of local vector fields in $\IX(U)$ forming a basis of each tangent space $T_{x}M$, $x \in U$. Given a path 
$$(-\varepsilon, \varepsilon) \mapsto T_{\gamma(t)}M : t \mapsto  X_t = \sum_{j=1}^n a_j(t) X^j(\gamma(t))$$ 
in $TM$. Its velocity vector $\displaystyle{\X = \frac{dX_t}{dt}|_{t=0}}$ admits the following expression~:
$$\displaystyle{ \X = \sideset{}{_{*}}\sum_{j=1}^n m_{a_j(0)*} \bigl(X^j_{*_x} Y_x\bigr) + \Bigl[i\Bigl(a_j(0) X^j_x \Bigr) +_* i^p_{0_M} \Bigl(\frac{da_j}{dt}\Bigl|_{t=0}X^j_x \Bigr) \Bigr] },$$
where $\sideset{}{_{*}}\sum$ indicates that the addition is $+_*$ and where $\displaystyle{Y_x = \frac{d \gamma(t)}{dt}|_{t=0}}$.
\end{lem}

\Pf It is essentially the same statement as \rref{Leibniz}. A proof is nevertheless included mainly to prepare for the proof of the corresponding statement for $T^3M$ (\lref{111-jetsasmaps}). Let us abbreviate $a_j(0)$ by $a_j$.
$$\begin{array}{lll}
\X & = & \displaystyle{\frac{d}{dt}\sum_{j=1}^n m\Bigl(a_j(t), X^j(\gamma(t))\Bigr) \Bigl|_{t=0} 
= \sideset{}{_{*}}\sum_{j=1}^n m_{*_{(a_j, X^j_x)}} \Bigl( \partial_ta_j(0), X^j_{*_x} Y_x\Bigr) }\\
& = & \displaystyle{ \sideset{}{_{*}}\sum_{j=1}^n m_{*_{(a_j, X^j_x)}} \Bigl( 0_{a_j}, X^j_{*_x}Y_x \Bigr) + m_{*_{(a_j, X^j_x)}} \Bigl( \frac{da_j(t)}{dt}\Bigl|_{t=0}, 0_{X^j_x}\Bigr) }\\
& = & \displaystyle{ \sideset{}{_{*}}\sum_{j=1}^n m_{a_j*} \bigl(X^j_{*_x} Y_x \bigr) + \frac{d}{dt} m \Bigl(a_j + \frac{da_j}{dt}\Bigl|_{t=0} t, X^j_x\Bigr)\Bigr|_{t=0}}\\
& = & \displaystyle{ \sideset{}{_{*}}\sum_{j=1}^n m_{a_j*} \bigl(X^j_{*_x} Y_x \bigr) + \Bigl[ \frac{d}{dt} m\Bigl(a_j, X^j_x \Bigr)\Bigr|_{t=0} +_* m \Bigl(t, \frac{da_j}{dt}\Bigl|_{t=0} X^j_x \Bigr)\Bigr|_{t=0}}\Bigl]\\
& = & \displaystyle{ \sideset{}{_{*}}\sum_{j=1}^n m_{a_j*} \bigl( X^j_{*_x} Y_x \bigr) + \Bigl[ i\Bigl(a_j X^j_x \Bigr) +_* i^p_{0_M} \Bigl(\frac{da_j}{dt}\Bigl|_{t=0} X^j_x \Bigr) \Bigr]}
\end{array}$$ 
\cqfd

\noindent
{\bf Proof of \lref{11-jetsasmaps}} It is quite obvious that a $(1,1)$-jet does enjoy the properties $(a)$ and $(b)$. A detailed proof of the converse is provided  mainly because it will lighten up the proof of \lref{111-jetsasmaps}. Consider a map $\ell : T^2_xM  \to T^2_yM$ that satisfies $(a)$ and $(b)$. Let $\{X^1, ..., X^n\}$ be a basis of $T_xM$ and let $\{Y^1 = p(\ell)(X^1), ..., Y^n = p(\ell)(X^n)\}$ be the corresponding basis of $T_yM$. For each $1 \leq i \leq n$, let $H^i$ be a $n$-dimensional subspace in $T_{X^i}TM$ complementary to $T^p_{X^i}TM$ and let $X^i : U \to TM$ (\rp $Y^i : V \to TM$) be a local section of $TM$ defined over a neighborhood $U$ of $x$ (\rp $V$ of $y$) in $M$, passing through $X^i$ (\rp $Y^i$) and tangent to $H^i$ (\rp $\ell(H^i)$). We may assume that the set $\{X^1(x'), ..., X^n(x')\}$ (\rp $\{Y^1(y'), ..., Y^n(y')\}$) is a basis of $T_{x'}M$ (\rp $T_{y'}M$) for all $x' \in U$ (\rp $y' \in V$) and that there exists a local diffeomorphism $g : U \to V$ such that $g_{*_x} = p_*(\ell)$. Now define a local bisection $b$ of $\b^{(1)}(\PP(M))$ over $U$ as follows~:

$$b(x') : T_{x'}M \to T_{g(x')}M : \sum_{i = 1}^n a_i X^i(x') \mapsto \sum_{i=1}^n a_i Y^i(g(x')).$$
It is now trivial to verify, by means of \lref{decomposition}, that the action of the $(1,1)$-jet $j^1_xb$ on $T^2M$ coincides with $\ell$. Indeed, it amounts to verifying that the image of the vectors $X^j_{*_x}X^k_x$, $i(X^j_x)$, $i^p_{0_M}(X^j_x)$ under the action of $j^1_xb$ and that of $\ell$ agree. This is implemented in the construction of $b$.
\cqfd

\begin{rmk}\label{hypoth-b} About the importance of hypothesis $(b)$. Consider the map
$$m_{-1} \circ m_{-1*} : T^2M \to T^2M.$$
It is a homomorphism for both vector bundle structures on $T^2M$ and it preserves the vertical sub-bundle $T^pTM$ but restricts to $\id$ on $T^pTM$ instead of $-\id$. Hence it is not induced by a $(1,1)$-jet and in particular does not yield a (necessarily canonical) affine connection (cf.~\pref{symm-jet<-->conn}). 
\end{rmk}

\begin{dfn} A homomorphism of $T^2M$ is a bijective morphism of vector bundles $\ell : (T^2_xM, p_o) \to (T^2_yM, p_o)$, for some $x, y \in M$ and  for both $p_o = p$ and $p_o = p_*$ over linear isomorphisms denoted by $p(\ell)$ and $p_*(\ell) : T_xM \to T_yM$ respectively. The set of homomorphisms of $T^2M$ is denoted by $\EL(T^2M)$ and the map $\b^{(1,1)}(\PP(M)) \to \EL(T^2M) : \xi \mapsto \ell_\xi$ by $\fL$. The set $\EL(T^2M)$, which is endowed with a Lie groupoid structure for which $\fL$ is a groupoid morphism, has the following distinguished Lie subgroupoids~:
\begin{enumerate}
\item[-] $\EL^{(1,1)}(T^2M) = \fL(\b^{(1,1)}(\PP(M)))$,
\item[-] $\EL^{(1,1)}_{sh}(T^2M) = \fL(\b^{(1,1)}_{sh}(\PP(M)))$,
\item[-] $\EL^{(2)}(T^2M) = \fL(\b^{(2)}(\PP(M)))$.
\end{enumerate}
\end{dfn}

Notice that he difference between $\EL(T^2M)$ and $\EL^{(1,1)}(T^2M)$ is that the action of an element in $\EL(T^2M)$ on the vertical subbundle $\im i^p_{0_M}$ is through any linear map, generally unrelated to $p(\ell)$ or $p_*(\ell)$.

\begin{dfn}\label{kappa-definition}
The natural involution $\kappa$ on $\b_{sh}^{(1,1)}(M)$ is defined through the expression~:
\begin{equation}\label{kappa}
\kappa(\xi) \cdot \X = \kappa \bigl( \xi \cdot \kappa(\X)\bigr).
\end{equation}
It consists basically in exchanging the order of the two derivatives involved in a $(1,1)$-jet. 
\end{dfn}

We have not been able to find a definition of a canonical involution on $(1,1)$-jets in the literature. The common definition 
is a (non-canonical) involution on $J^1(J^1(Y \to X))$, for a fibration $Y \to X$, depending on the choice of a connection on $X$.

\begin{rmk}\label{kappa-rem}
Notice that if $\xi$ lies in $\b^{(1,1)}(\PP(M)) - \b^{(1,1)}_{sh}(\PP(M))$, the right-hand side of (\ref{kappa}) does not define anymore a $(1,1)$-jet as the condition {\it (b)} in \pref{11-jetsasmaps} fails. Indeed, 
$$\kappa(\xi) \cdot i^p_{0_M}(U) = \kappa \bigl( \xi \cdot \kappa(i^p_{0_M}(U))\bigr) = \kappa \bigl( \xi \cdot i^p_{0_M}(U)\bigr) = \kappa \bigl( i^p_{0_M}(p(\xi) \cdot U)\bigr) = i^p_{0_M}(p(\xi) \cdot U),$$
while a $(1,1)$-jet acts on a vertical vector via its $p$-component, which, in the case of $\kappa(\xi)$, must be $p_*(\xi)$. Nevertheless $\kappa$ is defined on the entire space $\EL(T^2M)$ of homomorphisms of $T^2M$.
\end{rmk}

\begin{lem}\label{kappa-prop} The involution $\kappa$ is a groupoid automorphism whose fixed point set is  $\b^{(2)}(\PP(M))$ and such that $p\circ \kappa = p_*$ and $p_* \circ \kappa = p$. 
\end{lem}


\begin{rmk}\label{affine-str} A $(1,1)$-jet $\xi$ is also a morphism between the affine bundles $P : T^2_xM \to T_xM \oplus T_xM$ and $P : T^2_yM \to T_yM \oplus T_yM$ over the map $p(\xi) \oplus p_*(\xi) : T_xM \oplus T_x M \to T_yM \oplus T_yM$. Its {\bf pure second order part} in the direction of two vectors $X_x$ and $Y_x$ in $T_xM$ is the affine map 
$$\xi(X_x, Y_x) : T_{X_x}^{Y_x}TM \to T_{p(\xi)(X_x)}^{p_*(\xi)(Y_x)}TM,$$
where $T_{X_x}^{Y_x}TM$ denotes the intersection $p^{-1}(X_x) \cap p_*^{-1}(Y_x)$. The set of $(1,1)$-jets over a fixed map $\xi_1 \oplus \xi_2 : T_x M \oplus T_xM \to T_yM \oplus T_yM$ is an affine space modeled on the space of bilinear maps from $T_xM \times T_xM$ to $T_yM$. Indeed, let $\xi_0, \xi \in \b^{(1,1)}(\PP(M))$ be such that $p(\xi_0) = p(\xi)$ and $p_*(\xi_0) = p_*(\xi)$, then
\begin{equation}\label{xi-xi0}
\xi - \xi_0 : T_xM \times T_xM \to T_yM : (X_x, Y_x) \mapsto \pi \Bigl(\xi \cdot Y_{*_x} X_x, \xi_0 \cdot Y_{*_x} X_x \Bigr),
\end{equation}
where the map $\pi$ has been defined by (\ref{i}) in \aref{second-t-b} and where $Y$ is a local vector field extending $Y_x$, whose choice is irrelevant. Indeed, formula (\ref{scalar-mult}) implies that the right hand side of (\ref{xi-xi0}) is $C^\infty(M)$-bilinear. Notice that if both $\xi_0$ and $\xi$ are holonomic, or $\kappa$-invariant then $\xi - \xi_0$ is symmetric. Indeed, let $X, Y$ be local vector fields extending $X_x, Y_x$ and satisfying $[X, Y]_x = 0$. Then
$$\begin{array}{ccl}
\bigl(\xi - \xi_0\bigr) (X_x, Y_x) & = & \pi \Bigl(\xi \cdot Y_{*_x} X_x, \xi_0 \cdot Y_{*_x} X_x \Bigr) \\
& = & \pi \Bigl(\kappa \bigl(\xi \cdot Y_{*_x} X_x\bigr), \kappa \bigl(\xi_0 \cdot Y_{*_x} X_x\bigr) \Bigr) \;\;\mbox{implied by (\ref{i-kappa})} \\
& = & \pi \Bigl(\xi \cdot X_{*_x} Y_x, \xi_0 \cdot X_{*_x} Y_x \Bigr)  \;\;\mbox{implied by \pref{involution-prop}} \\
& = &  \bigl(\xi - \xi_0\bigr) (Y_x, X_x).
\end{array}$$ 
\end{rmk}

We finally prove a lemma about conjugation of $2$-jets that reveals useful when dealing with torsionless affine connections (cf.~\sref{Kobayashi}). Let $f : U \subset M \to M$ be a local diffeomorphism of $M$ such that for some $x \in U$, we have $f_{*_x} = I_x : T_xM \to T_xM$. Then the map
$$f_{**_{X_x}} - I : T_{X_x}TM \to T_{X_x}TM : \X \mapsto f_{**_{X_x}}(\X) - \X$$
\begin{enumerate}
\item[-] vanishes on $T^p_{X_x}TM$ since $({f_{**}}_{X_x} - I) \circ i^p_{X_x} = 0$,
\item[-] takes value in $T^p_{X_x}TM$ since $p_{*_{X_x}}\circ ({f_{**}}_{X_x} - I) = 0$ (cf.~(\ref{rel-f**})). 
\end{enumerate}
Whence there is a linear map $d^2f(X_x) : T_xM \to T_xM$, such that $({f_{**}}_{X_x} - I)$ coincides with the composition 

$$T_{X_x}TM \stackrel{p_{*_{X_x}}}{\longrightarrow} T_xM \stackrel{d^2f(X_x)}{\longrightarrow} T_xM \stackrel{i^p_{X_x}}{\longrightarrow} T_{X_x}TM.$$

\begin{lem}\label{sec-der-commutation}
Consider $2$-jets $j^2_xf, j^2_xg, j^2_yh \in \b^{(2)}(\PP(M))$ such that $p(j^2_xf) = I_x$ and $p(j^2_xg) = \iota \circ p(j^2_yh)$. Then
$$d^2 (g \circ f \circ h)(X_y) = g_{*_{x}} \circ d^2f({h_*}_y X_y) \circ h_{*_y} + d^2 (g \circ h)(X_y).$$ In particular, if $g_{*_x} = -I_x = h_{*_x}$, then 
$$d^2 (g \circ f \circ h)(X_x) = - d^2f(X_x) + d^2 (g \circ h)(X_x).$$ 
\end{lem}

\Pf The proof is just a short verification.
\begin{multline*}
\Bigl({(g \circ f \circ h)_{**}}_{X_x} - I\Bigr) \\
 \begin{array}{ccl}
& = & \Bigl. {g_{**}}_{{h_*}_x(X_x)} \circ {f_{**}}_{{h_*}_x(X_x)} \circ {h_{**}}_{X_x} - I \\
& = & {g_{**}}_{{h_*}_x(X_x)} \circ \Bigl({f_{**}}_{{h_*}_x(X_x)} - I \Bigr) \circ {h_{**}}_{X_x} + {g_{**}}_{{h_*}_x(X_x)} \circ {h_{**}}_{X_x} - I \\
& = & {g_{**}}_{{h_*}_x(X_x)} \circ \Bigl(i^p_{{h_*}_x(X_x)} \circ d^2f({h_*}_x(X_x)) \circ {p_*}_{{h_*}_x(X_x)} \Bigr) \circ {h_{**}}_{X_x} \\
& & \Bigl. + i^p_{X_x} \circ d^2 (g \circ h)(X_x) \circ {p_*}_{X_x} \\
& = & i^p_{X_x} \circ \Bigl( {g_*}_{h(x)} \circ d^2f({h_*}_x(X_x)) \circ {h_*}_x +  d^2 (g \circ h)(X_x)\Bigr) \circ {p_*}_{X_x}
\end{array}
\end{multline*}
\cqfd
 
\section{Third tangent bundle}\label{third-der}

Consider the third tangent bundle, denoted $T^3M$, and defined to be the tangent bundle of the total space of $T^2M$. 
Its element will be denoted by ``frak" letters like ${\mathfrak X}$. As described below, it is endowed with three canonical projections onto $T^2M$, giving $T^3M$ three distinct vector bundle structures over $T^2M$, three natural projections onto $TM$ and one onto $M$. It contains three vertical inclusions of $T^2M$ and six vertical inclusions of $T^2M \oplus T^2M$ and admits three involutions, each permuting two of the vector bundle structures and permuting the inclusions of $T^2M$ and $T^2M \oplus T^2M$. \\

\noindent
{\bf Vector bundle structures~:} \\

$\bullet$ $p : T^3M \to T^2M : \IX = \frac{d\ZE_t}{dt}|_{t = 0} \mapsto p(\IX) = \ZE_0$. The fiberwise addition is denoted by $+ : T^3M \times_{(p,p)}T^3M \to T^3M$ and the scalar multiplication by a real $a$ by $m_a : T^3M \to T^3M : \IX \mapsto m_a(\IX)= a \IX$. \\

$\bullet$ $p_* : T^3M \to T^2M : \IX = \frac{d\ZE_t}{dt}|_{t = 0} \mapsto p_*(\IX) =  \frac{dp(\ZE_t)}{dt}|_{t = 0}$. The $p_*$-fiberwise addition is the differential of the fiberwise addition $+$ on $T^2M$~:
$$+_* : T^3M\times_{(p_*, p_*)}T^3M \to T^3M : \Bigl(\frac{d\ZE_t}{dt}\Bigl|_{t=0}, \frac{d\ZE'_t}{dt} \Bigr|_{t=0}\Bigr) \mapsto \frac{d(\ZE_t + \ZE'_t)}{dt} \Bigr|_{t=0},$$ 
where we assume without loss of generality that $p(\ZE_t) =p(\ZE'_t)$ for all $t$'s. Similarly, the scalar multiplication by a real $a$ is the differential of the scalar multiplication $m_a$ on $T^2M$~:
$$m_{a*} : T^3M \to T^3M : \IX = \frac{d\ZE_t}{dt}\Bigl|_{t=0} \mapsto m_{a*} (\IX) = \frac{dm_a(\ZE_t)}{dt}\Bigl|_{t=0}.$$

$\bullet$ $p_{**} : T^3M \to T^2M : \IX = \frac{d\ZE_t}{dt}|_{t = 0} \mapsto p_{**}(\IX) =  \frac{dp_*(\ZE_t)}{dt}|_{t = 0}$. The $p_{**}$-fiberwise addition and scalar multiplication by a real $a$ are denoted respectively by $+_{**}$ and $m_{a**}$ and are the differential of the $p_*$-fiberwise addition $+_*$ and scalar multiplication $m_{a*}$ on $T^2M$. \\

Let ${\mathfrak X} = \frac{d\ZE_t}{dt}|_{t = 0}$ in $T^3M$, with $t \mapsto \ZE_t$ a path in $T^2M$. Set $X_t = p(\ZE_t)$, $Y_t = p_*(\ZE_t)$, $x_t = p(X_t) = p(Y_t)$. Then 

$$\begin{array}{lll}
\Bigl.p({\mathfrak X}) = \ZE_0 &p \circ p ({\mathfrak X}) = X_0 & p_* \circ p ({\mathfrak X}) = Y_0 \\
\Bigl.p_*({\mathfrak X}) = \frac{dX_t}{dt}|_{t=0} \stackrel{\rm not}{=} \Y \quad & p \circ p_* ({\mathfrak X}) = p(\Y) = X_0 \quad & p_* \circ p_* ({\mathfrak X}) = p_*(\Y) \stackrel{\rm not}{=} Z \\
\Bigl.p_{**}({\mathfrak X}) = \frac{dY_t}{dt}|_{t=0} \stackrel{\rm not}{=} \X & p \circ p_{**} ({\mathfrak X}) = p(\X) = Y_0 & p_* \circ p_{**} ({\mathfrak X}) = p_*(\X) \stackrel{\rm not}{=} Z \\
\end{array}$$

\begin{figure}[h!]  
\begin{center}
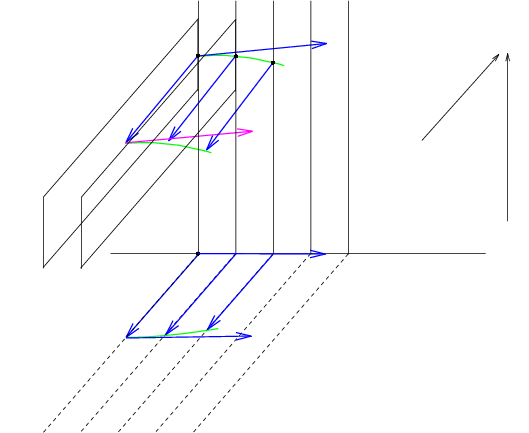 
\caption{A picture of $T^3M$}
\end{center}
\end{figure}
A comment about this picture. We think of $T^3M$ as the set of tangent vectors to the union over all $X_x \in TM$ of the tangent spaces $T_{X_x}TM$ that are represented on the picture above as $2$-planes sticking out ``horizontally''. The dotted lines indicates that that copy of $TM$ does not really lie in $T^3M$. It is added to the picture in order to better represent the $p_*$ projections of the vectors in $T^2M$. \\

The projections satisfy the following relations 
$$p \circ p = p \circ p_* \qquad p_* \circ p = p \circ p_{**} \qquad p_* \circ p_* = p_* \circ p_{**}.$$ 
For conciseness we denoted those three projections $p \circ p, p_* \circ p, p_* \circ p_* : T^3M \to TM$ by $p_1$, $p_2$, $p_3$ respectively. These commuting relations are summarized in the following diagram where the maps form the edges of a cube resting on a vertex (the altitude of a vertex depends on the number of derivatives)~:

\hspace{3cm}
\xymatrix{
& T^3M \ar[ld]_{p_*} \ar@{-->}[d]^{p} \ar[rd]^{p_{**}} & \\
T^2M \ar[d]_p \ar[rd]_<<<<{p_*} & T^2M \ar@{-->}[ld]^<<<{p} \ar@{-->}[rd] _<<<{p_*}& T^2M \ar[ld]^<<<{p_*} \ar[d]^{p}\\
TM \ar[rd]_p & TM \ar[d]^p & TM \ar[ld]^p \\
&M&
}
\vspace{.7cm}

Observe also that each projection $p, p_*, p_{**} : T^3M \to T^2M$ is linear for the vector bundle structures associated to the two other projections. More precisely, there are the following homomorphisms of vector bundles~:

\begin{equation}\label{peti}
\begin{array}{ll}
\Bigl.p : (T^3M, p_*) \to (T^2M, p) & p : (T^3M, p_{**}) \to (T^2M, p_*) \\
\Bigl.p_* : (T^3M, p) \to (T^2M, p) & p_* : (T^3M, p_{**}) \to (T^2M, p_*)\\
\Bigl.p_{**} : (T^3M, p) \to (T^2M, p) & p_{**} : (T^3M, p_*) \to (T^2M, p_*)
\end{array}
\end{equation}

\noindent
{\bf Horizontal inclusions and projections~:} Dual to the three projections $p, p_*, p_{**} : T^3 M \to T^2M$, there are three injections $i, i_*, i_{**} : T^2M \to T^3M$ whose images are the three different zero-sections for the vector bundle structures associated to $p, p_*, p_{**}$ respectively. Since $i : N \to TN$ denotes the canonical injection of a manifold in its tangent bundle as the zero section, then $i_*$ is the differential of $i : TM \to T^2M$ and $i_{**}$ is the second differential of $i : M \to TM$. Each inclusion $i_o$ is a vector bundle morphism for the same structures as for the corresponding projection $p_o$. The images of $i, i_*, i_{**}$ are denoted respectively by $0_{T^2M}, {0_*}_{T^2M}, {0_{**}}_{T^2M}$ and the image of a vector $\X$ by $0_\X, {0_*}_\X, {0_{**}}_\X$. Moreover, we have the following relations expressing that the cube is also commutative if one adds the arrows corresponding to $i$, $i_*$ and $i_{**}$~:
$$\begin{array}{lllll}
i \circ i = i_* \circ i && i \circ i_* = i_{**} \circ i && i_* \circ i_* = i_{**} \circ i_*
\end{array}$$
on $TM$ and 
$$\begin{array}{lllll}
p \circ i = id && p \circ i_* = i \circ p && p \circ i_{**} = i_* \circ p \\
p_* \circ i = i \circ p && p_* \circ i_* = id && p_* \circ i_{**} = i_* \circ p_* \\
p_{**} \circ i = i \circ p_* && p_{**} \circ i_* = i_* \circ p_* && p_{**} \circ i_{**} = id.
\end{array}$$
on $T^2M$. Let us introduce some more notation~:
\begin{enumerate}
\item[-] $i \circ i = i_* \circ i \stackrel{\rm not}{=} I_1$,
\item[-] $i \circ i_* = i_{**} \circ i \stackrel{\rm not}{=} I_2$,
\item[-] $i_* \circ i_* = i_{**} \circ i_* \stackrel{\rm not}{=} I_3$.
\item[-] $i \circ i \circ i \stackrel{\rm not}{=}{\bf i}$
\end{enumerate}
Notice that ${\bf i}$ coincides with any other inclusion of $M$ into $T^3M$ built from the various $i, i_*, i_{**}$'s. For convenience we will denote by $e$, $e_*$, $e_{**}$ the projections $i \circ p$, $i_* \circ p_*$, $i_{**} \circ p_{**}$ respectively which send a vector onto the zero vector in its $p$, $p_*$ or $p_{**}$-fiber respectively. The following relations are an easy consequence of (\ref{peti})~:
$$\begin{array}{lllll}
p \circ e = p && p \circ e_* = e \circ p && p \circ e_{**} = e_* \circ p \\
p_* \circ e = e \circ p_* && p_* \circ e_* =  p_* && p_* \circ e_{**} = e_* \circ p_* \\
p_{**} \circ e = e \circ p_{**} && p_{**} \circ e_* = e_* \circ p_{**} && p_{**} \circ e_{**} = p_{**}. \\
\end{array}$$
Furthermore, the projections $e$, $e_*$ and $e_{**}$ commute and the compositions of two such map is a new projection onto the image of some inclusion $I_j$ of $TM$. More precisely, set 
\begin{enumerate}
\item[-] $e \circ e_* \stackrel{\rm not}{=} E_1$,
\item[-] $e \circ e_{**} \stackrel{\rm not}{=} E_2$,
\item[-] $e_* \circ e_{**} \stackrel{\rm not}{=} E_3$.
\item[-] $e \circ e_* \circ e_{**} \stackrel{\rm not}{=} {\bf e}$.
\end{enumerate}
Then $E_j$ (\rp ${\bf e}$) is the projection of $T^3M$ onto the image of $I_j$ (\rp ${\bf i}$) and it satisfies~:
\begin{equation}\label{Eandp}
\begin{array}{lll}
p \circ E_1 = e \circ p & p_* \circ E_1 = e \circ p_* & p_{**} \circ E_1 = \ee \circ p_{**} \\
p \circ E_2 = e_* \circ p & p_* \circ E_2 = \ee \circ p_* & p_{**} \circ E_2 = e \circ p_{**}. \\
p \circ E_3 = \ee \circ p & p_* \circ E_3 = e_* \circ p_* & p_{**} \circ E_3 = e_* \circ p_{**}
\end{array}
\end{equation}
\noindent
{\bf Vertical inclusions of $T^2M$~:} There are three distinct vertical inclusions of $T^2M$ that parameterize the three transverse intersections 
$$\begin{array}{lll}
V_1 = p^{-1}(0_{TM}) \cap p_*^{-1}(0_{TM}) \\
V_2 = p_*^{-1}({0_*}_{TM}) \cap p_{**}^{-1}({0_*}_{TM}) \\
V_3 =  p^{-1}({0_*}_{TM}) \cap p_{**}^{-1}(0_{TM})
\end{array}$$
$$\begin{array}{lll}
\Biggl. i^p_{0_{TM}} : T^2M \stackrel{\sim}{\longrightarrow} V_1 = T^p_{0_{TM}}(T^2M) & : & \V \mapsto \displaystyle{\frac{d (t \V)}{dt}\Bigl|_{t=0}} \\
\Biggl. (i^p_{0_M})_{*} : T^2M \stackrel{\sim}{\longrightarrow} V_2 = T(T^p_{0_M}TM) & : & \displaystyle{\V = \frac{d V_t}{dt}\Bigl|_{t=0} \mapsto
\frac{d(i^p_{0_M}(V_t) )}{dt} \Bigl|_{t=0}} \\
\Biggl. i^{p_*}_{{0_*}_{TM}} : T^2M \stackrel{\sim}{\longrightarrow} V_3 = T^{p_*}_{{0_*}_{TM}}(T^2M) & : & \V \mapsto  \displaystyle{\frac{d (m_{t*} (\V))}{dt}\Bigl|_{t=0},} 
\end{array}$$
where indeed, the second one is the differential of the vertical inclusion $i^p_{0_M}$ of $TM$ into $T^2M$. These maps satisfy the following relations~:
$$\begin{array}{lll}
\Bigl.p \circ i^p_{0_{TM}} = e & p_* \circ i^p_{0_{TM}} = e & p_{**} \circ i^p_{0_{TM}} = i^p_{0_M} \circ p_* \\
\Bigl.p \circ (i^p_{0_M})_{*} = i^p_{0_M} \circ p & p_* \circ (i^p_{0_M})_{*} = e_* & p_{**} \circ (i^p_{0_M})_{*} = e_* \\
\Bigl.p \circ i^{p_*}_{{0_*}_{TM}} = e_* & p_* \circ i^{p_*}_{{0_*}_{TM}} = i^p_{0_M} \circ p & p_{**} \circ i^{p_*}_{{0_*}_{TM}} = i \circ p_*, 
\end{array}$$
and are vector bundle morphism in different ways~:
\begin{equation}\label{vert-incl}
\begin{array}{lll}
i^p_{0_{TM}} : (T^2M, p) \to (T^3M, \left\{\begin{array}{c} p \\ p_*\end{array}\right.) & i^p_{0_{TM}} : (T^2M, p_*) \to (T^3M, p_{**}) \\
(i^p_{0_{M}})_* : (T^2M, p) \to (T^3M, p) & (i^p_{0_{M}})_* : (T^2M, p_*) \to (T^3M, \left\{\begin{array}{c} p_* \\ p_{**} \end{array}\right.) \\
i^{p_*}_{{0_*}_{TM}} : (T^2M, p) \to (T^3M, p_*) & i^{p_*}_{{0_*}_{TM}} : (T^2M, p_*) \to (T^3M, \left\{\begin{array}{c} p \\ p_{**} \end{array}\right.).
\end{array}
\end{equation}
The presence of the bracket indicates that on the image of the inclusion at hand, the two linear structures coincide.
\begin{rmk} The other three intersections $p^{-1}(0_{TM}) \cap p_*^{-1}({0_*}_{TM})$, $p_*^{-1}({0_*}_{TM}) \cap p_{**}^{-1}(0_{TM})$ and $p^{-1}({0_*}_{TM}) \cap p_{**}^{-1}(0_{TM})$ could also be considered, but we do not need them here. A difference is that a vector $\IX$ in $p^{-1}(0_{TM}) \cap p_*^{-1}({0_*}_{TM})$ automatically belongs to $p^{-1}(0_{0_M})$.
\end{rmk}

\noindent
{\bf Vertical inclusion of $TM$~:} $T^3M$ supports also a vertical inclusion of $TM$~:
$$I : TM \stackrel{\sim}{\longrightarrow} T^3M : V_x \mapsto \frac{d}{dt} \Bigl(t \frac{d(sV_x)}{ds}\Bigl|_{s=0}\Bigr)\Bigl|_{t=0},$$
defined by pre-composing any vertical inclusions $i^p_{0_{TM}}$, $i^{p_*}_{{0_*}_{TM}}$ or $(i^p_{0_M})_{*}$ with the vertical inclusion $i^p_{0_M}$. It is a vector bundle morphism between $TM$ and all three vector bundle structures on $T^3M$. Moreover~:
$$
p \circ I = p_* \circ I = p_{**} \circ I = {\bf i}.
$$

\noindent
{\bf Vertical inclusions of $T^2M \oplus T^2M$~:} The various ``kernels" $p^{-1}(0_X)$, $p^{-1}({0_*}_X)$, $p_*^{-1}(0_X)$, $p_*^{-1}({0_*}_X)$, $p_{**}^{-1}(0_X)$, $p_{**}^{-1}({0_*}_X)$ admit the following parameterization by $T^2M \oplus T^2M$ (where the direct sum is either with respect to the projection $p$ or to the projection $p_*$ as indicated below)~: 
\begin{enumerate}
\item[1)] $\I_p = I_p^{TM} : T^2M \times_{(p, p)}T^2M \to p^{-1}(0_{TM})$,   
$$\I_p (\Y_{X}, \V_{X}) = i_{*_{X}}(\Y_{X}) + i^p_{0_{TM}}(\V_{X}).$$ 
\item[2)] $\I_p^* : T^2M \times_{(p, p_*)}T^2M \to p^{-1}({0_*}_{TM})$, 
$$\I_p^* (\X_{X}, \V^{X}) = i_{**_{X}}(\X_{X}) + i^{p_*}_{{0_*}_{TM}}(\V^{X}).$$
\item[3)] $\I_{p_*} = I_{p_*}^{TM} : T^2M \times_{(p,p)} T^2M \to p_*^{-1}(0_{TM})$,
$$\I_{p_*} (\ZE_{X}, \V_{X}) =  i(\ZE_X) +_* i^p_{0_{TM}}(\V_X).$$
\item[4)] $\I_{p_*}^* = (I_p)_* : T^2M \times_{(p_*,p_*)} T^2M \to p_*^{-1}({0_*}_{TM})$, 
$$\I_{p_*}^*(\X^{X}, \V^{X}) = i_{**} (\X^X) +_* (i^p_{0_M})_* (\V^X).$$
\item[5)] $\I_{p_{**}} : T^2M \times_{(p_*,p_*)} T^2M \to p_{**}^{-1}(0_{TM})$,
$$\I_{p_{**}} (\ZE^{X}, \V^{X}) = i(\ZE^X) +_{**} i^{p_*}_{{0_*}_{TM}}(\V^X).$$
\item[6)] $\I^*_{p_{**}} = (I_{p_*})_* : T^2M \times_{(p_*, p_*)} T^2M \to p_{**}^{-1}({0_*}_{TM})$,
$$\I^*_{p_{**}} : (\Y^X, \V^X) = i_* (\Y^X) +_{**} (i^p_{0_M})_* (\V^X).$$
\end{enumerate}
The rule to form these vertical inclusion can be described as follows~: if $j$ stands either for nothing or for $*$ or for $**$ and $k$ stands for either nothing or $*$, then
$$p_j^{-1}\Bigl(\im (i_k)\Bigr) = \im (i_{l(j,k)}) +_{j} \Bigl(p_j^{-1}(\im i_k) \cap p_{l(j,k)}^{-1}(\im i_{l(j,k)})\Bigr),$$ 
where $l(j, k)$ is so that $p_j$ realizes a morphism between $(T^3M, p_{l(j,k)})$ and $(T^2M, p_k)$. \\

Observe that for $p_i = p$ or $p_*$, $i=1,2$, the direct sum $T^2M \times_{(p_1, p_2)} T^2M$ is naturally a vector bundle of rank $4n$ over $TM$ for the projection but also a vector bundle of rank $3n$ over $TM \times_{(p,p)}TM$ for the projection $\hat{p}_1 \times \hat{p}_2$, where $\hat{p}_1$ denotes $p_*$ if $p_1 = p$ and $p$ otherwise. With respect to these bundle structures, each inclusion of a direct sum $T^2M \oplus T^2M$ is a vector bundle morphism in two fashions~: 
\begin{enumerate}
\item[-] $\I_p : (T^2M \times_{(p,p)} T^2M, \left\{\begin{array}{l} p=p \\ p_* \times p_*\end{array}\right.) \to (T^3M, \left\{\begin{array}{l} p \\ p_{**}\end{array}\right.)$
\item[-] $\I_p^* : (T^2M \times_{(p,p_*)} T^2M, \left\{\begin{array}{l} p=p_* \\ p_* \times p \end{array}\right.) \to (T^3M, \left\{\begin{array}{l} p \\ p_* \end{array}\right.)$
\item[-] $\I_{p_*} : (T^2M \times_{(p,p)} T^2M, \left\{\begin{array}{l} p=p \\ p_* \times p_* \end{array}\right.) \to (T^3M, \left\{\begin{array}{l} p_* \\ p_{**} \end{array}\right.)$
\item[-] $\I_{p_*}^*  : (T^2M \times_{(p_*,p_*)} T^2M, \left\{\begin{array}{l} p_*=p_* \\ p \times p \end{array}\right.) \to (T^3M, \left\{\begin{array}{l} p_* \\ p \end{array}\right.)$
\item[-] $\I_{p_{**}} : (T^2M \times_{(p_*,p_*)} T^2M, \left\{\begin{array}{l} p_*=p_* \\ p \times p \end{array}\right.) \to (T^3M, \left\{\begin{array}{l} p_{**} \\ p_* \end{array}\right.)$
\item[-] $\I^*_{p_{**}} : (T^2M \times_{(p_*, p_*)} T^2M, \left\{\begin{array}{l} p_*=p_* \\ p \times p \end{array}\right.) \to (T^3M, \left\{\begin{array}{l} p_{**} \\ p \end{array}\right.)$ 
\end{enumerate}

\noindent
{\bf Affine structures over $T^2M \oplus T^2M$~:} Any choice of two projections $p_1, p_2$ amongst $p, p_*, p_{**} : T^3M \to T^2M$, yields an affine fibration $p_1 \times p_2 : T^3M \to T^2M \oplus T^2M$. Altogether this provides $T^3M$ with three affine fibration structures over some fiber-product of $T^2M$ with itself~:
$$\begin{array}{lllll}
\p_1 \stackrel{\rm not}{=} p \times p_* & : & T^3M & \to & T^2M \times_{(p, p)} T^2M \\
\p_2 \stackrel{\rm not}{=} p_* \times p_{**} & : & T^3M & \to & T^2M \times_{(p_*, p_*)} T^2M \\
\p_3 \stackrel{\rm not}{=} p_{**} \times p & : & T^3M & \to & T^2M \times_{(p, p_*)} T^2M
\end{array}$$
whose respective fibers $\p_i^{-1}(\X_1, \X_2)$ admit two distinct affine structures (one for each factor of the projection $\p_i$) modeled on the fiber of either $p$ or $p_* : T^2M \to TM$. A fiber $(p_1 \times p_2)^{-1}(\X_1, \X_2)$, endowed with its affine structures induced by $p_i$ will be denoted by $((p_1 \times p_2)^{-1}(\X_1, \X_2), p_i)$. It is modeled on the vector space $p_1^{-1}(\X_1) \cap p_2^{-1}(p_2(e_o(\IX)))$, where $e_o$ coincides with $e$, $e_*$ or $e_{**}$ depending on whether $p_1$ is $p$, $p_*$ or $p_{**}$. Thus

\begin{enumerate}
\item[-] $\bigl(\p_1^{-1}(\ZE, \Y \bigr), p)$ is modeled on $p^{-1}(\ZE) \cap p_*^{-1}(0_X)$,
\item[-] $\bigl(\p_1^{-1}(\ZE, \Y  \bigr), p_*)$ is modeled on $p_*^{-1}(\Y)\cap p^{-1}(0_X)$,
\item[-] $\bigl(\p_2^{-1}(\Y, \X \bigr), p_*)$ is modeled on $ p_*^{-1}(\Y)\cap p_{**}^{-1}({0_*}_Z)$,
\item[-] $\bigl(\p_2^{-1}(\Y, \X \bigr), p_{**})$ is modeled on $p_*^{-1}(\X) \cap p_*^{-1}({0_*}_Z)$,
\item[-] $\bigl(\p_3^{-1}(\X, \ZE \bigr), p_{**})$ is modeled on $ p_{**}^{-1}(\X) \cap p^{-1}({0_*}_Y)$,
\item[-] $\bigl(\p_3^{-1}(\X, \ZE \bigr), p)$ is modeled on $p^{-1}(\ZE)\cap p_{**}^{-1}(0_Y)$. 
\end{enumerate}

Let us to describe  explicitly the three canonical inclusions of $T^2M$ parameterizing the various affine fibers passing through a given element ${\mathfrak X} \in T^3M$. They are denoted by the symbols $A^{\IX}_{(\p_i,p_j)}$, where $p_j$ is one of the factors of $\p_i$, indicating that we consider the $\p_i$-fiber of $\IX$ endowed with the affine structure given by the first projection $p_j$. Notice that each such inclusion is also obtained by translation of one of the vertical inclusions $i^p_{0_{TM}}, i^{p_*}_{{0_*}_{TM}}, (i^p_{0_M})_*$ of $T^2M$ along one of the inclusions $i, i_*, i_{**}$. More precisely,

$$\begin{array}{lllllllll}
A^{\IX}_{(\p_1, p)} & : & T_XTM & \to & T^3M & : & \V_X & \mapsto & \IX + \Bigl( e(\IX) +_* i^p_{0_X}(\V_X)\Bigr) \\
A^{\IX}_{(\p_1, p_*)} & : & T_XTM  & \to & T^3M & : & \V_X & \mapsto & \IX +_* \Bigl( e_*(\IX) + i^p_{0_X}(\V_X)\Bigr) \\
A^{\IX}_{(\p_2, p_*)} & : & T^ZTM  & \to & T^3M & : & \V^{Z} & \mapsto & \IX +_* \Bigl(e_* (\IX) +_{**} (i^p_{0_M})_* (\V^Z) \Bigr) \\
A^{\IX}_{(\p_2, p_{**})} & : & T^ZTM  & \to & T^3M & : & \V^{Z} & \mapsto & \IX +_{**} \Bigl(e_{**} (\IX) +_* (i^p_{0_M})_* (\V^Z)\Bigr) \\
A^{\IX}_{(\p_3, p_{**})} & : & T^YTM  & \to & T^3M & : & \V^{Y} & \mapsto & \IX +_{**} \Bigl( e_{**}(\IX) + i^{p_*}_{i_*(X)}(\V^{Y})\Bigr) \\
A^{\IX}_{(\p_3, p)} & : & T^YTM  & \to & T^3M & : & \V^Y & \mapsto & \IX + \Bigl( e(\IX) +_{**} i^{p_*}_{i_*(X)}(\V^Y)\Bigr) 
\end{array}$$

\noindent
\begin{rmks} \
\begin{enumerate} 
\item[-] The place of the parentheses above is important, only because shifting them makes generally appear sums of elements of $T^3M$ not belonging to a same fiber of any of the vector bundle structures on $T^3M$. Nevertheless when displacing parentheses yields a sensible expression, it is guaranteed to agree with the initial one. 
\item[-] Since $A^{\IX}_{(\p_1, p)} = A^{\IX}_{(\p_1, p_*)}$, $A^{\IX}_{(\p_2, p_*)} = A^{\IX}_{(\p_2, p_{**})}$, $A^{\IX}_{(\p_3, p_{**})} = A^{\IX}_{(\p_3, p)}$, we may remove the second subscript $p_2$ from the notation and talk about the three maps $A^{\IX}_{\p_i}$, $i = 1, 2, 3$.
\end{enumerate}
\end{rmks}

Dually, there are three projection maps $\Pi_i$ from $T^3M \times_{(\p_i, \p_i)} T^3M$ to $T^2M$ defined by 
\begin{equation}\label{proj123}
\Pi_i(\IX^1, \IX^2) = \U \iff \IX^1 =A^{\IX_2}_{\p_i} (\U)
\end{equation}
\ \\
\noindent
{\bf Affine structure modeled on $TM$~:} There is yet another structure of affine fibration on $T^3M$, obtained by considering all three projections $p$, $p_*$ and $p_{**}$. Denote by ${\mathcal P}(M)$ the image of the map $\p = p \times p_* \times p_{**} : T^3M \to T^2M \times T^2M \times T^2M$, that is the set of triples $(\X_1, \X_2, \X_3)$ of vectors in $T^2M$ such that 
$$p(\X_1) = p(\X_2) \qquad p_*(\X_1) = p(\X_3) \qquad p_*(\X_2) = p_*(\X_3).$$
The map 
$$\begin{array}{llcll}
\p & : & T^3M & \to & {\mathcal P}(M) \\
&&{\mathfrak X} & \mapsto & (p({\mathfrak X}), p_*({\mathfrak X}), p_{**}({\mathfrak X}))
\end{array}$$
in an affine fibration with typical fiber modeled on $T_xM$. For each element $\IX$ in $T^3M$, with $p \circ p \circ p (\IX) = x$, there is a parameterization of its $\p$-fiber by $T_xM$ that admits six different expressions~:

\begin{equation}\label{param-P-fiber}
A^{\IX}_{\p}(V_x) = \left\{
\begin{array}{l}
\IX + \Bigl( e(\IX) +_* \bigl(e_*(e(\IX)) +_{**} I(V_x)\bigr) \Bigr) \\
\IX + \Bigl( e(\IX) +_{**} \bigl(e_{**}(e(\IX)) +_{*} I(V_x)\bigr)\Bigr)\\
\IX +_* \Bigl( e_*(\IX) + \bigl(e(e_*(\IX)) +_{**} I(V_x)\bigr) \Bigr) \\
\IX +_* \Bigl( e_*(\IX) +_{**} \bigl( e_{**}(e_*(\IX)) + I(V_x) \bigr)\Bigr)\\
\IX +_{**} \Bigl( e_{**}(\IX) + \bigl(e(e_{**}(\IX)) +_* I(V_x)\bigr) \Bigr) \\
\IX +_{**} \Bigl( e_{**}(\IX) +_* \bigl( e_*(e_{**}(\IX)) + I(V_x) \bigr)\Bigr)\\
\end{array}\right.
\end{equation}
Let us explain the first equality of (\ref{param-P-fiber}). First of all, $\p(I(V_x)) = (\ii(x), \ii(x), \ii(x))$ and $p_{**} \circ e_* \circ e = \ee \circ p_{**} = \ii(x)$ (see (\ref{Eandp})) imply that the sum $+_{**}$ makes sense. Adding $I(V_x)$ does not change the $\p$-fiber, so 
$$p_*\Bigl(e_*(e(\IX)) +_{**} I(V_x)\Bigr) = p_* (e_* (e (\IX))) = p_* ( e (\IX)).$$
Hence the sum $+_*$ makes sense as well. Furthermore, 
$$p\Bigl(e(\IX) +_* \bigl(e_*(e(\IX)) +_{**} I(V_x) \bigr)\Bigr) = p(e(\IX)) + p(e_*(e(\IX))) = p(\IX) + e(p(\IX)) = p(\IX),$$
so that the third sum $+$ is well-defined. The other expressions are treated similarly.\\ 

In particular, there is a map~: 
\begin{equation}\label{Pi}
\Pi: T^3M \times_{(\p, \p)} T^3M \to TM : (\IX^1, \IX^2) \mapsto \Pi(\IX^1, \IX^2),
\end{equation}
such that $\Pi(\IX^1, \IX^2) = V_x$ if $\IX^1 = A^{\IX^2}_{\p}(V_X)$. It satisfies 
$$\Pi(\IX^1, \IX^2) = \Pi(\IX^1 +_o \IX, \IX^2 +_o \IX)$$ 
when $+_o$ denotes either $+$, $+_*$ or $+_{**}$ and $p_o(\IX) = p_o(\IX^i) \in T^3M$ for the corresponding projection $p_o$. For an element $\IX \in T^3M$ that admits either one of the following descriptions
$$\IX = \left\{ \begin{array}{ll}
e(\IX) +_* \bigl(e_*(e(\IX)) +_{**} I(V_x)\bigr) & e(\IX) +_* \bigl(e_*(e(\IX)) +_{**} I(V_x)\bigr) \\
e_*(\IX) + \bigl(e(e_*(\IX)) +_{**} I(V_x)\bigr) & e_*(\IX) +_{**} \bigl( e_{**}(e_*(\IX)) + I(V_x) \bigr) \\
e_{**}(\IX) + \bigl(e(e_{**}(\IX)) +_* I(V_x)\bigr) & e_{**}(\IX) +_* \bigl( e_*(e_{**}(\IX)) + I(V_x) \bigr), 
\end{array}\right.$$
we set 
\begin{equation}\label{pi-bis}
\Pi(\IX) = \Pi (\IX, e(\IX)) = V_x.
\end{equation}

\noindent{\bf Involutions.} On $T^3M$, there are three natural involutive automorphisms, each permuting two of the three vector bundle structures. The first one is the natural involution $\kappa^{TM}$ of the second tangent bundle $T^2N$ of the manifold $N=TM$. It is denoted by either $\kappa_1$ or $\kappa$. The second one is the differential $\kappa_*^{M}$ of the involution $\kappa^M$ of $T^2M$. It is denoted by $\kappa_2$ or $\kappa_*$. The third one is the conjugate of $\kappa$ by $\kappa_*$ and is denoted by either $\kappa_3$ or $\kappa'$. Thus $\kappa_3 = \kappa_* \circ \kappa \circ \kappa_*$. These three involutions correspond to the three involutive automorphisms of the cube obtained by reflexion relative to the planes that contain the two vertices $T^3M$ and $M$. They generate the group --- isomorphic to $S_3$ --- of ``level-preserving" automorphisms of the cube resting on its vertex $M$. More precisely, 

\begin{equation}\label{kappa-proj}
\begin{array}{lll}
p \circ \kappa = p_* & p_* \circ \kappa = p & p_{**} \circ \kappa = \kappa \circ p_{**}\\
p \circ \kappa_* = \kappa \circ p & p_* \circ \kappa_* = p_{**} & p_{**} \circ \kappa_* = p_* \\
p \circ \kappa' = \kappa \circ p_{**} & p_* \circ \kappa' = \kappa \circ p_*& p_{**} \circ \kappa' = \kappa \circ p.
\end{array}
\end{equation} 

The first line follows directly from the corresponding properties (\ref{prop-kappa}) of the involution $\kappa^{TM}$ and \rref{rem-f**-kappa}. The second line consists in differentiating the relations (\ref{prop-kappa}) for $\kappa^M$. The third line follows from the first two. More is true~: $\kappa$ is a vector bundle isomorphism between $(T^3M, p)$ and $(T^3M, p_*)$

\begin{enumerate}
\item[-] $\kappa : (T^3M, p) \stackrel{\sim}{\longrightarrow} (T^3M, p_*)$ over $\id_{T^2M}$, 
\item[-] $\kappa : (T^3M, p_{**}) \stackrel{\sim}{\longrightarrow} (T^3M, p_{**})$ over $\kappa$, 
\item[-] $\kappa_* : (T^3M, p_*) \stackrel{\sim}{\longrightarrow} (T^3M, p_{**})$ over $\id_{T^2M}$,  
\item[-] $\kappa_* : (T^3M, p) \stackrel{\sim}{\longrightarrow} (T^3M, p)$ over $\kappa$, 
\item[-] $\kappa' : (T^3M, p) \stackrel{\sim}{\longrightarrow} (T^3M, p_{**})$ over $\kappa$, 
\item[-] $\kappa' : (T^3M, p_*) \stackrel{\sim}{\longrightarrow} (T^3M, p_*)$ over $\kappa$, .
\end{enumerate}
Moreover $\kappa$ is an isomorphism of each $(T^3M, p_o)$ over $\kappa: T^2M \to T^2M$. Whence follows a series of equalities relating $\kappa$ with the various inclusions and projections. \\
$$\begin{array}{lll}
\kappa \circ i = i_* & \kappa \circ i_* = i & \kappa \circ i_{**} = i_{**} \circ \kappa \\
\kappa_* \circ i = i \circ \kappa & \kappa_* \circ i_* = i_{**} & \kappa_* \circ i_{**} = i_* \\
\kappa' \circ i = i_{**} \circ \kappa & \kappa' \circ i_* = i_* \circ \kappa & \kappa' \circ i_{**} = i \circ \kappa. 
\end{array}$$
Whence 
$$\begin{array}{lll}
\kappa \circ e = e_* \circ \kappa  & \kappa \circ e_* = e \circ \kappa & \kappa \circ e_{**} = e_{**} \circ \kappa \\
\kappa_* \circ e = e \circ \kappa_* & \kappa_* \circ e_* = e_{**} \circ \kappa_* & \kappa_* \circ e_{**} = e_*\circ \kappa_* \\
\kappa' \circ e = e_{**} \circ \kappa'  & \kappa' \circ e_* = e_* \circ \kappa' & \kappa' \circ e_{**} = e \circ \kappa'. 
\end{array}$$
The involutions also permutes the vertical inclusions of $T^2M$~: 
$$\begin{array}{lll}
\kappa \circ i^p_{0_{TM}} = i^p_{0_{TM}} & \kappa \circ (i^p_{0_M})_* = i^{p_*}_{{0_*}_{TM}} & \kappa \circ i^{p_*}_{{0_*}_{TM}} = (i^p_{0_M})_* \\
\kappa_* \circ i^p_{0_{TM}} = i^{p_*}_{{0_*}_{TM}} \circ \kappa & \kappa_* \circ (i^p_{0_M})_* =  (i^p_{0_M})_* & \kappa_* \circ i^{p_*}_{{0_*}_{TM}} = i^p_{0_{TM}} \circ \kappa \\
\kappa' \circ i^p_{0_{TM}} = (i^p_{0_M})_* \circ \kappa & \kappa' \circ (i^p_{0_M})_* = i^p_{0_{TM}} \circ \kappa & \kappa' \circ i^{p_*}_{{0_*}_{TM}} = i^{p_*}_{{0_*}_{TM}}. 
\end{array}$$
We may now easily deduce the action of the involutions on the vertical inclusions of $T^2M \oplus T^2M$~:
$$\begin{array}{llllll}
\kappa \circ \I_p =  \I_{p_*} & \kappa_* \circ \I_p = \I_p^* \circ (\id \times \kappa) & \kappa' \circ \I_p = \I^*_{p_{**}} \circ (\kappa \times \kappa) \\
\kappa \circ \I_p^* = \I^*_{p_*} \circ (\kappa \times \id) & \kappa_* \circ \I_p^* = \I_p \circ (\id \times \kappa) &  \kappa' \circ \I_p^* =  \I_{p_{**}} \circ (\kappa \times \id) \\
\kappa \circ \I_{p_*} = \I_p &  \kappa_* \circ \I_{p_*} = \I_{p_{**}} \circ (\kappa \times \kappa) & \kappa' \circ \I_{p_*} = \I_{p_*}^* \circ (\kappa \times \kappa) \\
\kappa \circ \I_{p_*}^* = \I_p^* \circ (\kappa \times \id) & \kappa_* \circ \I_{p_*}^* = \I^*_{p_{**}} & \kappa' \circ \I_{p_*}^* = \I_{p_*} \circ (\kappa \times \kappa) \\
\kappa \circ \I_{p_{**}} = \I^*_{p_{**}} & \kappa_* \circ \I_{p_{**}} = \I_{p_*} \circ (\kappa \times \kappa) & \kappa' \circ \I_{p_{**}} = \I_p^* \circ (\kappa \times \id) \\
\kappa \circ \I_{p_{**}}^*= \I_{p_{**}} & \kappa_* \circ \I_{p_{**}}^*= \I^*_{p_*} & \kappa' \circ \I_{p_{**}}^*= \I_p \circ (\kappa \times \kappa) 
\end{array}$$
Finally, the involutions fix the vertical inclusion of $TM$~:
$$\kappa_o \circ I = I,$$ 
for $\kappa_o = \kappa$, $\kappa_*$ or $\kappa'$. As a consequence, the relations (\ref{param-P-fiber}) imply

\begin{equation}\label{i-kappa-bis}
\Pi\Bigl(\IX_1, \IX_2\Bigr) = \Pi\Bigl(\kappa_o(\IX_1), \kappa_o(\IX_2)\Bigr).
\end{equation}

\section{Structure of $\b^{(1,1,1)}_{nh}(\PP(M))$}\label{b111}

As to $\b^{(1,1,1)}_{nh}(\PP(M))$, its elements are called $(1,1,1)$-jets and are of the type $\xi = j^1_x j^1_\centerdot b_\centerdot$, where $(b_{x'})_{x' \in U_x}$ is a smooth family of local bisections of $\b^{(1)}(\PP(M))$ parameterized by the elements $x'$ of a neighborhood $U_x$ of $x$ in $M$. We consider the local bisection $j^1_\centerdot b_\centerdot : x' \mapsto j^1_{x'} b_{x'}$ of $\b^{(1,1)}_{nh}(\PP(M))$ and its first order jet $j^1_xj^1_{\centerdot} b_{\centerdot}$ at $x$. 

\begin{figure}[h!]  
\begin{center}
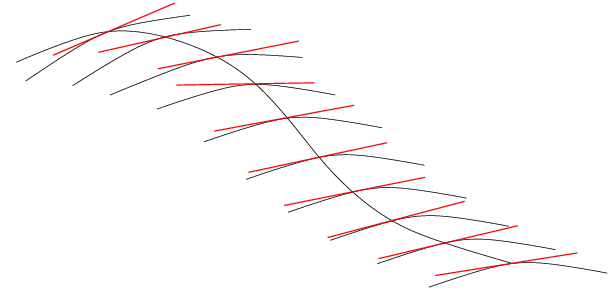 
\caption{The family $(b_{x'})_{x'\in U_x}$ and $j^1_{x'}b_{x'}$.}
\end{center}
\end{figure}

There are three natural projections from $\b^{(1,1,1)}_{nh}(\PP(M))$ to $\b^{(1,1)}_{nh}(\PP(M))$ denoted by $p$, $p_*$ and $p_{**}$ (cf.~\nref{bkl}). They admit the following description~:

\begin{enumerate}
\item[-] $p : \xi = j^1_x j^1_\centerdot b_\centerdot \mapsto (j^1_\centerdot b_\centerdot)(x) = j^1_x b_x$,
\item[-] $p_* : \xi \mapsto j^1_x(p \circ j^1_\centerdot b_\centerdot) = j^1_x(b_\centerdot(\centerdot))$,
\item[-] $p_{**} : \xi \mapsto j^1_x(p_* \circ j^1_\centerdot b_\centerdot) = j^1_x (j^1_\centerdot (p \circ b_\centerdot)) = j^1_x j^1_\centerdot (b^0_\centerdot)$.
\end{enumerate}
Remembering that the data of an element $\xi = j^1_x j^1_\centerdot b_\centerdot$ of $\b^{(1,1,1)}_{nh}(\PP(M))$ is equivalent to that of the plane 
$$D(\xi) = {(j^1_\centerdot b_\centerdot)}_{*_x} (T_xM) \subset T_{j^1_x b_{x}}\b^{(1,1)}_{nh}(\PP(M)),$$
the projections $p_*$ and $p_{**}$ are just the differential of the projections $p, p_* : \b^{(1,1)}_{nh}(\PP(M)) \to \b^{(1)}(\PP(M))$, i.e.~:
$$D(p_*(\xi)) = p_*(D(\xi)) \qquad D(p_{**}(\xi)) = (p_*)_*(D(\xi)).$$
Furthermore, the projections $p$, $p_*$ and $p_{**}$ satisfy the same relations as the corresponding projections from $T^3M$ to $T^2M$~: 
\begin{enumerate}
\item[-] $p \circ p = p \circ p_* : \xi = j^1_x j^1_\centerdot b_\centerdot \mapsto b_x(x)$,
\item[-] $p_* \circ p = p \circ p_{**} : \xi = j^1_x j^1_\centerdot b_\centerdot \mapsto j^1_x b^0_x$,
\item[-] $p_* \circ p_* = p_* \circ p_{**} : \xi = j^1_x (j^1_\centerdot b_\centerdot) \mapsto j^1_x (b^0_\centerdot(\centerdot))$.
\end{enumerate}
Altogether we obtain a cube resting on a vertex whose edges consist of groupoid morphisms~: \\

\hspace{0.5cm}
\xymatrix{
& \b^{(1,1,1)}_{nh}(\PP(M)) \ar[ld]_{p_*} \ar@{-->}[d]^{p} \ar[rd]^{p_{**}} & \\
\b^{(1,1)}_{nh}(\PP(M)) \ar[d]_p \ar[rd]_<<<<{p_*} & \b^{(1,1)}_{nh}(\PP(M)) \ar@{-->}[ld]^<<<{p} \ar@{-->}[rd] _<<<{p_*}& \b^{(1,1)}_{nh}(\PP(M)) \ar[ld]^<<<{p_*} \ar[d]^{p}\\
\b^{(1)}(\PP(M)) \ar[rd]_p & \b^{(1)}(\PP(M)) \ar[d]^p & \b^{(1)}(\PP(M)) \ar[ld]^p \\
&M&
}
\vspace{.7cm}

\begin{dfn}\label{b111h} Denote by $\b^{(1,1,1)}(\PP(M))$ the set of $(1,1,1)$-jets for which the three projections onto $\b^{(1,1)}_{nh}(\PP(M))$ coincide as well as the three projections onto $\b^{(1)}(\PP(M))$. In other terms 
\begin{multline}
\b^{(1,1,1)}(\PP(M)) = \Bigl\{\xi^{(1,1,1)} \in \b^{(1,1,1)}_{nh}(\PP(M)) \;\Bigr|\; \\
p(\xi) = p_*(\xi) = p_{**}(\xi) \in \b^{(1,1)}(\PP(M)) \Bigr\}.
\end{multline}

\end{dfn}
When $\xi = j^1_xj^1_\centerdot b_\centerdot$ lies in $\b^{(1,1,1)}(\PP(M))$, we may assume that $b_{x'}(x') = b_x(x')$ and that $b_{x'}$ is tangent to $\e$ at $x'$. 

\begin{figure}[h!]  
\begin{center}
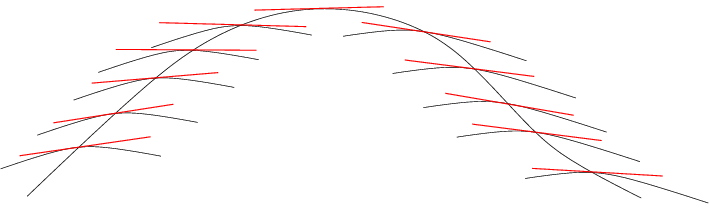 
\caption{Here $b_x (\centerdot) = (b_{\centerdot}(\centerdot))$.}
\end{center}
\end{figure}

Of course a genuine $3$-jet belongs to $\b^{(1,1,1)}(\PP(M))$ but $\b^{(1,1,1)}(\PP(M))$ contains $(1,1,1)$-jets that are not $3$-jets. \\

Remember the holonomic distribution on $\b^{(1,1)}_{nh}(\PP(M))$ introduced in \dref{def-e}~:

$$\e^{\b^{(1)}(\PP(M))}_\xi \stackrel{\rm not}{=} \e^{(1,1)}_\xi = p_{*_\xi}^{-1}(D(\xi)).$$

\begin{lem} An element $\xi$ in $\b^{(1,1,1)}_{nh}(\PP(M))$ belongs to $\b^{(1,1,1)}(\PP(M))$ if and only if $p(\xi) \in \b^{(1,1)}(\PP(M))$ and 
$$D(\xi) \subset \e^{(1,1)}_{p(\xi)} \cap T_{p(\xi)}\b^{(1,1)}(\PP(M)).$$
\end{lem}

\Pf The inclusion
$$D(\xi) \subset \e^{(1,1)}_{p(\xi)},$$
is equivalent to $p_*(D(\xi)) = D(p(\xi))$. Since $p_*(D(\xi)) = D(p_*(\xi))$, this amounts to $p_*(\xi) = p(\xi)$. Now, since $(1,1)$-jets $\xi$ belonging to $\b^{(1,1)}(\PP(M))$ are characterized by the relation $p(\xi) = p_*(\xi)$, the inclusion 
$$D(\xi) \subset T_{p(\xi)}\b^{(1,1)}(\PP(M))$$ 
is equivalent to $p_*(D(\xi)) = (p_{*})_*(D(\xi))$, that is $D(p_*(\xi)) = D(p_{**}(\xi))$ or $p_*(\xi) = p_{**}(\xi)$. 
\cqfd

\begin{lem}
A bisection $b$ of $\b^{(1,1)}(\PP(M))$ is everywhere tangent to $\e^{(1,1)}$ if and only if it is a $2$-jet extension~: $b = j^2f$.
\end{lem}

\Pf \pref{e} already implies that if a local bisection $b$ of $\b^{(1,1)}_{nh}(\PP(M))$ is tangent to $\e^{(1,1)}$ then it is a holonomic bisection of $\b^{(1)}(\b^{(1)}(\PP(M)))$, that is $b = j^1b'$ for $b' = p \circ b$. Now the local bisection $b'$ is necessarily tangent to $\e$. Indeed,
$$T_{b'(x')}b' = (p \circ b)_{*_{x'}}(T_{x'}M) = p_{*_{b(x')}}\Bigl(b_{*_{x'}}(T_{x'}M)\Bigr) \subset p_{*_{b(x')}}(\e^{(1,1)}_{b(x')}) = D(b(x')) \subset \e.$$
Whence the bisection $b'$ is holonomic as well~: $b' = j^1(p \circ b') = j^1b'^0$ and thus $b = j^2b'^0$.
\cqfd

Recall the natural action $\rho^{(1,1,1)}$ of the groupoid $\b^{(1,1,1)}_{nh}(\PP(M)) \rightrightarrows M$ on the fibration $T^3M \to M$~:
$$\b^{(1,1,1)} \times_{(\alpha, p^3)} T^3M : (\xi, {\mathfrak X}) \mapsto \xi \cdot {\mathfrak X},$$
where for $\xi = j^1_x b$, with $b \in \b_\ell(\b^{(1,1)}_{nh}(\PP(M)))$, and ${\mathfrak X} = \frac{d\X_t}{dt}\bigl|_{t=0}$, 
$$\xi \cdot {\mathfrak X} = \frac{d(b \cdot \X_t)}{dt}\Bigr|_{t=0}.$$
We will now characterize, as has been done for $\rho^{(1,1)}$ in a previous section, the partial maps $T^3M \to T^3M$ arising from the action of a $(1,1,1)$-jet.

\begin{dfn}\label{el-def} A homomorphism of $T^3M$ is a bijective map $\ell : T^3_xM \to T^3_yM$, $x$, $y \in M$ which is a vector bundle isomorphism $\ell : (T^3_xM, p_o) \to(T^3_yM, p_o)$ over a homomorphism of $T^2M$ denoted by $p_o(\ell)$ where $p_o$ stands for either $p$, $p_*$ or $p_{**}$. Moreover,
\begin{enumerate}
\item[-] $p \circ p(\ell) = p \circ p_*(\ell)$, 
\item[-] $p_* \circ p(\ell) = p \circ p_{**}(\ell)$,
\item[-] $p_* \circ p_*(\ell) = p_* \circ p_{**}(\ell)$.
\end{enumerate}
Let $\EL(T^3M)$ denote the set of homomorphisms of $T^3M$. It is naturally endowed with a groupoid structure.
\end{dfn}

\begin{rmk}\label{cons-el-def} As a consequence of this definition, an element $\ell$ of $\EL(T^3M)$ preserves the various zero sections in $T^3M$. More precisely,
$$\ell \circ i = i \circ p(\ell) \quad \ell \circ i_* = i_* \circ p_*(\ell) \quad \ell \circ i_{**} = i_{**} \circ p_{**}(\ell).$$
In particular, a homomorphism $\ell$ of $T^3M$ preserves the three vertical copies $p^{-1}(0_{TM}) \cap p_*^{-1}(0_{TM})$, $p_*^{-1}({0_*}_{TM}) \cap p_{**}^{-1}({0_*}_{TM})$ and $p^{-1}({0_*}_{TM}) \cap p_{**}^{-1}(0_{TM})$ of $T^2M$. Moreover, the fact that the vertical inclusion $i^p_{0_{TM}}$, $i^{p_*}_{{0_*}_{TM}}$, $(i^p_{0_M})_*$ are vector bundle morphisms as specified in (\ref{vert-incl}) implies that a homomorphism of $T^3M$ acts on their images through homomorphisms of $T^2M$. Similarly, a homomorphism of $T^3M$ preserves the vertical inclusion $I$ of $TM$ and acts linearly on its image.
\end{rmk}

\begin{lem}\label{111-jetsasmaps}
Via the action $\rho^{(1,1,1)}$, the groupoid $\b^{(1,1,1)}_{nh}(\PP(M))$ is canonically identified with the subset of $\EL(T^3M)$, denoted by $\EL^{(1,1,1)}(T^3M)$, of homomorphisms $\ell : T^3_xM \to T^3_yM$ assuming the following specific values on the images of the vertical inclusions $i^p_{0_{TM}}, (i^p_{0_M})_*, i^{p_*}_{{0_*}_{TM}}$~:
\begin{equation}\label{values-on-vert}
\begin{array}{rll}
\ell \circ i^p_{0_{TM}} & = & i^p_{0_{TM}}  \circ p(\ell) \\
\ell \circ (i^{p}_{0_M})_* & = & (i^{p}_{0_M})_* \circ p_*(\ell) \\
\ell \circ i^{p_*}_{{0_*}_{TM}} & = & i^{p_*}_{{0_*}_{TM}} \circ p(\ell). 
\end{array}
\end{equation}
Given a $(1,1,1)$-jet $\xi$, the associated linear map $\ell = \ell_\xi = \rho^{(1,1,1)}(\xi, \centerdot) : T^3_xM \to T^3_yM$ satisfies $p(\ell) = \EL(p(\xi))$, $p_*(\ell) = \EL(p_*(\xi))$ and $p_{**}(\ell) = \EL(p_{**}(\xi))$.
\end{lem}

\begin{rmk}
A homomorphism $\ell : T^3_xM \to T^3_yM$ acts on the vertical inclusion  $I : T_xM \to T^3_xM$ through $p(p(\ell)) = p(p_*(\ell))$, that is
$$\ell \circ I = I \circ p(p(\ell)).$$
This implies in particular that $\ell$ acts on the image of $A^{\IX}_{\p} : T_xM \to T^3M$, $\IX \in T^3M$ via $p \circ p(\ell)$ as well~:
\begin{equation}\label{action-sur-TM}
\ell \circ A^{\IX}_{\p} = A^{\IX}_{\p} \circ \; p(p(\ell)). 
\end{equation}
\end{rmk}

\begin{notas}
Let $\fL : \b^{(1,1,1)}_{nh}(\PP(M)) \to \EL(T^3M)$ denote the map $\xi \mapsto \ell_\xi$ and set 
\begin{enumerate}
\item[-] $\EL^{(1,1,1)}_h(T^3M) = \fL( \b^{(1,1,1)}(\PP(M)))$,
\item[-] $\EL^{(3)}(T^3M) = \fL( \b^{(3)}_h(M))$.
\end{enumerate}
\end{notas}
The following extension to $T^3M$ of \lref{decomposition} is useful in order to prove \lref{111-jetsasmaps}.
\begin{lem}\label{decomposition-ordre-3} Let $\{X^1, ..., X^n\}$ be a local basis of sections of $TM$ and write a given $\IX \in T^3M$ as follows~
$$\IX = \frac{d}{dt} \frac{d}{ds}\sum_{j=1}^na_j(t,s) X^j(\gamma(t,s)) \Bigl|_{s=0} \Bigr|_{t=0}.$$ 
Then $\IX$ admits the following expression as a linear combination of horizontal and vertical vectors~: 
$$\begin{array}{c}
{\displaystyle \sideset{}{_{**}}\sum_{j=1}^n } {\displaystyle \Biggl\{\Biggl[m_{a_j**} X^j_{**_{Y_x}}Y_{*_x}Z_x + }
\Bigl[m_{a_j**} \Bigl( i \bigl(X^j_{*_x}Y_x \bigr)\Bigr) +_{**}  m_{\partial_ta_j(0)**} \Bigl(i^{p_*}_{{0_*}_{TM}} \bigl( X^j_{*_x}Y_x \bigr)\Bigr) \Bigr] \Biggr]\\
 +_* \Biggl[\Bigl\{ m_{a_j**} \Bigl( i_{*_{X^j_x}} \bigl(X^j_{*_{x}} Z_x \bigr) \Bigr) + \Bigl[ m_{a_j**} \Bigl({\bf i} \bigl( X^j_x \bigr) \Bigr) +_{**} m_{\partial_t a_j(0)} \Bigl(i_{*_{0_x}} \bigl(i^p_{0_M} (X^j_x) \bigl) \Bigr) \Bigr] \Bigr\} \\
 \qquad +_{**} \; m_{\partial_s a_j(0)*} \Bigl( (i^p_{0_M})_* \bigl( X^j_{*_x} Z_x \bigr) \Bigr) + \Bigl[m_{\partial_s a_j(0)*}  i \Bigl( i^p_{0_M} (X^j_x ) \Bigr) +_* \partial^2_{ts}a_j(0) \bigl( I ( X^j_x)\bigr) \Bigr] \Bigr\} \Biggr] \Biggr\} .
\end{array}$$

\end{lem}

\Pf Set 
\begin{enumerate}
\item[-] $a_j(t) = a_j(t,0)$, 
\item[-] $a_j = a_j(0)$,
\item[-] $\partial_sa_j(t) = \frac{\partial a_j}{\partial s}(t,0)$, 
\item[-] $\partial^2_{ts}a_j(0) = \frac{\partial^2 a_j}{\partial t \partial s}(0,0)$, 
\item[-] $\gamma(t) = \gamma(t,0)$. 
\end{enumerate}
Using \lref{decomposition} we compute $\X_t = \frac{d}{ds}\sum_{j=1}^n a_j(t,s)X^j(\gamma(t,s))|_{s=0}$~:
$$\begin{array}{c}
\displaystyle{\X_t =  \sideset{}{_{*}}\sum_{j=1}^n \Bigl\{m_{a_j(t)*}X^j_{*_{\gamma(t)}} Y_t + \Bigl[ i\Bigl(a_j(t) X^j(\gamma(t)) \Bigr) +_* i^p_{0_M} \Bigl(\partial_sa_j(t) X^j(\gamma(t)) \Bigr) \Bigr] \Bigr\}},
\end{array}$$ 
where $Y_t = \frac{\partial \gamma(t,s)}{\partial s}|_{s=0}$. Now the vector $\IX$ is a sum of three types of vectors~:
\begin{multline}\label{IX}
\IX = \displaystyle{ \sideset{}{_{**}}\sum_{j=1}^n \Bigg\{ \frac{d}{dt} m_{a_j(t)*}X^j_{*_{\gamma(t)}} Y_t \Bigr|_{t=0} }\\
\displaystyle{ +_* \Biggl[ \frac{d}{dt} i\Bigl(a_j(t) X^j(\gamma(t)) \Bigr) \Bigr|_{t=0} } \\
\displaystyle{ +_{**} \frac{d}{dt} i^p_{0_M} \Bigl(\partial_sa_j(t) X^j(\gamma(t)) \Bigr) \Bigr|_{t=0} \Biggr] \Biggr\}. }
\end{multline}
The first term of (\ref{IX}) yields~:
$$\begin{array}{l}
\displaystyle{\frac{d}{dt} m_{a_j(t)*}X^j_{*_{\gamma(t)}} Y_t \Bigr|_{t=0} } \\
= \displaystyle{ \frac{d}{dt} m_{*} \Bigl( a_j(t), X^j_{*_{\gamma(t)}} Y_t \Bigr) \Bigr|_{t=0}}\\
= \displaystyle{ (m_*)_{*_{(a_j,X^j_{*_x}Y_x)}} \Bigl( \partial_ta_j(0), X^j_{**_{Y_x}}Y_{*_x}Z_x\Bigr)}\\
= \displaystyle{ (m_*)_{*} \Bigl( 0_{a_j}, X^j_{**_{Y_x}}Y_{*_x}Z_x\Bigr) + (m_*)_{*} \Bigl( \frac{d}{dt} a_j(0) + t\partial_ta_j(0)\Bigr|_{t=0}, 0_{X^j_{*_x}Y_{x}}\Bigr)}\\
= \displaystyle{ (m_*)_{*} \Bigl( 0_{a_j}, X^j_{**_{Y_x}}Y_{*_x}Z_x\Bigr) + \frac{d}{dt } m_* \Bigl( a_j + t \partial_ta_j(0), X^j_{*_x}Y_{x}\Bigr) \Bigr|_{t=0}}\\
= \displaystyle{ m_{a_j**} X^j_{**_{Y_x}}Y_{*_x}Z_x + \Bigl[ i \Bigl(m_{a_j*}X^j_{*_x}Y_x \Bigr) +_{**}  i^{p_*}_{{0_*}_{TM}} \Bigl(m_{\partial_ta_j(0)*} X^j_{*_x}Y_x \Bigr) \Bigr]} \\
= \displaystyle{ m_{a_j**} X^j_{**_{Y_x}}Y_{*_x}Z_x + \Bigl[m_{a_j**} \Bigl( i \bigl(X^j_{*_x}Y_x \bigr)\Bigr) +_{**}  m_{\partial_ta_j(0)**} \Bigl(i^{p_*}_{{0_*}_{TM}} \bigl( X^j_{*_x}Y_x \bigr)\Bigr) \Bigr], } \\
\end{array}$$
where $Y_{*_x}Z_x = \frac{dY_t}{dt}|_{t=0}$. In particular $Z_x = \frac{d\gamma(t)}{dt}|_{t=0}$. The second term yields~:
$$\begin{array}{l}
\displaystyle{ \frac{d}{dt} i\Bigl(a_j(t) X^j(\gamma(t)) \Bigr) \Bigr|_{t=0} }\\
= \displaystyle{ i_{*_{a_jX^j_x}} \Bigl[ m_{a_j*}X^j_{*_{x}} Z_x+ \Bigl( i(a_j X^j_x) +_* i^p_{0_M} \bigl(\partial_ta_j(0) X^j_x \bigr)\Bigr) \Bigr] }\\ 
= \displaystyle{ i_{*_{a_jX^j_x}} \Bigl(m_{a_j*}X^j_{*_{x}} Z_x \Bigr) + i_{*_{a_jX^j_x}} \Bigl( i(a_j X^j_x) +_* i^p_{0_M} \bigl(\partial_ta_j(0) X^j_x \bigr)\Bigr)}\\
= \displaystyle{ m_{a_j**} \Bigl( i_{*_{X^j_x}} \bigl(X^j_{*_{x}} Z_x \bigr) \Bigr) + i_{*_{a_jX^j_x}} \Bigl( m_{a_j*} \bigl(i(X^j_x)\bigr) +_* \partial_ta_j(0) i^p_{0_M} \bigl(X^j_x \bigr)\Bigr)}\\
= \displaystyle{ m_{a_j**} \Bigl( i_{*_{X^j_x}} \bigl(X^j_{*_{x}} Z_x \bigr) \Bigr) + \Bigl[ i_{*_{a_jX^j_x}} \Bigl(m_{a_j*} \bigl(i(X^j_x)\bigr) \Bigr) +_{**} i_{*_{0_x}} \Bigl(\partial_ta_j(0) i^p_{0_M} \bigl( X^j_x\bigr) \Bigr)\Bigr]}\\
= \displaystyle{ m_{a_j**} \Bigl( i_{*_{X^j_x}} \bigl(X^j_{*_{x}} Z_x \bigr) \Bigr) + \Bigl[ m_{a_j**} \Bigl({\bf i} \bigl( X^j_x \bigr) \Bigr) +_{**} \partial_t a_j(0) \Bigl(i_{*_{0_x}} \bigl(i^p_{0_M} (X^j_x) \bigl) \Bigr) \Bigr], }\\
\end{array}$$
and the third term yields~:
$$\begin{array}{l}
\displaystyle{ \frac{d}{dt} i^p_{0_M} \Bigl(\partial_sa_j(t) X^j(\gamma(t)) \Bigr) \Bigl|_{t=0}}\\
= \displaystyle{ \Bigl( (i^p_{0_M})_* \bigl( m_{\partial_s a_j(0)*}  X^j_{*_x} Z_x \bigr) \Bigr) + (i^p_{0_M})_* \Bigl[ i\Bigl(\partial_s a_j(0) X^j_x \Bigr) +_* i^p_{0_M} \Bigl(\partial^2_{ts}a_j(0) X^j_x\Bigr) \Bigr] }\\
= \displaystyle{ m_{\partial_s a_j(0)*} \Bigl( (i^p_{0_M})_* \bigl( X^j_{*_x} Z_x \bigr) \Bigr) + \Bigl[ (i^p_{0_M})_* \Bigl(m_{\partial_s a_j(0)*} i \bigl( X^j_x \bigr) \Bigr) +_* (i^p_{0_M})_* \Bigl(\partial^2_{ts}a_j(0) i^p_{0_M} \bigl( X^j_x \bigr) \Bigr) \Bigr] }\\
= \displaystyle{ m_{\partial_s a_j(0)*} \Bigl( (i^p_{0_M})_* \bigl( X^j_{*_x} Z_x \bigr) \Bigr) + \Bigl[m_{\partial_s a_j(0)*}  i \Bigl( i^p_{0_M} (X^j_x ) \Bigr) +_* \partial^2_{ts}a_j(0) \bigl( I ( X^j_x)\bigr) \Bigr]. }\\
\end{array}$$
\cqfd

\noindent
{\bf Proof of \lref{111-jetsasmaps}} The first part  of the proof consists in showing that the action of an element $\xi = j^1_xb = j^1_x j^1_{x'}b_{x'}$ of $\b^{(1,1,1)}_{nh}(\PP(M))$ on $T^3M$ is a homomorphism of $T^3M$. To handle the $p$-linearity, let $\IX_1$, $\IX_2$ belongs to some $p$-fiber of $T^3_xM$ and let $a$ be a real number. Then if $Z_i$ denotes $p_* \circ p_*(\IX_i)$, $i = 1, 2$, we have 
$$\begin{array}{cll}
j^1_xb \cdot \Bigl( a \IX_1 + \IX_2\Bigr) & = & \rho^{(1,1)}_* \Bigl( b_{*_x} (a Z_1 + Z_2), a\IX_1 + \IX_2\Bigr) \\
& = & a \rho^{(1,1)}_* \Bigl( b_{*_x} (Z_1), \IX_1 \Bigr) +  \rho^{(1,1)}_* \Bigl( b_{*_x} (Z_2), \IX_2\Bigr).
\end{array}$$
Supposing instead that $\IX_1$ and $\IX_2$ belong to the same $p_*$-fiber, implying in particular that $Z_1 = Z_2$, we have
$$\begin{array}{cll}
j^1_xb \cdot \Bigl( m_{a*} \IX_1 +_* \IX_2\Bigr) & = & \rho^{(1,1)}_* \Bigl( b_{*_x} (Z_1), m_{a*} \IX_1 +_* \IX_2\Bigr) \\
& = & \displaystyle{\frac{d}{dt} \rho^{(1,1)}\Bigl( b (\gamma(t)), a\X_{1t} + \X_{2t}\Bigr) \Bigr|_{t=0}, }
\end{array}$$
where $\frac{d}{dt}\gamma(t)|_{t=0} = Z_1$. The $p_*$-linearity follows from the linearity of $\rho^{(1,1)}$. Supposing now that $\IX_1$ and $\IX_2$ belong to the same $p_{**}$-fiber, we see that
$$\begin{array}{cll}
j^1_xb \cdot \Bigl( m_{a**} \IX_1 +_{**} \IX_2\Bigr) & = & \rho^{(1,1)}_* \Bigl( b_{*_x} (Z_1), m_{a**} \IX_1 +_{**} \IX_2\Bigr) \\
& = & \displaystyle{\frac{d}{dt} \rho^{(1,1)}\Bigl( j^1_{\gamma(t)}b_{\gamma(t)}, m_{a*} \X_{1t} +_* \X_{2t}\Bigr) \Bigr|_{t=0}, } \\
& = & \displaystyle{\frac{d}{dt} \frac{d}{ds} \rho^{(1)}\Bigl( b_{\gamma(t)}(\gamma(t,s)), aX_{1ts} + \X_{2ts}\Bigr) \Bigr|_{s=0}\Bigr|_{t=0}, }
\end{array}$$
The $p_{**}$-linearity follows thus from the linearity of $\rho^{(1)}$. \\

Now, let us show that the action of a $(1,1,1)$-jet takes the specific values (\ref{values-on-vert}) on the vertical copies of $T^2_xM$. Let $\xi = j^1_xb \in \b^{(1,1,1)}_{nh}(\PP(M))$ and $\V = \frac{dV_t}{dt}|_{t=0} \in T^2_xM$, then~:

$$\begin{array}{lll}
\xi\cdot i^p_{0_{TM}}(\V) & = & \displaystyle{j^1_x b \cdot \frac{d t\V}{dt} \Bigl|_{t=0}} = \displaystyle{\frac{d}{dt} b(x) \cdot t\V \Bigl|_{t=0}} = \displaystyle{\frac{d}{dt} t \bigl(b(x) \cdot \V \bigr) \Bigl|_{t=0}} \\
& = & \displaystyle{i^p_{0_{TM}} \Bigl(p(\xi) \cdot \V \Bigr)}, \\
\displaystyle{\xi\cdot (i^p_{0_{M}})_*(\V) } & = & \displaystyle{j^1_x b \cdot \frac{d \bigl(i^p_{0_M}(V_t)\bigr)}{dt} \Bigl|_{t=0}} = \displaystyle{\frac{d}{dt} b \cdot i^p_{0_M}(V_t) \Bigl|_{t=0}} =  \displaystyle{\frac{d}{dt} i^p_{0_M}(p(b) \cdot V_t) \Bigl|_{t=0}} \\
& = & \displaystyle{(i^p_{0_{M}})_* \Bigl(p_*(\xi) \cdot \V \Bigr)}, \\
\displaystyle{\xi\cdot i^{p_*}_{{0_*}_{TM}}(\V) } & = & \displaystyle{j^1_x b \cdot \frac{d\bigl(m_{t*} \V\bigr)}{dt} \Bigl|_{t=0}} = \displaystyle{\frac{d}{dt} \bigl(b(x) \cdot m_{t*} \V \bigr)\Bigl|_{t=0}} = \displaystyle{\frac{d}{dt} m_{t*} \bigl(b(x) \cdot \V \bigr) \Bigl|_{t=0}} \\
& = & \displaystyle{i^{p_*}_{{0_*}_{TM}} \Bigl(p(\xi) \cdot \V \Bigr)}.
\end{array}$$

The main part of the proof consists in showing the surjectivity of $\fL$ onto $\EL^{(1,1,1)}(T^3M)$. So let $\ell : T^3M \to T^3M$ in $\EL$ be a homomorphism that satisfies (\ref{values-on-vert}), we will show that it coincides with the action of a $(1,1,1)$-jet. In order to prove this, we need to construct a family of linear maps 
$$b_{x'}(x'') : T_{x''}M \to T_{b^0_{x'}(x'')}M,$$ 
with $(x', x'') \in \U = \cup_{x' \in U}\{x'\} \times U_{x'}$, where $U$ is a neighborhood of $x$ in $M$ and for $x' \in U$, $U_{x'}$ is a neighborhood of $x'$ in $M$, that ``integrates" $\ell$ in the sense that the action of the $(1,1,1)$-jet $\xi = j^1_x j^1_{x'} b_{x'}$ on $T^3M$ coincides with $\ell$. Let $\{X^1, ..., X^n\}$ be a local basis of vector fields defined on a neighborhood $U$ of $x$ in $M$. In view of the decomposition of any element $\IX$ in $T^3_xM$ described in \lref{decomposition-ordre-3}, the \lref{decomposition} and the \rref{cons-el-def}, it is sufficient to construct a family $b_{x'}(x'')$ for which the associated $(1,1,1)$-jet $\xi$ satisfies $\EL(p(\xi)) = p(\ell)$, $\EL(p_*(\xi)) = p_*(\ell)$ and $\EL(p_{**}(\xi)) = p_{**}(\ell)$ as well as 
\begin{equation}\label{horiz-action}
\xi \cdot X^j_{**_{X^k_x}}X^k_{*_x}(T_xM) = \ell \Bigl(X^j_{**_{X^k_x}}X^k_{*_x}(T_xM)\Bigr),
\end{equation}
for all $j, k = 1, ..., n$. To achieve these conditions, first integrate the homomorphism $p_{**}(\ell)$ of $T^2M$ into a local bisection $b : U \to \b^{(1)}(\PP(M))$ $: x' \mapsto b(x') = j^1_{x'}\varphi_{x'}$ such that $\EL(j^1_xb) = p_{**}(\ell)$ (cf.~\lref{11-jetsasmaps}).  Then consider a family of linear isomorphisms $b_{x}(x') : T_{x'}M \to T_{\varphi_x(x')}M$, $x' \in U$ such that $\EL(j^1_xb_x) = p(\ell)$ (once more \lref{11-jetsasmaps}). Now extend $b_x(x)$ into a family $b_{x'}(x') : T_{x'}M \mapsto T_{\varphi_{x'}(x')}M$, $x' \in U$ such that $\EL (j^1_x(b_{\centerdot}(\centerdot))) = p_*(\ell)$. So far, we have ensured that if $b_x(x')$ and $b_{x'}(x')$ is further extended to a family $b_{x'}(x'') : T_{x''}M \to T{\varphi_{x'}}(x'')M$, then the corresponding $(1,1,1)$-jet $\xi$ satisfies $\EL(p_o(\xi)) = p_o(\ell)$ for $p_o = p, p_*, p_{**}$. In order to guaranty the condition (\ref{horiz-action}), set 
$$H^j_{x'} = X^j_{*_{x'}}(T_{x'}M) \qquad \mbox{and} \qquad \h^j = \cup_{x' \in U}H^j_{x'}.$$ 
Then $\h^j$ is a submanifold of $T^2M$ which is completely determined by that data of the family $X^j_{*_{x'}}X^k_{x'}$, $x' \in U$, $k=1, ..., n$ of bases of the various horizontal spaces $H^j_{x'}$. Consider now the image of the differential of $X^j_{*}X^k$, that is 
$$X^j_{**_{X^k_x}}X^k_{*_x}(T_xM) \stackrel{\rm not}{=} T^j_k.$$ 
It is a horizontal $n$-plane in $T^3_xM$ tangent to $\h^j$ whose image under $\ell$ is a horizontal $n$-plane $\ell(T^j_k)$ in $T^3_yM$. For each pair $j,k = 1, ..., n$, let $E^j_k : V \to T^2M$ be a smooth section of $p^2 : T^2M \to M$ defined on a neighborhood $V$ of $y$ such that
\begin{enumerate}
\item[-] $(E^j_k)_{*_y} (T_yM) = \ell(T^j_k)$,
\item[-] $p \circ E^j_k(y') = b_{x'}(x') X^j(x')$, for $y' = \varphi_{x'}(x')$,
\item[-] $p_* \circ E^j_k(y') = j^1_{x'}(\varphi_\centerdot(\centerdot)) X^k_{x'}$, for $y' = \varphi_{x'}(x')$.
\end{enumerate}
The first condition implies in particular that $E^j_k(y) = p(\ell) (X^j_{*_x}X^k_x)$. Now for each $y'$ in $V$, let 
$$J^j_{y'} = {\rm span} \Bigl\{E^j_k(y') ; k = 1, ..., n \Bigr\}\subset T_{b_{x'}(x') X^j_{x'}}TM.$$ 
It is a horizontal space because the vectors $p_*(E^j_k(y'))$ are linearly independent and we may choose a smooth family $Y^j_{y'}$, $y' \in V$ of local vector fields such that 
$$(Y^j_{y'})_{*_{y'}}(T_{y'}M) = J^j_{y'}.$$ 
For $y'$ sufficiently close to $y$, and $y''$ sufficiently close to $y'$, the vector fields $Y^j_{y'}(y'')$ form a basis of $T_{y''}M$. This allows us to define $b_{x'}(x'')$, for $x' \in U$, $x'' \in U_{x'}$ via~:
$$b_{x'}(x'') X^j_{x''} = Y^j_{b^0_x(x')}(b^0_{x'}(x'')).$$

This proves surjectivity of $\EL$. Injectivity follows from the effectiveness of the action $\rho^{(1,1,1)}$ (cf.~\rref{derived-action-of-pair-groupoid}).
\cqfd

\begin{figure}[h!]
\begin{center}
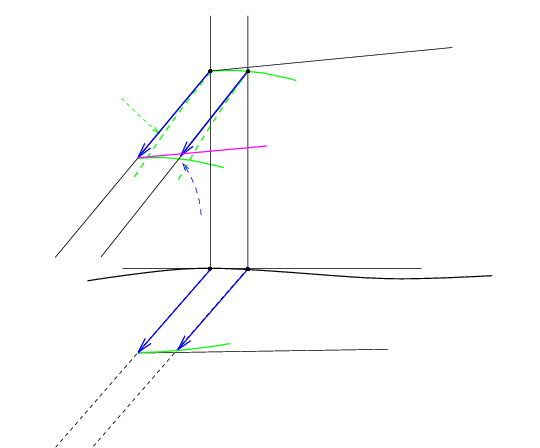 
\caption{A picture of the various actors of the proof}
\end{center}
\end{figure}

\begin{rmk}\label{vert-part} Denote by $\p(M)$ the subset of $\b^{(1,1)}_{nh}(\PP(M)) \times \b^{(1,1)}_{nh}(\PP(M)) \times \b^{(1,1)}_{nh}(\PP(M))$ consisting of the triples $(\xi_{1}, \xi_{2}, \xi_{3})$ for which $p(\xi_1) = p(\xi_2)$, $p_*(\xi_1) = p(\xi_3)$, $p_*(\xi_2) = p_*(\xi_3)$. It is the image of the projection $\p=p \times p_* \times p_{**} : \b^{(1,1,1)}_{nh}(\PP(M)) \to \b^{(1,1)}_{nh}(\PP(M)) \times \b^{(1,1)}_{nh}(\PP(M)) \times \b^{(1,1)}_{nh}(\PP(M))$. Now the map 
$$\p : \b^{(1,1,1)}_{nh}(\PP(M)) \to \p(M)$$ 
is an affine bundle whose fiber over any triple $(\xi_1, \xi_2, \xi_3)$ is modeled on the set of trilinear maps $T_xM \times T_xM \times T_xM \to T_yM$, where $x = \alpha(\xi_i)$ and $y = \beta(\xi_i)$. Indeed, let $\xi_0, \xi$ be two elements in $\p^{-1}(\xi_{1}, \xi_{2}, \xi_{3})$, then
$$\xi - \xi_0 : T_xM \times T_xM \times T_xM \to T_yM : (X_x, Y_x, Z_x) \mapsto \Pi (\xi \cdot \IX, \xi_0 \cdot \IX),$$
where $\IX \in T^3M$ satisfies $p \circ p(\IX) = X_x$, $p_* \circ p(\IX) = Y_x$ and $p_* \circ p_*(\IX) = Z_x$ defines a trilinear map independent on the choice of $\IX$. 
\end{rmk}

\noindent
{\bf Canonical Involutions~:} \lref{111-jetsasmaps} allows us to transport on $\b^{(1,1,1)}(\PP(M))$ the involutions $\kappa$, $\kappa_*$ and $\kappa'$ on $T^3M$. 
\begin{cor}
The expression
\begin{equation}\label{kappa-def}
\kappa_o(\xi) \cdot {\mathfrak X} = \kappa_o (\xi \cdot \kappa_o({\mathfrak X}))
\end{equation}
defines, for $\kappa_o = \kappa$, $\kappa_*$ or $\kappa'$ an involutive automorphism of the groupoid $\b^{(1,1,1)}(\PP(M))$ permuting two of the three fibrations. As is the case for $T^3M$, we have the relations~: 
$$\begin{array}{lll}
p \circ \kappa = p_* & p_* \circ \kappa = p & p_{**} \circ \kappa = \kappa \circ p_{**} \\
p \circ \kappa_* = \kappa \circ p & p_* \circ \kappa_* = p_{**} & p_{**} \circ \kappa_* = p_* \\
p \circ \kappa' = \kappa \circ p_{**} & p_* \circ \kappa' = \kappa \circ p_*& p_{**} \circ \kappa' = \kappa \circ p.\\
\end{array}$$
\end{cor}

\Pf It suffices to prove that the right-hand side of (\ref{kappa-def}) defines a map $T^3_xM \to T^3_yM$ satisfying the hypotheses of \lref{111-jetsasmaps}. This follows from the properties of the various involutions on $T^3M$. In addition,
$$\begin{array}{ccl}
\kappa_o(\xi_1 \cdot \xi_2) \cdot \IX & = & \kappa_o\bigl((\xi_1 \cdot \xi_2) \cdot \kappa_o(\IX)\bigr) \\
& = & \kappa_o\bigl(\xi_1 \cdot (\xi_2 \cdot \kappa_o(\IX)\bigr)\\
& = & \kappa_o(\xi_1) \cdot \kappa_o(\xi_2 \cdot \kappa_o(\IX)) \\
& = & \kappa_o(\xi_1) \cdot \kappa_o(\xi_2) \cdot \IX.
\end{array}$$
\cqfd

\begin{lem} The fixed point set of $\kappa$ (\rp $\kappa_*$) is $\b^{(2,1)}_h(M)$ (\rp $\b^{(1,2)}_h(M)$).
\end{lem}

\begin{rmk}\label{kappa-rem-bis} As for $(1,1,1)$-jets not in $\b^{(1,1,1)}(\PP(M))$, the expression (\ref{kappa-def}) does not in general define a $(1,1,1)$-jet (cf.~\rref{kappa-rem}). More precisely, we may define 
\begin{enumerate}
\item[-] $\kappa(\xi)$ when $p(\xi) = p_*(\xi)$, which implies that $p_{**}(\xi) \in \b_h^{(1,1)}$.
\item[-] $\kappa_*(\xi)$ when $p_*(\xi) = p_{**}(\xi)$, which implies that $p(\xi) \in \b_h^{(1,1)}$. 
\end{enumerate}
For an arbitrary element $\xi \in \b^{(1,1,1)}_{nh}(\PP(M))$ or even in $\EL(T^3M)$, one may define $\kappa(\xi)$ as the element in $\EL(T^3M)$ that satisfies~:
$$\kappa(\xi) \cdot \IX = \kappa(\xi \cdot \kappa(\IX)).$$
\end{rmk}

\begin{rmk}\label{kappa*-diff} Given a $(1,1,1)$-jet $\xi = j^1_xb$ for which $\kappa_*$ is defined, that is $p_*(\xi) = p_{**}(\xi)$, the corresponding tangent plane $D(\xi)$ is contained in the tangent space to the subgroupoid $\b^{(1,1)}(\PP(M))$ and $\kappa_*$ coincides with the differential of $\kappa^M$, that~is~:
$$D(\kappa_*(\xi)) = (\kappa^M)_*(D(\xi)).$$
Equivalently,
$$\kappa_*(j^1_xb) = j^1_x (\kappa \circ b).$$
\end{rmk}

\end{document}